\documentclass[11pt]{report}
\usepackage[letterpaper,hmargin=1in,vmargin=1.25in]{geometry}

\usepackage{amssymb, amsmath, amsthm, graphicx}
\usepackage[hypertex]{hyperref}

\def\~{\tilde }
\def\1{{\sf 1}}
\def\C{{\cal C}}
\def\L{{\cal L}}
\def\Q{{\sf Q}}
\def\I{{\sf I}}
\def\P{{\sf P}}

\def\D{{\sf D}}
\def\E{{\mathbb E}}
\def\R{{\mathbb R}}
\def\K{{\sf K}}

\def\MM{{\cal M}}

\newcommand{\charf}[1]{\mbox{\raise.48ex\hbox{$\chi$}$_{#1}$}}

\theoremstyle{plain}
\newtheorem{theorem}{Theorem}[section]

\newtheorem{lemma}[theorem]{Lemma}
\newtheorem{corollary}[theorem]{Corollary}

\newtheorem{remark}[theorem]{Remark}

\theoremstyle{definition}
\newtheorem{definition}[theorem]{Definition}
\newtheorem{example}[theorem]{Example}

\begin{document}

\title{Duality and evolving set bounds on mixing times}

\author {Ravi Montenegro \thanks{Department of Mathematical Sciences, University of Massachusetts Lowell, Lowell, MA 01854, ravi\_montenegro@uml.edu; partially supported by a VIGRE grant at the Georgia Institute of Technology.}}


\maketitle
\begin{abstract}
We sharpen the Evolving set methodology of Morris and Peres and extend it to study convergence in total variation, relative entropy, $L^2$ and other distances. Bounds in terms of a modified form of conductance are given which apply even for walks with no holding probability. These bounds are found to be strictly better than earlier Evolving set bounds, may be substantially better than conductance profile results derived via Spectral profile, drastically sharpen Blocking Conductance bounds if there are no bottlenecks at small sets, and give intuition into the workings of Canonical Path methods. 

\vspace{3ex}\noindent
This paper is intended solely to develop theoretical underpinnings, and as such we focus on two points : proving the sharpest most general results we can, and showing the Evolving Set methods to be better than previous isoperimetric methods. In order to learn about Evolving Sets we recommend the relevant chapter in our book with Tetali \cite{MT06.1}, and of course the original paper of Morris and Peres \cite{MP05.1}. To learn about some applications please see our paper on Cheeger Inequalities \cite{Mon07.1}, that on Canonical Path bounds for non-lazy walks \cite{Mon07.2}, our alternate interpretation of Morris' study of the Thorp shuffle \cite{MT06.1,Mor06.2}, Morris' paper on the Exclusion process \cite{Mor06.1}, and the paper of Diaconis and Fill on the duality method \cite{DF90.1}. 

\vspace{3ex}
\end{abstract}

\newpage
\tableofcontents

\newpage
\chapter{Introduction}

An isoperimetric bound on mixing time uses a geometric quantity, such as conductance, to bound the rate of convergence of a Markov chain. Such bounds have played a key role in proving mixing time results, beginning with Jerrum and Sinclair's \cite{JS88.1} proof that a random walk for approximating the permanent of a dense matrix converges in polynomial time. Their idea has been extended to apply to non-reversible non-lazy walks \cite{Mih89.1,Fill91.1}, to continuous state spaces \cite{LS93.1}, to walks with low conductance on small sets \cite{LS93.1}, and to walks with high conductance on small sets \cite{LK99.1}. 

Three recent papers have built on the Average Conductance idea of Lov\'asz and Kannan \cite{LK99.1}. Morris and Peres \cite{MP05.1} develop the Evolving Set methodology to show very strong results in terms of $L^2$ distance. Kannan, Lov\'asz and Montenegro \cite{KLM06.1} show similar results for total variation distance of a reversible, lazy walk through the method of Blocking Conductance. Finally, Goel, Montenegro, and Tetali \cite{GMT06.1} use the notion of Spectral Profile to extend an approach of Fill \cite{Fill91.1} and bound $L^2$ mixing of finite Markov chains. Each of these were shown by very different methods: by using a duality based approach, by considering the $n$-step average distribution, and by direct examination of the drop in variance, respectively.

The goal of this paper is to develop a general framework under which these isoperimetric results are unified as much as possible. This will be done by strengthening the Evolving Set methodology. Our improved argument leads to bounds on any convex notion of distance: including total variation, relative entropy, $L^2$, Hellinger, and Wasserstein distances. These are the first isoperimetric bounds on most of these distances, and even when past bounds are known these are the first which are sharp. For each of these distances we can also derive bounds in terms of an extension of the conductance method, known as modified conductance, which is consistent with past bounds when applied to lazy walks but which also applies in the setting of walks with no holding probability. 

How do our new Evolving Set results compare to previous isoperimetric bounds? We find that our new $L^2$ mixing bound is slightly better than earlier Evolving Set results, our conductance bounds on $L^2$ mixing may be substantially better than those derived from Spectral Profile bounds, and our mixing bounds are significantly sharper than those of Blocking Conductance except when the worst bottleneck is at a small set. Moreover, our results explain the curious existence of three total variation mixing bounds in the Blocking Conductance paper \cite{KLM06.1}. We find these are in fact total variation extensions of a bound on $L^2$ mixing, a bound on relative entropy mixing, and a direct bound on total variation mixing. An Evolving set approach to canonical paths also suggests that previous forms in terms of edge-expansion, vertex-expansion, or path lengths can all be considered  to be bounds of the form (edge-expansion)*(vertex-expansion).

This paper is focused on developing a rich theoretical framework, and comparing it to past methods. As such it is not so much as a text on applying Evolving Sets, as a text developing theory. Some examples are, however, included at the end in the Examples section. The interested reader can find additional theoretical developments in \cite{Mon07.1}, where we show a version of Cheeger's inequality which bounds (complex-valued) eigenvalues of non-reversible chains, a version to bound the smallest eigenvalue of a reversible chain, and we also sharpen Cheeger inequalities of Jerrum and Sinclair, Alon, and Stoyanov for bounding the spectral gap in terms of isoperimetric measures of edge and/or vertex expansion of a non-reversible walk. 
In \cite{Mon07.2} we develop a canonical path bound for non-reversible non-lazy walks, and use this to extend past results on mixing times of reversible walks on Cayley graphs into the general setting. Finally, together with Tetali \cite{MT06.1} we substantially improve on mixing time bounds of Morris for the Thorp shuffle \cite{Mor06.2}, by use of a conductance-profile bound based on ideas developed in this paper for walks with no holding probability.

The paper proceeds as follows. In Section \ref{sec:main} we introduce the notion of Evolving sets, and use this to show isoperimetric bounds on distances and mixing times. This is followed in Section \ref{sec:conductance} by conductance and modified conductance, an extension of conductance to non-lazy walks. These new results are compared to previous isoperimetric methods in Section \ref{sec:congestion}. 


\chapter{Set bounds on distance and Mixing Times} \label{sec:main}

In this section the main development of this paper is given: isoperimetric methods for bounding several notions of distance and mixing time. The arguments are based on the evolving set process of Morris and Peres \cite{MP05.1} which was also described in the context of duality by Diaconis and Fill \cite{DF90.1}.

A little notation is required. Let $\P$ be a finite irreducible Markov kernel on state space $V$ with stationary distribution $\pi$, that is, $\P$ is a $|V|\times|V|$ matrix with entries in $[0,1]$, row sums are one, $V$ is connected under $\P$ (i.e. $\forall x,y\in V\,\exists t:\,\P^t(x,y)>0$), and $\pi$ is a distribution on $V$ with $\pi\P=\pi$. The time-reversal $\P^*$ is given by $\P^*(x,y)=\frac{\pi(y)\P(y,x)}{\pi(x)}$ and is a Markov chain with stationary distribution $\pi$ as well. If $A,B\subset V$ the ergodic flow from $A$ to $B$ is given by $\Q(A,B)=\sum_{x\in A,y\in B} \pi(x)\P(x,y)$. Given initial distribution $\sigma$, the $n$-step discrete time distribution is given by $\sigma\P^n$, and if the walk is aperiodic then $\sigma\P^n\xrightarrow{n\to\infty}\pi$. 


\section{Duality and Evolving sets}

\begin{figure}[ht]
\begin{center}
\includegraphics[height=1.0in]{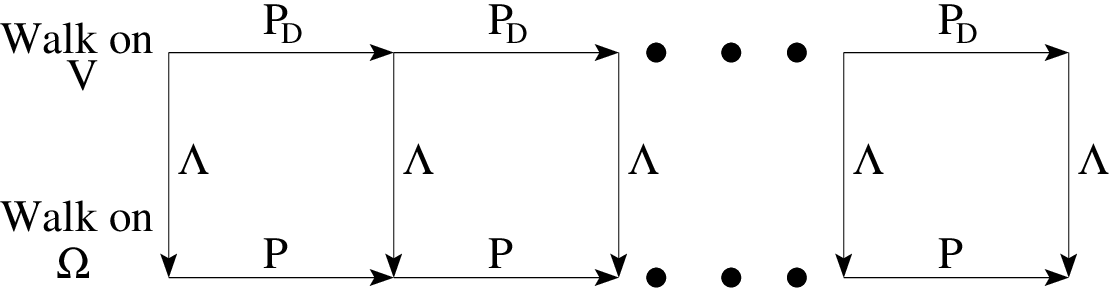}
\caption{The dual walk $\P_D=\hat\P$ projects onto the original chain $\P$.} \label{fig:duality}
\end{center}
\end{figure}

In order to relate a property of sets (conductance) to a property of the original walk (mixing time) we construct a walk on sets that is a dual to the original Markov chain. Given a Markov chain on $V$ with transition matrix $\P$, a {\em dual process} consists of a walk $\hat\P$ on some state space $V$ and a {\em link}, or transition matrix, $\Lambda$ from $V$ to $V$ such that
$$
\P\Lambda = \Lambda\hat\P\,.
$$
In particular, $\P^n\Lambda = \Lambda\hat\P^n$ and so the evolution of $\P^n$ and $\hat\P^n$ will be closely related. This relation is given visually by Figure \ref{fig:duality}.

The projection $\Lambda(S,y) = \frac{\pi(y)}{\pi(S)}\,\1_{S}(y)$ is a natural candidate to link a walk on sets to a walk on states. Diaconis and Fill \cite{DF90.1} have shown that for certain classes of Markov chains that the walk $\hat\K$ below is the unique dual process with link $\Lambda$, so this is the walk on sets that should be considered. We use notation of Morris and Peres \cite{MP05.1}.

\begin{definition}
Given set $A\subset V$ a step of the {\em evolving set process} is given by choosing $u\in[0,1]$ uniformly at random, and transitioning to the set $A_u = \{y\in V : \Q(A,y) \geq u\,\pi(y)\}$. The walk is denoted by $S_0$, $S_1$, $S_2$, $\ldots$, $S_n$, with transition kernel $\K^n(A,S)=Prob(S_n=S|S_0=A)$.
\end{definition}

\begin{definition}
The {\em Doob transform} of the Evolving set process is the Markov chain $\hat\K$ on sets with transition probabilities 
$$
\hat\K(S,S')=\frac{\pi(S')}{\pi(S)}\,\K(S,S')\,.
$$
The $n$-step transition probabilities are $\hat\K^n(S,S') = \frac{\pi(S')}{\pi(S)}\,\K^n(S,S')$.
\end{definition}

The Doob transform produces a Markov chain because of a Martingale property.

\begin{lemma} \label{lem:martingale} If $A\subset V$ then
$$
\int_0^1 \pi(A_u)\,du = \pi(A)\,.
$$
\end{lemma}

\begin{proof}
$$
\int_0^1 \pi(A_u)\,du = \sum_{y\in V} \pi(y)Prob(y\in A_u) = \sum_{y\in V} \pi(y)\frac{\Q(A,y)}{\pi(y)} = \pi(A)\,.
$$
\end{proof}

The walk $\hat\K$ is a dual process of $\P$.

\begin{lemma} \label{lem:duality}
If $S\subset V$, $y\in V$ and $\Lambda(S,y)=\frac{\pi(y)}{\pi(S)}\1_S(y)$ is the projection linkage, then
$$
\P\Lambda(S,y) = \Lambda\hat\K(S,y)\,.
$$
\end{lemma}

\begin{proof}
\begin{eqnarray*}
\P\Lambda(S,y) &=& \sum_{z\in S} \frac{\pi(z)}{\pi(S)}\P(z,y) = \frac{\Q(S,y)}{\pi(S)} \\
\Lambda\hat\K(S,y) &=& \sum_{S'\ni y} \hat\K(S,S')\,\frac{\pi(y)}{\pi(S')}
  = \frac{\pi(y)}{\pi(S)}\,\sum_{S'\ni y} \K(S,S') = \frac{\Q(S,y)}{\pi(S)}
\end{eqnarray*}
The final equality is because $\sum_{S'\ni y} \K(S,S')=Prob(y\in S')=\Q(S,y)/\pi(y)$.
\end{proof}

With duality it becomes easy to write the $n$ step density in terms of the walk $\hat\K$.

\begin{lemma} \label{lem:Doob}
Let $\hat\E_n$ denote expectation under $\hat\K^n$. If $x\in V$ and $S_0=\{x\}$ then
$$
\P^n(x,y) = \hat\E_n \pi_{S_n}(y)\,,
$$
where $\pi_S(y) = \frac{\1_S(y)\pi(y)}{\pi(S)}$ denotes the probability distribution induced on set $S$ by $\pi$.
\end{lemma}

\begin{proof}
$$
\P^n(x,y) = \P^n\Lambda(\{x\},y) = \Lambda\hat\K^n(\{x\},y) = \hat\E_n \pi_{S_n}(y)
$$
The final equality is because $\Lambda(S,y)=\pi_S(y)$.
\end{proof}


\section{Evolving set bounds on distances}

It is now a short hop from Lemma \ref{lem:Doob} to a bound on mixing times. First, however, note that if a distance $dist(\mu,\pi)$ is convex in $\mu$ (i.e. $dist(\alpha\mu_1+(1-\alpha)\mu_2,\pi)\leq \alpha dist(\mu_1,\pi)+(1-\alpha)dist(\mu_2,\pi)$), then for any distribution $\sigma$ 
$$
dist(\sigma\P^n,\pi)
    =  dist\left(\sum_{x\in V}\sigma(x)\P^n(x,\cdot),\pi\right) 
  \leq \sum_{x\in V} \sigma(x)dist\left(\P^n(x,\cdot),\pi\right)
  \leq \max_{x\in V} dist(\P^n(x,\cdot),\pi)\,.
$$
In this case distance is maximized when the initial distribution is a point mass, i.e. $\sigma(y)=\delta_{y=x}$ for some $x\in V$. Given the preceding lemmas it is easy to show an evolving set bound for all convex distances.

\begin{lemma} \label{lem:main_convex}
Consider a finite Markov chain with stationary distribution $\pi$. Any distance $dist(\mu,\pi)$ which is convex in $\mu$ satisfies
$$
dist(\P^n(x,\cdot),\pi) \leq \hat\E_n dist(\pi_{S_n},\pi)
$$
whenever $x\in V$ and $S_0=\{x\}$.
\end{lemma}

\begin{proof}
By Lemma \ref{lem:Doob} and convexity,
$$
dist(\P^n(x,\cdot),\pi) = dist(\hat\E_n \pi_{S_n},\pi) \leq \hat\E_n dist(\pi_{S_n},\pi)\,.
$$
\end{proof}

Many distances are used in studying mixing times. These include:
\begin{itemize}
\item {\em Separation distance: } $s(\mu,\pi)=\max_{x\in V} 1-\frac{\mu(x)}{\pi(x)}$
\item {\em Total variation distance: } $\|\mu-\pi\|_{TV} = \frac 12\sum_{x\in V} |\mu(x)-\pi(x)|$
\item {\em Relative Entropy: } $D(\mu\|\pi) = \sum_{x\in V} \mu(x)\log \frac{\mu(x)}{\pi(x)}$
\item {\em $L^2$ distance: } $\left\|\frac{\mu}{\pi}-1\right\|_{2,\pi} = \sqrt{\sum_{x\in V} \pi(x)\,\left|\frac{\mu(x)}{\pi(x)}-1\right|^2}$ \\
\item {\em Relative Pointwise distance ($L^{\infty}$): } $\left\|\frac{\mu}{\pi}-1\right\|_{\infty} = \max_{x\in V} \left|\frac{\mu(x)}{\pi(x)}-1\right|$
\item {\em Hellinger distance: } $H(\mu,\pi) = \sum_{y\in V} \left(\sqrt{\frac{\mu(y)}{\pi(y)}}-1\right)^2\,\pi(y)$
\item {\em Wasserstein distance $W_p(\mu,\pi)$: } Given metric $d:V\times V\rightarrow\R^+$, let
$$
W_p^p(\mu,\pi) 
  = \sup_{\substack{f,\,g:\,V\rightarrow \R,\\\forall y,z\in V:\,f(y)+g(z)\leq d(y,z)^p}} 
        \E_{\mu} f + \E_{\pi} g 
$$
\end{itemize}

Each of these distances can be bounded easily with Lemma \ref{lem:main_convex}.

\begin{theorem} \label{thm:main}
Given a finite, ergodic Markov chain, $x,y\in V$ and $S_0=\{x\}$, then
$$
\begin{array}{lcl}
\vspace{0.5ex}\displaystyle s(\P^n(x,\cdot),\pi) &\leq& \displaystyle Prob_{\hat\K^n}(S_n\neq V) \\
\vspace{0.5ex}\displaystyle \|\P^n(x,\cdot)-\pi\|_{TV} &\leq& \displaystyle \hat\E_n (1-\pi(S_n)) \\
\vspace{0.5ex}\displaystyle \D(\P^n(x,\cdot)\|\pi) &\leq& \displaystyle \hat\E_n \log \frac{1}{\pi(S_n)} \\
\vspace{0.5ex}\displaystyle \|\P^n(x,\cdot)-\pi\|_{2,\pi} &\leq& \displaystyle \hat\E_n \sqrt{\frac{1-\pi(S_n)}{\pi(S_n)}} \\
\vspace{0.5ex}\displaystyle \left|\frac{\P^n(x,y)}{\pi(y)}-1\right|
   &\leq& \displaystyle \max\left\{\frac{1-\pi(y)}{\pi(y)},\,1\right\}\,Prob_{\hat\K^n}(S_n\neq V) \\
\vspace{0.5ex}\displaystyle H(\P^n(x,\cdot),\pi) &\leq& \displaystyle 2\hat\E_n (1-\sqrt{\pi(S_n)}) \\
\displaystyle W_p(\P^n(x,\cdot),\pi) &\leq& \displaystyle \sqrt[p]{\hat\E_n W_p^p(\pi(S_n),\pi)}
\end{array}
$$
\end{theorem}

Most of these are immediate from the lemma and computation of $dist(\pi_S,\pi)$. For instance, in the total variation case $\|\pi_S-\pi\|_{TV}=1-\pi(S)$.

A few cases are worth mentioning further. The relative pointwise bound is because $dist(\mu,\pi)=\left|\frac{\mu(y)}{\pi(y)}-1\right|$ is convex, with
$$
dist(\pi_S,\pi)=\left| \frac{\pi_S(y)}{\pi(y)}-1\right| 
  = \1_S(y)\frac{1-\pi(S)}{\pi(S)}+\1_{S^c}(y) \leq \max\left\{\frac{1-\pi(y)}{\pi(y)},\,1\right\}\delta_{S\neq V}\,.
$$
The Hellinger distance is a special case of $dist(\mu,\pi) = \L_{\pi}\left(\frac{\mu}{\pi}\right)$ for a convex functional $\L_{\pi}:(\R^+)^V\to\R$. Wasserstein distance is a case of $\L_{\pi}(f) = \sup_{h\in H} \sum_{y\in V} h(y)\,f(y)\,\pi(y)$ for some class of functions $H$, by rewriting as
$$
W_p^p(\P^n(x,\cdot),\pi) = \sup_{\substack{f,\,g:\,V\rightarrow \R,\\\forall y,z\in V:\,f(y)+g(z)\leq d(y,z)^p}} 
        \sum_{y\in V} \left(f(y) + \E_{\pi} g\right) \left(\frac{\mu(y)}{\pi(y)}\right)\,\pi(y)\,.
$$

The Wasserstein distance $W_p^p$ is just the total variation distance when $d(y,z) = \delta_{y\neq z}$. It is easily checked that $W_p^p(\pi_S,\pi)=1-\pi(S)$ in this case, and so
$$
\|\P^n(x,\cdot)-\pi\|_{TV} = W_p^p(\P^n(x,\cdot),\pi) \leq \hat\E_n W_p^p(\pi_{S_n},\pi) = \hat\E_n (1-\pi(S_n))\,,
$$
which shows the Wasserstein bound generalizes the total variation bound.


\begin{remark}
When the initial distribution $\sigma$ is not a point mass then $S_0$ should be chosen from a distribution. Set $\sigma_0=\sigma$. Inductively define $A_i=\{x:\,\sigma_i(x)>0\}$, let $Prob(S_0=A_i)=\pi(A_i)\,\min_{x\in A_i} \frac{\sigma_i(x)}{\pi(x)}$, and $\sigma_{i+1}(x)=\sigma_i(x)-\1_{A_i}(x)\frac{\pi(x)}{\pi(A_i)}\,Prob(S_0=A_i)$. Then 
Lemma \ref{lem:Doob} generalizes to
$$
\frac{\sigma\P^n(y)}{\pi(y)}
   = \sum_{S\subset V} Prob(S_0=S)\,\frac{Prob(y\in S_n | S_0=S)}{\pi(S)} 
   = \hat\E_n \frac{\1_{S_n}(y)}{\pi(S_n)}\,.
$$
The results in Theorem \ref{thm:main} generalize to this case as well, whereas those in the next section will replace $\pi_*$ with $f^{-1}(\E f(\pi(S_0)))$.
\end{remark}


\section{Mixing times} \label{sec:mixing}

Throughout this section assume that the distance to be studied is of the form
$$
dist(\P^n(x,\cdot),\pi) \leq \hat\E_n f(\pi(S_n))
$$
for a decreasing function $f:\;[0,1]\to\R_+$. For instance, the total variation, $L^p$ and relative entropy bounds in Theorem \ref{thm:main} are all of this form. Let $\tau(\epsilon)$ denote the mixing time in this distance, that is, the minimum number of steps to guarantee that this distance is at most $\epsilon$.

Mixing time will be bounded using the $f$-congestion.
\begin{definition}
Given a finite Markov chain, and function $f:[0,1]\to\R_+$ non-zero except possibly at $0$ and $1$, then the $f$-congestion $\C_f$ and $f$-congestion profile $\C_f(r)$ are given by
$$
\forall A\subset V:\,\C_f(A) = \frac{\int_0^1 f(\pi(A_u))\,du}{f(\pi(A))}\,,
\quad
\forall r>0:\,\C_f(r)=\max_{\substack{\pi(A)\leq r,\\ A\neq \emptyset,\,V}} \C_f(A)\,,
\quad
\C_f = \C_f(1)\,.
$$
\end{definition}

The starting point for our calculations will be the following discrete analog of differentiation.

\begin{lemma} \label{lem:discreteCongestion}
\begin{eqnarray*}
\hat\E_{n+1} f(\pi(S_{n+1})) - \hat\E_n f(\pi(S_n)) &=& - \hat\E_n f(\pi(S_n))\,(1-\C_{af(a)}(S_n)) \\
  &\leq& -(1-\C_{af(a)})\,\hat\E_n f(\pi(S_n))
\end{eqnarray*}
\end{lemma}

\begin{proof} The inequality is because $\forall S\subset V:\,1-\C_{af(a)} \leq 1-\C_{af(a)}(S)$, by definition of $\C_{af(a)}$. For the equality, 
\begin{eqnarray*}
\lefteqn{\hat\E_{n+1} f(\pi(S_{n+1})) 
  = \hat\E_n \sum_S \hat\K(S_n,S) f(\pi(S)) }\\
  &=& \hat\E_n f(\pi(S_n))\frac{\sum_S \K(S_n,S) \pi(S)f(\pi(S))}{\pi(S_n)f(\pi(S_n))}
  = \hat\E_n f(\pi(S_n))\C_{af(a)}(S_n)
\end{eqnarray*}
\end{proof}

A basic mixing time bound follows easily:

\begin{corollary} \label{cor:main}
In discrete time
$$
\tau(\epsilon) \leq \left\lceil \frac{1}{1-\C_{af(a)}}\,\log \frac{f(\pi_*)}{\epsilon} \right\rceil
$$
\end{corollary}

\begin{proof}
By Lemma \ref{lem:discreteCongestion} $\hat\E_{n+1} f(\pi(S_{n+1})) \leq \C_{af(a)}\,\hat\E_n f(\pi(S_n))$.
Applying induction to this yields the relation $\hat\E_n f(\pi(S_n)) \leq \C_{af(a)}^n\,f(\pi(S_0))$. Solving for when this drops to $\epsilon$ and using the approximation $\log\C_{af(a)}\leq-(1-\C_{af(a)})$, gives the corollary.
\end{proof}

This can be generalized to take into consideration set sizes. A stronger bound holds under a fairly weak convexity condition, with about a factor of two lost in the general case.

\begin{theorem} \label{thm:profile-mixing}
If $x\left(1-\C_{af(a)}(f^{-1}(x))\right)$ is convex then
$$
\tau(\epsilon) \leq 
\left\lceil 
       \int_{\pi_*}^{f^{-1}(\epsilon)} \frac{-f'(x)\,dx}{f(x)(1-\C_{af(a)}(x))} 
\right\rceil\,,
$$
while in general
$$
\tau(\epsilon) \leq
\left\lceil 
       \int_{f^{-1}(f(\pi_*)/2)}^{f^{-1}(\epsilon/2)} \frac{-2f'(x)\,dx}{f(x)(1-\C_{af(a)}(x))} 
\right\rceil\,.
$$
\end{theorem}

\begin{proof} First consider the convex case.

By Lemma \ref{lem:discreteCongestion} and Jensen's inequality for the convex function $x\left(1-\C_{af(a)}(f^{-1}(x))\right)$,
\begin{eqnarray}
\hat\E_{n+1} f(\pi(S_{n+1})) &-& \hat\E_n f(\pi(S_n))
    = -\hat\E_n f(\pi(S_n))\,(1-\C_{af(a)}(S_n))\nonumber \\
   &\leq& -\hat\E_n f(\pi(S_n))\,\left[1-\C_{af(a)}\left(f^{-1}\circ f(\pi(S_n))\right)\right] \nonumber \\
   &\leq& -\left[\hat\E_n f(\pi(S_n))\right]\,
           \left[1-\C_{af(a)}\left(f^{-1}(\hat\E_n f(\pi(S_n)))\right)\right] \label{eqn:convexity} \,.
\end{eqnarray}
Since $I(n)=\hat\E_n f(\pi(S_n))$ and $1-\C_{af(a)}(f^{-1}(x))$ are non-increasing, the piecewise linear extension of $I(n)$ to $t\in\R_+$ satisfies
\begin{equation} \label{eqn:derivative}
I'(t) \leq -I(t)\,\left[1-\C_{af(a)}(f^{-1}(I(t)))\right]
\end{equation}
At integer $t$ the derivative can be taken from either right or left.

Then,
$$
\int_{I(0)}^{I(t)} \frac{dI}{I\left(1-\C_{af(a)}(f^{-1}(I))\right)} \leq -\int_0^t dt = -t\,.
$$
A change of variables to $v=f^{-1}(I)$ implies that
$$
\int_{f^{-1}(I(0))}^{f^{-1}(I(t))} \frac{f'(v)\,dv}{f(v)(1-\C_{af(a)}(v))} \leq -t\,.
$$
By continuity of $I(t)$ there exists $T$ such that $I(T)=\epsilon$. The theorem follows from $f^{-1}(I(0))=f^{-1}(f(\pi_*))=\pi_*$ and $f^{-1}(I(T))=f^{-1}(\epsilon)$.

For the general case, use Lemma \ref{lem:nonconvex} instead of convexity at \eqref{eqn:convexity}.
\end{proof}

\begin{lemma} \label{lem:nonconvex}
If $Z\geq 0$ is a nonnegative random variable and $g$ is a nonnegative increasing function, then
$$
E \left(Z\,g(Z)\right) \geq \frac{E Z}{2}\,g(EZ/2)\,.
$$
\end{lemma}

\begin{proof} See \cite{MP05.1}.
Let $A$ be the event $\{Z\geq EZ/2\}$. Then $E(Z\,\1_{A^c})\leq E Z/2$, so $E(Z\1_A)\geq EZ/2$. Therefore,
$$
E\left(Z\,g(2Z)\right) \geq E \left(Z\1_A\,g(E Z)\right) \geq \frac{E Z}{2}\,g(E Z)\,.
$$
Let $U=2Z$ to get the result.
\end{proof}

It is fairly easy to translate these to mixing time bounds. For instance, if $f(a)=\sqrt{\frac{1-a}{a}}$ then by Theorem \ref{thm:main}, Corollary \ref{cor:main} and Theorem \ref{thm:profile-mixing} the $L^2$-mixing times (denoted by $\tau_2(\epsilon)$) are:
$$
\tau_2(\epsilon) \leq
\begin{cases}
\vspace{2ex}\displaystyle 
\left\lceil \frac{1}{2(1-\C_{\sqrt{a(1-a)}})}\,\log \frac{1-\pi_*}{\pi_*\epsilon^2}  \right\rceil    & 
  \textrm{in general} \\
\vspace{2ex}\displaystyle 
\left\lceil \int_{\pi_*}^{\frac{1}{1+\epsilon^2}} \frac{dr}{2r(1-r)(1-\C_{\sqrt{a(1-a)}}(r))} \right\rceil & 
  \textrm{if $r\left(1-\C_{\sqrt{a(1-a)}}\left(\frac{1}{1+r^2}\right)\right)$ is convex} \\
\displaystyle \left\lceil \int_{\frac{4\pi_*}{1+3\pi_*}}^{\frac{1}{1+\epsilon^2/4}} \frac{dr}{r(1-r)(1-\C_{\sqrt{a(1-a)}}(r))} \right\rceil &
  \textrm{in general}
\end{cases}
$$
By making the change of variables $x=\frac{r}{1+r}$
and applying a few pessimistic approximations one obtains a result
more strongly resembling average conductance bounds:
$$
\tau_2(\epsilon) \leq \left\{
\begin{array}{ll}
\vspace{2ex}\displaystyle
\left\lceil \frac{1}{1-\C_{\sqrt{a(1-a)}}}\,\log \frac{1}{\epsilon\sqrt{\pi_*}}  \right\rceil 
  & \textrm{in general} \\
\vspace{2ex}\displaystyle
\left\lceil \int_{\pi_*}^{1/\epsilon^2} \frac{dx}{2x(1-\C_{\sqrt{a(1-a)}}(x))} \right\rceil
  & \textrm{if $x\left(1-\C_{\sqrt{a(1-a)}}\left(\frac{1}{1+x^2}\right)\right)$ is convex} \\
\displaystyle \left\lceil \int_{4\pi_*}^{4/\epsilon^2}
\frac{dx}{x(1-\C_{\sqrt{a(1-a)}}(x))} \right\rceil
  &  \textrm{in general} 
\end{array}
\right.
$$

It is often unnecessary to compute $\C_f(r)$ for $r>1/2$. Observe that $u$ almost everywhere $(A^c)_u = (A_{1-u})^c$. It follows that if $f(a)=f(1-a)$ then 
\begin{equation} \label{eqn:symmetry}
\C_f(A)=\int_0^1 \frac{f(\pi(A_{1-u}))}{f(\pi(A))}\,du 
  = \int_0^1 \frac{f(\pi((A_{1-u})^c))}{f(\pi(A^c))}\,du
  = \int_0^1 \frac{f(\pi((A^c)_u))}{f(\pi(A^c))}\,du = \C_f(A^c)\,.
\end{equation}
In particular, $\forall r\geq 1/2:\,\C_f(r)=\C_f(1/2)=\max_{\pi(A)\leq 1/2}\C_f(A)$.

\begin{remark} \label{rmk:sizes}
Mixing time bounds implied by the theorems of this section follow easily for the other distances, but for instance with $\C_{a(1-a)}$ for total variation distance and $\C_{a\log(1/a)}$ for relative entropy. However, it is often better to work with a harder distance, such as bounding total variation mixing ($\tau_{TV}(\epsilon)$) by instead bounding $L^2$ mixing ($\tau_2(\epsilon)$) and applying the relation $\tau_{TV}(\epsilon) \leq \tau_2(2\epsilon)$. The quantities are related by $\C_{a\log(1/a)}(A) \leq \frac 12\,(1+\C_{a(1-a)}(A))$ (see remarks after Theorem \ref{thm:gradients}) and $\C_{\sqrt{a(1-a)}}(A) \leq \sqrt{\C_{a(1-a)}(A)}$ (Cauchy-Schwartz), so generally the relative entropy or $L^2$-mixing bounds are less than a factor two worse than the total variation bound. In contrast, the lazy walk on a binary cube $\{0,1\}^d$ has tiny $1-\C_{a(1-a)}(\{x\}) \approx \frac{d+1}{2}\,2^{-d}$, but huge $1-\C_{\sqrt{a(1-a)}}(\{x\}) \approx \frac 12-\frac{1}{2\sqrt{d}}$, so the $L^2$ bounds will give much better asymptotics for this example.
\end{remark}


\section{Continuous Time}

Not much need be changed for continuous time. Let $H_t = e^{-t(\I-\P)}$ denote the continuous time Markov chain at time $t$. It is easily verified that if $\hat\K_t=e^{-t(\I-\hat\K)}$ then
$$
H_t(x,y) = \hat\E_t \pi_{S_t}
$$
where $S_0=\{x\}$ and $\hat\E_t$ is the expectation under the walk $\hat\K_t$. Bounds involving $\P^n(x,y)$ then translate directly into bounds in terms of $H_t(x,y)$. Once Lemma \ref{lem:discreteCongestion} is replaced by
$$
\frac{d}{dt} \hat\E_t f(\pi(S_t)) = -\hat\E_t f(\pi(S_t))(1-\C_{af(a)}(S_t))
$$
then mixing time bounds also carry over to the continuous-time case, although it is no longer necessary to approximate by a derivative at \eqref{eqn:derivative} nor necessary to take the ceiling of the bounds.


\chapter{Conductance and Modified Conductance} \label{sec:conductance}

The most common geometric tool for studying mixing time is the
{\em conductance} $\Phi$, a measure of the chance of leaving a set
after a single step. Such bounds have been shown only for $L^2$ mixing time. In this section we show bounds on $f$-congestion in terms of conductance for lazy walks, the most common situation. The real innovation of this section, however, is the modified conductance, a new quantity which is equivalent to conductance for a lazy walk in $L^2$ distance, but which also applies to walks with no holding probability and to other distances as well. We finish the section with a discussion on how our evolving set bounds are effected by changes in edge or vertex-expansion, or through re-scaling the transition kernel when slowing down a walk to increase its laziness.


\section{Conductance}

Let us begin with a formal definition of conductance.

\begin{definition}
The conductance $\~\Phi$ and conductance profile $\~\Phi(r)$ are given by
$$
\forall A\subset V:\,\~\Phi(A)=\frac{\Q(A,A^c)}{\pi(A)\pi(A^c)}\,,
\quad
\forall r>0:\,\~\Phi(r)=\min_{\pi(A)\leq r}\~\Phi(A)\,,
\quad
\~\Phi = \~\Phi(1/2) = \min_{A\subset V} \~\Phi(A)\,.
$$
The conductance $\Phi$ and conductance profile $\Phi(r)$ are defined similarly, but in terms of $\Phi(A)=\frac{\Q(A,A^c)}{\min\{\pi(A),\pi(A^c)\}}$. When necessary, notation such as $\Phi_{\K}$ will be used to denote conductance for Markov chain $\K$.
\end{definition}

The conductance profile $\~\Phi(r)$ can also be used to
upper bound the various $f$-congestion quantities $\C_f$ when the
Markov chain is lazy. The argument is not hard (see also \cite{MP05.1}).

\begin{theorem} \label{thm:conductance-congestion}
Given a lazy Markov chain, and $f$ concave, then
$$
\C_f(A) \leq \frac{f(\pi(A)+2\Q(A,A^c))+f(\pi(A)-2\Q(A,A^c))}{2f(\pi(A))}\,.
$$
\end{theorem}

\begin{proof}
For a lazy chain, if $u>1/2$ then $A_u\subset A$, and so
$$
\int_{1/2}^1 \pi(A_u)\,du = \sum_{y\in A} \left(\frac{\Q(A,y)}{\pi(y)}-\frac 12\right)\pi(y) = \Q(A,A)-\frac{\pi(A)}{2} = \frac{\pi(A)}{2} - \Q(A,A^c)\,.
$$
By the Martingale property $\int_0^1\pi(A_u)du=\pi(A)$ it follows that
$$
\int_0^{1/2}\pi(A_u)\,du = \pi(A)-\int_{1/2}^1\pi(A_u)du = \frac{\pi(A)}{2}+\Q(A,A^c)\,.
$$
Recall Jensen's inequality, that $\int g\circ h(u)\,du\leq g(\int h(u)\,du)$ if $u$ is a probability distribution and $g$ is concave. By concavity of $f$,
$$
\C_f(A) = \frac{\int_0^{1/2} f(\pi(A_u))\,\frac{du}{1/2}+\int_{1/2}^1 f(\pi(A_u))\,\frac{du}{1/2}}{2f(\pi(A))}
  \leq \frac{f\left(\int_0^{1/2}\pi(A_u)\,\frac{du}{1/2}\right) + f\left(\int_{1/2}^1 \pi(A_u)\frac{du}{1/2}\right)}
            {2f(\pi(A))}\,.
$$
\end{proof}

For each choice of $f$ a bit of simplification leads to bounds on
$\C_f$. For instance,  a lazy Markov chain will have
\begin{equation} \label{eqn:P-mixing}
\C_{\sqrt{a(1-a)}}(A) \leq \sqrt{1-\~\Phi(A)^2}
\quad\textrm{and}\quad \tau_2(\epsilon) \leq \left\lceil
\frac{2}{\~\Phi^2}\log\frac{1}{\epsilon\sqrt{\pi_*}}\right\rceil\,.
\end{equation}
See proof of Theorem \ref{thm:profile} for a similar calculation.
A conductance bound for a non-lazy walk will be considered later.


\section{Modified conductance} \label{sec:modified-conductance}

While the conductance has proven useful for studying lazy walks, if the chain is not lazy then the conductance $\~\Phi(r)$ is not useful for studying mixing. Consider the simple random walk on the complete bipartite graph $K_{m,m}$, a periodic Markov chain. Every subset $A\subset K_{m,m}$ has many edges to $A^c$ so conductance is large, but if $A$ is one of the bipartitions then a Markov chain starting in $A$ will bounce from $A$ to $A^c$ and back again, but it will never mix. 

The problem here is that the Markov chain never grows into a larger set, but is always stuck in half of the space. Therefore, it seems more appropriate to consider how much flow from $A$ reaches a strictly larger set, that is the worst flow into a set $B$ where $\pi(B)=\pi(A^c)$. In particular, we consider $\Psi(A)=\Psi(A,\pi(A^c))$ where
\begin{equation} \label{eqn:psiA}
\Psi(A,t) = \min_{\substack{B\subset V,\,v\in V,\\\pi(B)\leq t,\,\pi(B\cup v)> t}}
              \Q(A,B)+(t-\pi(B))\,\frac{\Q(A,v)}{\pi(v)}
\end{equation}
is the smallest flow from $A$ to a set of size $t$. For a lazy chain the minimum in $\Psi(A)$ occurs at $B=A^c$, so $\Psi(A)=\Q(A,A^c)$. In general, if $\pi$ is uniform then $\Psi(A)$ simplifies to $\Psi(A)=\min_{\pi(B)=\pi(A^c)} \Q(A,B)$. 

It is now possible to define the set quantity that is the main innovation of this section.

\begin{definition}
The {\em modified conductance} $\~\phi$ and modified conductance profile $\~\phi(r)$ are given by
$$
\~\phi(A) = \frac{\Psi(A)}{\pi(A)\pi(A^c)}\,,\quad
\~\phi(r) = \min_{\pi(A)\leq r} \~\phi(A)\,,\quad
\~\phi = \~\phi(1/2) = \min_{A\subset V} \~\phi(A)\,.
$$
Define $\phi(A)$ similarly but without $\pi(A^c)$ in the denominator.
\end{definition}

For a lazy chain $\Psi(A)=\Q(A,A^c)$ and so $\~\phi(A)=\~\Phi(A)$, and modified conductance
extends conductance to the non-lazy case.  The modified conductance
captures important properties quite well. For instance, a connected reversible chain has $\Psi(A)=0$ if 
and only if $A$ is one of the bipartitions of a periodic walk; the minimum in $\Psi(A)$ is then 
achieved by $B=A$, and $\Psi(A)=\Q(A,A)=0$ rather than $\Q(A,A^c)>0$ as with conductance.

An alternate interpretation of $\Psi(A)$ is as follows. Given a set $A\subset V$ let $\wp_A\in[0,1]$ satisfy
$$
\inf\{y:\,\pi(A_y)\leq\pi(A)\}\leq \wp_A \leq \sup\{y:\,\pi(A_y)\geq\pi(A)\}\,.
$$
The set $V\setminus A_{\wp_A}$ contains the vertices with minimum flow from $A$, and so
if $u<\wp_A$ then $\pi(A_u)-\pi(A)=\pi(\{y\in V\setminus A_{\wp_A}:\,\Q(A,y)\geq u\pi(y)\})$. 
It follows that
\begin{equation} \label{eq:psi}
\Psi(A)=\int_0^{\wp_A} (\pi(A_u)-\pi(A))\,du=\int_{\wp_A}^1 (\pi(A)-\pi(A_u))\,du = \frac 12\,\int_0^1 |\pi(A)-\pi(A_u)|\,du \,,
\end{equation}
where the first equality is from the definition of $\Psi(A)$ and the second is from Lemma \ref{lem:martingale}.
Since $u$-almost everywhere $A_u=((A^c)_{1-u})^c$ the final equality shows that $\Psi(A)=\Psi(A^c)$,
a property which is also satisfied by conventional set expansion with $\Q(A,A^c)=\Q(A^c,A)$.


\section{An Inequality Prover}

With this formulation of $\Psi(A)$ it is possible to upper and lower bound each $\C_f(A)$ via Jensen's inequality, although the upper bounds require a careful setup. However, an argument based on Lemma \ref{lem:worst-case} is more appealing because it gives a general result for all concave $f$, and it immediately implies sharpness due to the explicit 
constructions \eqref{eqn:upper_sets} and \eqref{eqn:lower_sets}. We find it also to be more intuitive than Jensen, as it relates a graphical representation of $\pi(A_u)$ directly to the $f$-congestion.

%
%

\begin{lemma} \label{lem:worst-case}
Suppose that $f:[0,1]\rightarrow \R$ is concave, and $g,\,\hat{g}:\,[0,1]\rightarrow[0,1]$ are functions such that $\hat{g}$ is non-increasing and $(g-\hat{g})$ is continuous on a measure one open subset of $(0,1)$ (e.g. a step function). If 
$$
\forall t\in[0,1]:\,\int_0^t g(u)\,du \geq \int_0^t \hat{g}(u)\,du
$$
and $\int_0^1 g(u)\,du=\int_0^1 \hat{g}(u)\,du$, then
$$
\int_0^1 f\circ g(u)\,du \leq \int_0^1 f\circ \hat{g}(u)\,du\,.
$$
\end{lemma}


\begin{proof}
First we use concavity. Given $\delta,c>0$ and $x,y\in[0,1]$ with $1-\delta\geq x\geq y\geq\delta c^{-1}$, let $\lambda = \frac{\delta/c}{x-y+(1+c^{-1})\delta}\in(0,1)$. Then $x=\lambda c\,(y-\delta c^{-1})+(1-\lambda c)\,(x+\delta)$ and $y=(1-\lambda)\,(y-\delta c^{-1})+\lambda\,(x+\delta)$. By concavity, $f(x)\geq \lambda c\,f(y-\delta c^{-1})+(1-\lambda c)\,f(x+\delta)$ and $f(y)\geq (1-\lambda)\,f(y-\delta c^{-1})+\lambda\,f(x+\delta)$. It follows that
\begin{equation} \label{eqn:bigup_smalldown}
f(x)+c\,f(y) \geq f(x+\delta) + c\,f(y-\delta c^{-1})\,.
\end{equation}

Now, let $S_1=\{u\in(0,1)\,:\,g(u)>\hat{g}(u)\}$, $S_2=\{u\in(0,1)\,:\,g(u)=\hat{g}(u)\}$ and $S_3=\{u\in(0,1)\,:\,g(u)<\hat{g}(u)\}$. If $(g-\hat{g})$ is continuous on an open interval $I$ then $\{u\in I\,:\,g(u)>\hat{g}(u)\}$ and $\{u\in I\,:\,g(u)<\hat{g}(u)\}$ are open subsets of $I$. Hence $S_1$ and $S_3$ are open subsets of $(0,1)$, since $(g-\hat{g})$ is continuous on a measure one countable union of disjoint open intervals. More specifically, each is a countable union of disjoint open intervals. Define
$$
h(u) = \sup\left\{ s\in[0,1]\,:\, \int_0^u (g-\hat{g})^+(x)\,dx \geq \int_0^s (g-\hat{g})^-(x)\,dx \right\}
$$
where $F^{\pm}(x)=\max\{0,\pm F(x)\}$. 
Then $h$ is a bijection from $S_1\setminus S$ to $S_3$, where $S=\{u\in S_1\,:\,h(u)\in \overline{S_3}\setminus S_3\}$ is countable. Moreover, when $u\in S_1\setminus S$ then $h'(u)=\frac{(g-\hat{g})(u)}{(\hat{g}-g)(h(u))}>0$. Hence, via integration by substitution,
$$
\int_{S_3} (f\circ\hat{g})(x)\,dx = \int_{S_1\setminus S} (f\circ\hat{g})(h(u))\,h'(u)\,du\,,
$$
and likewise with $g$ in place of $\hat{g}$. The lemma then follows from integration:
\begin{eqnarray*}
\int_0^1 f\circ\hat{g}(u)\,du &=& \int_{S_1} + \int_{S_2} + \int_{S_3} (f\circ\hat{g})(u)\,du \\
  &=&    \int_{S_1\setminus S} \left( (f\circ \hat{g})(u) + h'(u)\;(f\circ\hat{g})(h(u))\right)\,du + \int_{S_2} (f\circ \hat{g})(u)\,du \\
  &\geq& \int_{S_1\setminus S} \left( (f\circ g)(u) + h'(u)\;(f\circ g)(h(u))\right)\,du + \int_{S_2} (f\circ g)(u)\,du \\ 
  &=&   \int_{S_1} + \int_{S_2} + \int_{S_3} (f\circ g)(u)\,du = \int_0^1 f\circ g(u)\,du
\end{eqnarray*}
The inequality was an application of \eqref{eqn:bigup_smalldown} with $x=\hat{g}(u)$, $y=\hat{g}(h(u))$, $c=h'(u)$, and $\delta=g(u)-\hat{g}(u)=c\,(\hat{g}(h(u))-g(h(u)))$.
\end{proof}

The lemma implies that for any set $A\subset V$, and for some initial conditions,
if there are non-increasing step functions $m,\,M:\,[0,1]\mapsto[0,1]$ such that 
\begin{eqnarray} 
\forall t\in[0,1]:\,\int_0^t M(u)\,du \geq \int_0^t \pi(A_u)\,du \geq \int_0^t m(u)\,du \label{eqn:extreme} \\
\textrm{and}\quad
\int_0^1 M(u)\,du = \int_0^1 \pi(A_u)\,du = \int_0^1 m(u)\,du \nonumber
\end{eqnarray}
then for every concave function $f(x)$ it follows that
$$
\frac{\int_0^1 f(M(u))\,du}{f(\pi(A))} \leq \C_f(A) \leq \frac{\int_0^1 f(m(u))\,du}{f(\pi(A))}\,.
$$

\begin{figure}[ht]
\begin{center}
\includegraphics[height=1.5in]{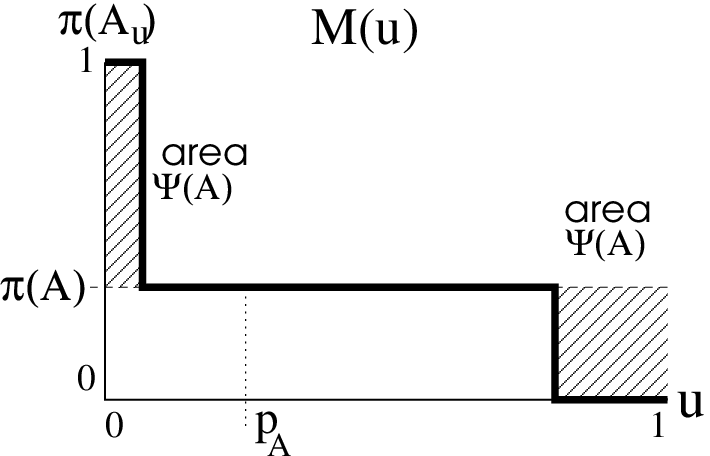}
\hspace{1in}
\includegraphics[height=1.5in]{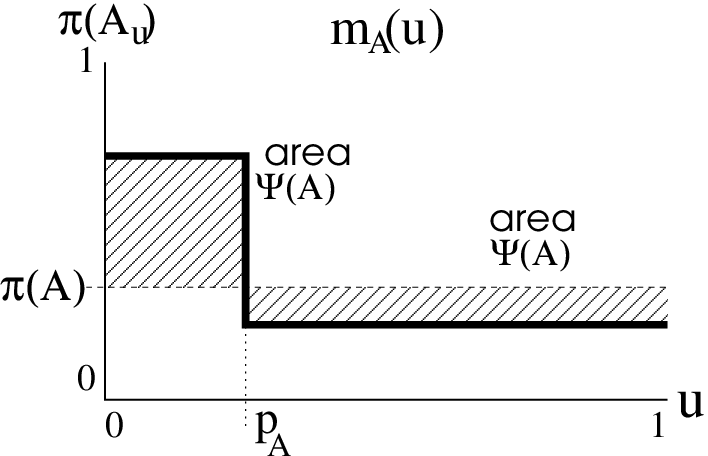}
\caption{Distributions such that $\int_0^t M(u)\,du \geq \int_0^t \pi(A_u)\,du \geq \int_0^t m(u)\,du$
    given $\Psi(A)$ and $\wp_A$.} \label{fig:extreme}
\end{center}
\end{figure}

In the problem at hand, $\pi(A_u)\in[0,1]$ is non-increasing and equation (\ref{eq:psi}) implies 
$\Psi(A)$ is the area below $\pi(A_u)$ and above $\pi(A)$,
and also above $\pi(A_u)$ and below $\pi(A)$. The extreme cases of $\pi(A_u)$ can be
drawn immediately, as in Figure \ref{fig:extreme}.

\section{Bounds on $f$-congestion $\C_f(A)$}

We now show modified conductance bounds on some of the $f$-congestion quantities of interest.

\begin{theorem} \label{thm:profile}
Given a subset $A\subset  V$ then
$$
\begin{array}{rclcl}
\vspace{1ex}
\~\phi(A) &\geq& 1-\C_{\sqrt{a(1-a)}}(A) &\geq& \displaystyle 1-\sqrt{1-\~\phi(A)^2} \geq \~\phi(A)^2/2 \\
\vspace{1ex}
\~\phi(A) &\geq& 1-\C_{a\log(1/a)}(A)      &\geq& \displaystyle \frac{2\phi(A)^2}{\log(1/\pi(A))} \\
\~\phi(A) &\geq & 1-\C_{a(1-a)}(A) &\geq& 4\~\phi(A)^2\pi(A)(1-\pi(A))
\end{array}
$$
\end{theorem}

\begin{proof}
For the upper bound, Figure \ref{fig:extreme} shows that, given $\Psi(A)$ then $\forall t\in[0,1]:\,\int_0^t M(u)\,du\geq\int_0^t\pi(A_u)\,du$ and $\int_0^1 M(u)\,du=\pi(A)=\int_0^1 \pi(A_u)\,du$, where
\begin{equation} \label{eqn:upper_sets}
M(u) =
\begin{cases}
0      & \textrm{if } u > 1 - \frac{\Psi(A)}{\pi(A)} \\
\pi(A) & \textrm{if } u\in \left(\frac{\Psi(A)}{1-\pi(A)},\,1-\frac{\Psi(A)}{\pi(A)}\right] \\
1      & \textrm{if } u \leq \frac{\Psi(A)}{1-\pi(A)}
\end{cases}
\end{equation}
By Lemma \ref{lem:worst-case} any choice of $f(z)$ which is concave and non-negative will therefore satisfy
\begin{eqnarray*}
\C_f(A) &\geq& \frac{\int_0^1 f(M(u))\,du}{f(\pi(A))} \\
        &=& \frac{\Psi(A)}{\pi(A)}\,\frac{f(0)}{f(\pi(A))}
           + \left(1-\frac{\Psi(A)}{\pi(A)\pi(A^c)}\right)\,\frac{f(\pi(A))}{f(\pi(A))}
           + \frac{\Psi(A)}{1-\pi(A)}\,\frac{f(1)}{f(\pi(A))} \\
        &\geq& 1-\~\phi(A)
\end{eqnarray*}
This shows all of the upper bounds.

To prove lower bounds, suppose $\wp_A$ and $\Psi(A)$ are known. Then Figure \ref{fig:extreme} demonstrates that $\forall t\in[0,1]:\,\int_0^t \pi(A_u)\,du \geq \int_0^t m(u)\,du$ and $\int_0^1 m(u)\,du=\pi(A)=\int_0^1 \pi(A_u)\,du$, where
\begin{equation} \label{eqn:lower_sets}
m(u) =
\begin{cases}
\pi(A) - \frac{\Psi(A)}{1-\wp_A} & \textrm{if } u > \wp_A \\
\pi(A) + \frac{\Psi(A)}{\wp_A}   & \textrm{if } u < \wp_A
\end{cases}
\end{equation}

All that remains is to substitute this into the formula for $1-\C_f(A)$ for the various $f(x)$ of interest,
and then minimize over all possible $\wp_A\in[0,1]$.

The bound on $1-\C_{a(1-a)}$ is the easiest. Apply Lemma \ref{lem:worst-case} with $f(z)=a(1-a)$ to obtain
\begin{eqnarray*}
\C_{a(1-a)}(A) &\leq& \frac{\int_0^1 f(\pi(A_u))\,du}{f(\pi(A))} \\
 &=& \wp_A\,\frac{\pi(A)+\frac{\Psi(A)}{\wp_A}}{\pi(A)}\,\frac{1-\pi(A)-\frac{\Psi(A)}{\wp_A}}{1-\pi(A)}
            + (1-\wp_A)\,\frac{\pi(A)-\frac{\Psi(A)}{1-\wp_A}}{\pi(A)}\,\frac{1-\pi(A)+\frac{\Psi(A)}{1-\wp_A}}{1-\pi(A)} \\
 &=& 1-\frac{\~\phi(A)^2\,\pi(A)\pi(A^c)}{\wp_A(1-\wp_A)} \leq 1-4\,\~\phi(A)^2\,\pi(A)\pi(A^c)
\end{eqnarray*}

For the lower bound on $\C_{a\log(1/a)}$ proceed similarly.
\begin{eqnarray*}
\lefteqn{\C_{a\log(1/a)}(A)} \\
 &\leq& -\wp_A    \,\frac{\pi(A)+\frac{\Psi(A)}{\wp_A}}{\pi(A)\log \frac{1}{\pi(A)}}\,\log \left(\pi(A)+\frac{\Psi(A)}{\wp_A}\right)        -(1-\wp_A)\,\frac{\pi(A)-\frac{\Psi(A)}{1-\wp_A}}{\pi(A)\log \frac{1}{\pi(A)}}\,\log \left(\pi(A)-\frac{\Psi(A)}{1-\wp_A}\right) \\
 &=&   1
       - \frac{\wp_A+\phi(A)}{\log \frac{1}{\pi(A)}}\,\log \frac{\wp_A+\phi(A)}{\wp_A}
       - \frac{1-\wp_A-\phi(A)}{\log \frac{1}{\pi(A)}}\,\log \frac{1-\wp_A-\phi(A)}{1-\wp_A}
\end{eqnarray*}
Then $(1-\C_{a\log(1/a)}(A))\,\log(1/\pi(A)) \geq g(\wp_A,\phi(A)) \geq 2\phi(A)^2$ by Lemma \ref{lem:log-inequality}.

Now for $\C_{\sqrt{a(1-a)}}(A)$. Applying Lemma \ref{lem:worst-case} and equation \eqref{eqn:lower_sets} as before,
\begin{eqnarray*}
\lefteqn{\C_{\sqrt{a(1-a)}}(A)} \\
&\leq& \mbox{$\wp_A\,\sqrt{\left(1 + \frac{\Psi(A)}{\wp_A\,\pi(A)}\right)    \left(1 - \frac{\Psi(A)}{\wp_A\,\pi(A^c)}\right)}
   + (1-\wp_A)\,\sqrt{\left(1 - \frac{\Psi(A)}{(1-\wp_A)\,\pi(A)}\right)\left(1 + \frac{\Psi(A)}{(1-\wp_A)\,\pi(A^c)}\right)}$}
\\
&=& \mbox{$\sqrt{\left(\wp_A+\~\phi(A)\,\pi(A^c)\right)    \left(\wp_A-\~\phi(A)\,\pi(A)\right)}
     + \sqrt{\left(1-\wp_A - \~\phi(A)\,\pi(A^c)\right)\left(1-\wp_A + \~\phi(A)\,\pi(A)\right)}$}
\end{eqnarray*}
The bound on $\C_{\sqrt{a(1-a)}}(A)$ follows from Lemma \ref{eqn:chi2-inequality} with $X = \displaystyle \wp_A + \~\phi(A)\,\pi(A^c)$ and $Y = \displaystyle \wp_A - \~\phi(A)\,\pi(A)$.
\end{proof}

It follows, for instance, that
\begin{equation} \label{eqn:L2-modified-conductance}
\tau_2(\epsilon) \leq \left\lceil \frac{2}{\~\phi^2}\log\frac{1}{\epsilon\sqrt{\pi_*}}\right\rceil
\quad\textrm{and}\quad
\tau_2(\epsilon) \leq \left\lceil \int_{4\pi_*}^{4/\epsilon^2} \frac{2\,dr}{r\,\~\phi(r)^2} \right\rceil\,.
\end{equation}

Conductance can be used to obtain a crude lower bound on the modified conductance.

\begin{lemma} \label{lem:non-lazy}
For an ergodic Markov chain, if $\forall x\in V:\,\P(x,x)\geq\gamma\in[0,1]$
and $A\subset V$ then
$$
\~\Phi(A) \geq \~\phi(A) \geq \min\left\{1,\,\frac{\gamma}{1-\gamma}\right\}\,\~\Phi(A)\,.
$$
\end{lemma}

\begin{proof}
The upper bound $\~\phi(A)\leq\~\Phi(A)$ is trivial because $\Psi(A)\leq\Q(A,A^c)$. The minimum in the lower bound is equal to $1$ exactly when $\gamma\geq 1/2$, but in this case $\~\phi(A)=\~\Phi(A)$, so this case is also trivial. It remains to consider the lower bound when $\gamma<1/2$.

In the definition of $\Psi(A)$ there
is a set $B$, and one vertex $v$ for which only a $\frac{\pi(A^c)-\pi(B)}{\pi(v)}$ fraction is
counted. Extend the state space $V$ to a space $V'$ by splitting $v$ into two vertices $v_1$ and $v_2$,
with $v_1$ of size $\pi(A^c)-\pi(B)$, $v_2$ with the remainder, and ergodic flows into $v_1$ and $v_2$
determined by their respective sizes.
Then let $C=B\cup v_1$ be the set where $\Psi(A)=\Q(A,C)$. It follows that
\begin{eqnarray*}
\Psi(A) &=& \Q(A,C\cap A^c) + \Q(A,C\cap A) \\
   &\geq& \Q(A,C\cap A^c) + \gamma\,\pi(C\cap A) \\    
   &\geq& \Q(A,C\cap A^c) + \gamma\,\frac{\Q(A,A^c\setminus C)}{1-\gamma} \\
   &\geq& \frac{\gamma}{1-\gamma}\,\Q(A,A^c)
\end{eqnarray*}
The first inequality uses the fact that $\forall v\in A:\,\Q(v,v)\geq \gamma\pi(v)$
and so $\Q(A,v)\geq \gamma\,\pi(v)$. The second inequality is because
$\pi(C\cap A)=\pi(A^c\setminus C)\geq  \frac{\Q((A^c\setminus C)^c,A^c\setminus C)}{1-\gamma}  \geq \frac{\Q(A,A^c\setminus C)}{1-\gamma}$.
\end{proof}

The $\gamma/(1-\gamma)$ factor is introduced when converting $C\cap A$ into a subset of $A^c$, in short primarily because $\Q(A,A^c)$ is not the correct quantity to work with for non-lazy chains. This induces a mixing bound in terms of conductance for non-lazy walks, but this will be substantially improved on later.


%

Although Lemma \ref{lem:worst-case} was just used to show the bounds of Theorem \ref{thm:profile}, the arguments are easily modified to use Jensen's Inequality instead (see \cite{Mon07.1} for just such an approach). However, the upper bound of $\~\phi(A)$ is somewhat more subtle, and so we give here a proof with Jensen's inequality as well.

\begin{proof}[Proof of upper bounds in Theorem \ref{thm:profile} using only Jensen's Inequality]
In the definition of $\Psi(A)$ there
is a set $B$, and one vertex $v$ for which only a $\frac{\pi(A)-\pi(B)}{\pi(v)}$ fraction is
counted. Extend the state space $V$ to a space $V'$ by splitting $v$ into two vertices $v_1$ and $v_2$,
with $v_1$ of size $\pi(A)-\pi(B)$, $v_2$ with the remainder and flows adjusted accordingly.
Then let $C=B\cup v_1$ be the set where $\Psi(A)$ is achieved in the space $V'$.

Observe that when $u>\wp_A$ then
$\pi(A_u) = \sum_{y\in C} \delta_{\Q(A,y)\geq u\cdot\pi(y)}\,\pi(A)\,\frac{\pi(y)}{\pi(A)}$.
Since $\pi(C)=\pi(A)$ and $f(x)=\sqrt{x(1-x)}$ is concave then by Jensen's inequality
\begin{eqnarray*}
\int_{\wp_A}^1 \sqrt{\frac{\pi(A_u)(1-\pi(A_u))}{\pi(A)(1-\pi(A))}}\,du
  &\geq& \int_{\wp_A}^1 \sum_{y\in C}
             \frac{\pi(y)}{\pi(A)}\,
                \sqrt{\frac{\delta_{\Q(A,y)\geq u\cdot\pi(y)}\pi(A)\,(1-\delta_{\Q(A,y)\geq u\cdot\pi(y)}\pi(A))}
             {\pi(A)(1-\pi(A))}}\,du \\
  &=& \int_{\wp_A}^1 \sum_{y\in C} \frac{\pi(y)}{\pi(A)}\,\delta_{\Q(A,y)\geq u\cdot\pi(y)}\,du \\
  &=& \int_{\wp_A}^1 \frac{\pi(A_u)}{\pi(A)}\,du = 1 - \int_0^{\wp_A} \frac{\pi(A_u)}{\pi(A)}\,du\,,
\end{eqnarray*}
where the final equality uses the Martingale property $\int_0^1 \pi(A_u)\,du=\pi(A)$.
Similarly, when $u\leq \wp_A$ then
$\pi(A_u) = \sum_{y\in V'\setminus C} \left(\pi(A)+\delta_{\Q(A,y)\geq u\cdot\pi(y)}\,\pi(A^c)\right)\,\frac{\pi(y)}{\pi(A^c)}$,
and because $\pi(V'\setminus C)=\pi(A^c)$ then Jensen's inequality can be applied once more to
obtain
$$
\int_0^{\wp_A} \sqrt{\frac{\pi(A_u)(1-\pi(A_u))}{\pi(A)(1-\pi(A))}}\,du \geq \int_0^{\wp_A} \frac{1-\pi(A_u)}{1-\pi(A)}\,du\,.
$$
Combining these inequalities and rearranging a bit gives the result
\begin{eqnarray*}
\C_{\sqrt{a(1-a)}}(A) &\geq& 1 - \int_0^{\wp_A} \frac{(1-\pi(A))\pi(A_u)-\pi(A)(1-\pi(A_u))}{\pi(A)\pi(A^c)}\,du \\
   &=& 1 - \int_0^{\wp_A} \frac{\pi(A_u)-\pi(A)}{\pi(A)\pi(A^c)}\,du = 1-\~\phi(A)\,.
\end{eqnarray*}

The only properties used in this argument were the concavity of $f(x)=\sqrt{x(1-x)}$
and that $f(0)=f(1)=0$. These also hold for $1-\C_{a(1-a)}(A)$ (with $f(x)=x(1-x)$) and
$1-\C_{a\log(1/a)}(A)$ (with $f(x)=x\log (1/x)$), giving the upper bounds for these.
\end{proof}

The following two inequalities were used in the proof of Theorem \ref{thm:profile}:

\begin{lemma} \label{eqn:chi2-inequality}
If $X,\,Y\in[0,1]$ then
$$
g(X,Y)=\sqrt{X\,Y}+\sqrt{(1-X)(1-Y)} \leq \sqrt{1-(X-Y)^2}\,.
$$
\end{lemma}

\begin{proof}
Observe that
$$
g(X,Y)^2 = 1-(X+Y)+2\,X\,Y + \sqrt{[1-(X+Y)+2\,X\,Y]^2 - [1-2(X+Y)+(X+Y)^2]}\,.
$$

Now, $\sqrt{A^2-B} \leq A-B$ if $A^2\geq B$, $A\leq\frac{1+B}{2}$ and $A\geq B$ (square both sides to show this). 
These conditions are easily verified with $A=1-(X+Y)+2\,X\,Y$ and $B=1-2(X+Y)+(X+Y)^2$, and so
\begin{eqnarray*}
g(X,Y)^2 &\leq& 2\left[1-(X+Y)+2\,X\,Y\right] - \left[1-2(X+Y)+(X+Y)^2\right] \\
       &=&    1+2\,X\,Y-X^2-Y^2 = 1-(X-Y)^2
\end{eqnarray*}
\end{proof}

\begin{lemma} \label{lem:log-inequality}
If $x\in(0,1)$ and $y\in[0,1-x)$ then
$$
g(x,y) = (x+y)\log \frac{x+y}{x} + (1-x-y)\log \frac{1-x-y}{1-x}\geq 2y^2\,.
$$
\end{lemma}

\begin{proof}
Start by seeing what can be shown by differentiation.
\begin{eqnarray*}
\frac{dg}{dx} &=& \log\left(1+\frac{y}{x}\right)-\frac{y}{x}
      -\log\left(1-\frac{y}{1-x}\right) - \frac{y}{1-x} \\
\frac{d^2g}{dx^2} &=& y^2\,\frac{(x+y)x^2+(1-x)^2(1-(x+y))}
                             {x^2(1-x)^2(x+y)(1-(x+y))} \geq 0
\end{eqnarray*}
It follows that $g(x,y)$ is convex with respect to $x$, and since $\left.\frac{dg}{dx}\right|_{x=(1-y)/2}\leq 0$ and $\left.\frac{dg}{dx}\right|_{x=1/2}\geq 0$ then the minimum occurs at some $x\in[(1-y)/2,\,1/2]$.

To lower bound the minimum we first lower bound $g(x,y)$. By the inequality $f(z)=z\log z + (1-z)\log (1-z)\geq -\log 2 + 2(z-1/2)^2$ when $z\in[0,1]$ it follows that
$$
g(x,y) = f(x+y) - \log(1-x) + (x+y)\log \frac{1-x}{x} \geq h(x,y)
$$
where $h(x,y) = -\log 2 + 2(x+y-1/2)^2 - \log(1-x) + (x+y)\log \frac{1-x}{x}$.
Now,
\begin{eqnarray*}
\frac{dh}{dx}   &=& y\,\left(4-\frac{1}{x(1-x)}\right) + 4x + \log \frac{1-x}{x} - 2 \\
\frac{d^2h}{dx} &=& y\,\frac{1-2x}{x^2(1-x)^2} + 4 - \frac{1}{x(1-x)} \\
\left.\frac{d^2h}{dx^2}\right|_{x=(1-c\,y)/2}
           &=& \frac{4\,c\,y^2\,(4-c) + 4\,c^4\,y^4}{(1-c^2\,y^2)^2}
\end{eqnarray*}
The second derivative is positive when $c\in[0,1]$, and so $h(x,y)$ is convex in $x$ when $x\in[(1-y)/2,1/2]$. However, $\left.\frac{dh}{dx}\right|_{x=1/2}=0$ and so $h(x,y)\geq h(1/2,y)=2y^2$ when $x\in[(1-y)/2,1/2]$.

It follows that
$\displaystyle g(x,y)\geq \min_{x\in[(1-y)/2,1/2]} g(x,y) \geq \min_{x\in[(1-y)/2,1/2]} h(x,y) \geq h(1/2,y)=2y^2$.
\end{proof}


\section{Flow distribution}

Now let us look at how flow distribution affects the $\C$ quantities. To do this we assume that we have two Markov chains which differ only in a single characteristic, related to either the level of ergodic flow or the way in which the ergodic flow is distributed among the vertices. The following applications of Lemma \ref{lem:worst-case} then give a good intuition into what governs mixing.

\begin{corollary}  \label{cor:expansion}
Suppose that $\MM$ and $\MM'$ are finite irreducible Markov chains with the
same stationary distribution $\pi$, that $A\subset V$ with $\Psi(A)=\Q(A,A^c)$ (e.g. any subset if $\MM$ is lazy), and that $f:\,[0,1]\rightarrow\R^+$ is a concave function. Then
$$
1-\C_f(A) \geq 1-\C_f'(A)
$$
if either of the following two conditions hold:
\begin{itemize}
\item {\em Edge expansion / flow:} $\MM'$ has smaller pointwise flow than $\MM$, that is
$$
\forall v\in A^c:\,\Q(A,v)\geq \Q'(A,v) \quad \textrm{and}\quad \forall v\in A:\,\Q(A^c,v)\geq \Q'(A^c,v)\,,
$$

\item {\em Vertex expansion:} $\MM'$ has less well distributed flow than $\MM$, that is
$$
\begin{array}{rclcl}
\vspace{1ex}
\forall u\in[0,1]&:&
\displaystyle \sum_{v\in A^c} \min\{u\pi(v),\,\Q(A,v)\} 
&\geq& \displaystyle \sum_{v\in A^c} \min\{u\pi(v),\,\Q'(A,v)\} \\
and \ \forall u\in[0,1]&:&
\displaystyle \sum_{v\in A} \min\{u\pi(v),\,\Q(A^c,v)\} 
&\geq& \displaystyle \sum_{v\in A}   \min\{u\pi(v),\,\Q'(A^c,v)\}
\end{array}
$$
and moreover $\Psi'(A)=\Q'(A,A^c)$ (e.g. if $\MM'$ is lazy) and $\Q(A,A^c)=\Q'(A,A^c)$ (i.e. equal ergodic flows). 
\end{itemize} 
\end{corollary}

The first relation says that, all other things being equal, if each vertex in $A^c$ gets less ergodic flow from $A$, and vice-versa, then the Evolving set mixing time bound will be worse. The converse does not hold, as the periodic walk on the two-point space has higher edge expansion than the lazy two-point walk, but does not mix.

In order to understand the second case we need to define exactly what is meant by vertex expansion. One reasonable definition is to say that the flow is well distributed if cutting it off at
some threshold does not cut off too much, that is if the threshold is $u$ then
$\sum_{v\in A^c} \min\{u\pi(v),\,\Q(A,v)\}$ is about the same size as $\Q(A,A^c)$, and likewise
with a sum over $v\in A$. The corollary then says that, all other conditions being equal, lower vertex expansion leads to a slower mixing time.

The requirement that $\Psi(A)=\Q(A,A^c)$ arises from the following lemma.

\begin{lemma} \label{lem:psi-flow}
Given an irreducible Markov chain and $A\subset V$,
then $\Psi(A)=\Q(A,A^c)$ if and only if  
$\bigcup_{u>\wp_A} A_u \subsetneq A \subseteq A_{\wp_A}$
for $\wp_A=\inf\{u:\pi(A_u)<\pi(A)\}$.
\end{lemma}

\begin{proof}
Observe that $\Psi(A)=\Q(A,A^c)$ if and only if
the set $B$ where the minimum occurs in the definition of $\Psi(A)$ can be taken as $B=A^c$.
This happens if and only if $\forall v\in A^c,\,v'\in A:\,\Q(A,v)/\pi(v)\leq \Q(A,v')/\pi(v')$,
which is in turn equivalent to $\pi(A_u)\geq\pi(A)$ if and only if $A\subseteq A_u$. 
This occurs if and only if $A\subseteq A_{\wp_A}$ (observe that $\pi(A_{\wp_A})\geq\pi(A)$) and
$\forall u>\wp_A:\, A_u\subsetneq A$.
\end{proof}

This shows that $\Psi(A)=\Q(A,A^c)$ if and only if the $A_u$ split into two types, everything at $u>\wp_A$ is in $A$ and everything dropped at $u<\wp_A$ is in $A^c$. Most properties of lazy Markov chains will hold for sets $A\subset V$ when $\Psi(A)=\Q(A,A^c)$. 

\begin{proof}[Proof of Corollary \ref{cor:expansion}]
Figure \ref{fig:expansion} gives a visual ``proof'' using Lemma \ref{lem:worst-case}. 
\begin{figure}[ht]
\begin{center}
\includegraphics[height=1.5in]{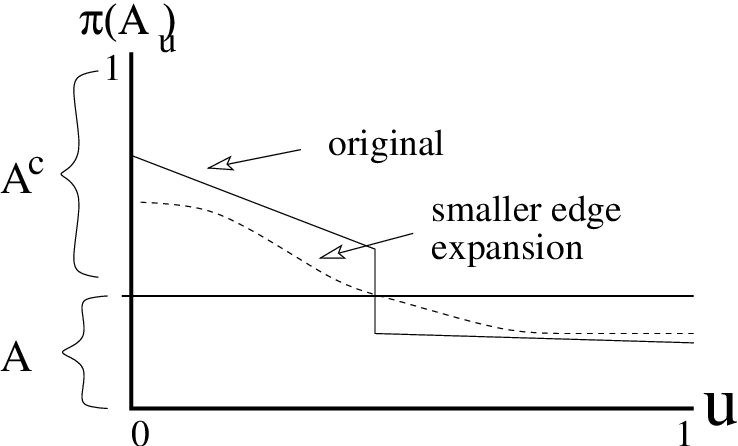}
\hspace{1in}
\includegraphics[height=1.5in]{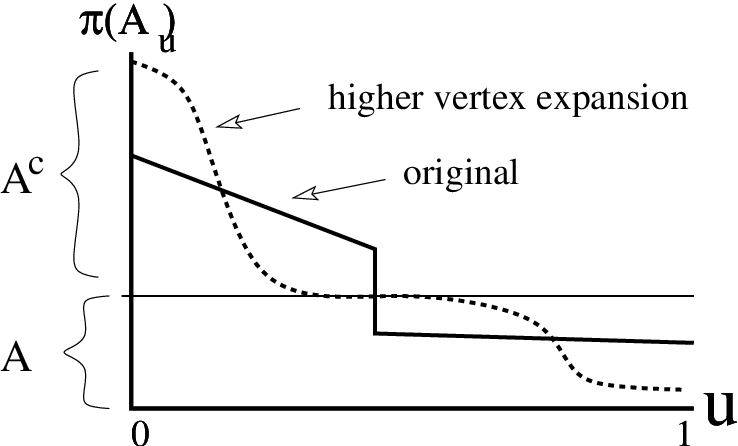}
\caption{The $\pi(A_u)$ for $\MM$ vs. smaller edge expansion $\MM'$, and vs. higher vertex expansion $\MM'$.}
\label{fig:expansion}
\end{center}
\end{figure}

Let us now show rigorously that Lemma \ref{lem:worst-case} can in fact be applied, as suggested by the pictorial representations.

First, edge-expansion. 

By Lemma \ref{lem:psi-flow}, if $u\leq\wp_A$ then $A_u = A \cup \left(\cup_{v\in A^c} \delta_{\Q(A,v)\geq u\pi(v)}\right)$,
while if $u>\wp_A$ then $A_u = \cup_{v\in A} \delta_{\Q(A,v)\geq u\pi(v)}$. The Markov chain $\MM'$ also splits
into cases of $u\leq\wp_A$ and $u>\wp_A$, because if $v\in A^c$ then $\Q'(A,v)\leq\Q(A,v)\leq\wp_A\pi(v)$, while
if $v\in A$ then $\Q'(A,v)\geq\Q(A,v)\geq\wp_A\pi(v)$, so we may assume $w(\pi(A))=\wp_A$ for $\MM'$ as well.

First consider the case that $t\leq\wp_A$. If $v\in A^c$ then $\Q(A,v)/\pi(v)\geq\Q'(A,v)/\pi(v)$, and so
$\pi(A_u)\geq\pi(A_u')$ if $u\leq\wp_A$, and in particular
$\forall t\in[0,\wp_A]:\,\int_0^t \pi(A_u)\,du\geq\int_0^t \pi(A_u')\,du$. 

Now consider the case when $t>\wp_A$. If $v\in A$ then $\Q(A,v)/\pi(v)\leq\Q'(A,v)/\pi(v)$ and so for $u\geq\wp_A$ it follows that
$\pi(A_u)\leq\pi(A_u')$, and therefore $\forall t\in[\wp_A,1]:\,\int_t^1 \pi(A_u)\,du \leq \int_t^1 \pi(A_u')\,du$.
But then, by the Martingale property Lemma \ref{lem:martingale}, 
$$
\int_0^t \pi(A_u)\,du 
  = \pi(A)-\int_t^1\pi(A_u)\,du
  \geq\pi(A)-\int_t^1\pi(A_u')\,du
  = \int_0^t\pi(A_u')\,du
$$

We have now established that $\forall t\in[0,1]:\,\int_0^t \pi(A_u)\,du \geq \int_0^t \pi(A_u')\,du$, and so the corollary follows from Lemma \ref{lem:worst-case}.
%

Now, vertex-expansion.

If $t\in[0,1]$ then 
\begin{eqnarray*}
\int_0^t \pi(A_u)\,du &=& 
\begin{cases} 
t - \sum_{v\in A^c} (\Q(A,v)-t\pi(v))^+ & \textrm{if } t\leq\wp_A \\
\pi(A)  - \sum_{v\in A}  (\Q(A^c,v)-(1-t)\pi(v))^+ & \textrm{if } t>\wp_A
\end{cases} \\
&=&
\begin{cases}
t - \sum_{v\in A^c} \Q(A,v)-\min\{\Q(A,v),t\pi(v)\} & \textrm{if } t\leq\wp_A \\
\pi(A)  - \sum_{v\in A} \Q(A^c,v)-\min\{\Q(A^c,v),(1-t)\pi(v)\} & \textrm{if } t>\wp_A
\end{cases} \\
&=&
\begin{cases}
t - \Q(A,A^c) + \sum_{v\in A^c} \min\{\Q(A,v),t\pi(v)\} & \textrm{if } t\leq\wp_A \\
\pi(A) - \Q(A,A^c) + \sum_{v\in A} \min\{\Q(A^c,v),(1-t)\pi(v)\} & \textrm{if } t>\wp_A
\end{cases}
\end{eqnarray*}
and likewise with $\int_0^t \pi(A_u')\,du$. By the conditions on vertex-expansion and the relation $\Q(A,A^c)=\Q'(A,A^c)$ it follows that $\int_0^t\pi(A_u)\,du\geq\int_0^\pi(A_u')\,du$, and so Lemma \ref{lem:worst-case} finishes the proof.
\end{proof}


\section{The effect of re-scaling on mixing time}

One feature of spectral gap/profile approaches to bounding mixing times is that they ``scale properly,'' in the sense that if the Markov chain is slowed by a factor of two by considering $\frac{\I+\P}{2}$ then the mixing time bound also changes by a factor of two, even for non-lazy walks. Conductance bounds don't immediately have this property, although they can be modified to behave accordingly. The following lemma shows that the $f$-congestion also ``scales properly'', in the sense that increasing the holding probability of a lazy walk also changes the $f$-congestion in an appropriate way. We also have a bound when the holding probability is decreased below $\gamma=1/2$, but at this point the walk may start to slow due to periodicity, and so our bound in this case is rather more complicated.

Let us start with the simplest case, re-scaling a walk where $\Psi(A)=\Q(A,A^c)$, such as a lazy walk.

\begin{lemma} \label{lem:rescaling}
Consider a finite, ergodic Markov chain such that $\forall x\in V:\,\P(x,x)\geq\gamma\in[0,1]$. If
$A\subset V$ with $\Psi(A)=\Q(A,A^c)$ then
$$
1-\C_f(A) = 2(1-\gamma)(1-\C_f'(A))
$$
where $\C_f'$ is the $f$-congestion for Markov kernel
$\P' = \frac{1-\gamma'}{1-\gamma}\P + \left(1-\frac{1-\gamma'}{1-\gamma}\right)\I$,
that is the Markov chain re-scaled to holding probability $1/2$.
\end{lemma}

The condition that $\Psi(A)=\Q(A,A^c)$ is necessary because, as will be seen in Example \ref{ex:complete}, a walk on the complete graph $K_m$ with holding probability $\gamma$ has
$$
1-\C_{\sqrt{a(1-a)}}(A) = \left| \frac{m\gamma-1}{m-1}\right|\,.
$$
The lemma will fail for $\gamma<1/m$, exactly the region for which $\Psi(A)\neq\Q(A,A^c)$.

\begin{proof}
Let $A_u'$ denote the evolving sets for $\P'$ and $A_u$ denote those for $\P$. 

\vspace{2ex}\noindent$\circ$ 
Suppose that $\gamma\geq 1/2$. Then
$$
A_u = 
\begin{cases}
A_{u/2(1-\gamma)}' & \textrm{if } u\leq 1-\gamma \\
A & \textrm{if } u\in[1-\gamma,\gamma] \\
A_{1-(1-u)/2(1-\gamma)}' & \textrm{if } u>\gamma
\end{cases}
$$
Applying these relations when integrating for $\C_f(A)$ leads to the lemma.

\vspace{2ex}\noindent$\circ$ 
Suppose that $\gamma<1/2$.  
Fix any $\wp_A\in[\inf\{u:\pi(A_u)\leq\pi(A)\},\sup\{u:\pi(A_u)\geq\pi(A)\}]$.

To begin with, if $y\in A$ then because $\P' = \frac{1}{2(1-\gamma)}\P + \left(1-\frac{1}{2(1-\gamma)}\right)\I$
it follows that
$$
\Q'(A,y) = \frac{1}{2(1-\gamma)}\,\Q(A,y) + \left(1-\frac{1}{2(1-\gamma)}\right)\,\pi(y)\,.
$$
Therefore, if $y\in A$ then $\Q'(A,y)\geq u\pi(y)$ if and only if
$$
\frac{\Q(A,y)}{\pi(y)} \geq \left(u -1 + \frac{1}{2(1-\gamma)}\right)\,2(1-\gamma)
  = 1 - 2(1-u)(1-\gamma)\,.
$$
It follows that $A_u'\cap A = A_{1 - 2(1-u)(1-\gamma)}\cap A$.

A similar argument holds for the case of $y\in A^c$ and shows that
$A_u'\cap A^c = A_{2u(1-\gamma)}\cap A^c$.

Combining the case of $A_u'\cap A^c$ and $A_u'\cap A$ shows that 
$$
A_u' = \left(A_{1 - 2(1-u)(1-\gamma)}\cap A\right) 
   \cup \left(A_{2u(1-\gamma)}\cap A^c\right)\,.
$$

If $u > \frac{\wp_A}{2(1-\gamma)}$ then $2u(1-\gamma)>\wp_A$ and
by Lemma \ref{lem:psi-flow} it follows that $A_u'\cap A^c=A_{2u(1-\gamma)}\cap A^c=\emptyset$
and so $A_u'=A_{1-2(1-u)(1-\gamma)}\cap A$.
If, moreover, $u> 1-\frac{1-\wp_A}{2(1-\gamma)}$ then 
$1 - 2(1-u)(1-\gamma)>\wp_A$ and so $A_{1 - 2(1-u)(1-\gamma)}\subseteq A$,
and $A_u'=A_{1-2(1-u)(1-\gamma)}$.
However, if $\frac{\wp_A}{2(1-\gamma)}<u\leq 1-\frac{1-\wp_A}{2(1-\gamma)}$ then
$1 - 2(1-u)(1-\gamma)\leq \wp_A$, and so by Lemma \ref{lem:psi-flow}
$A_u'\cap A = A_{1 - 2(1-u)(1-\gamma)}\cap A = A$ and $A_u'=A$.
Finally, if $u\leq \frac{\wp_A}{2(1-\gamma)}$ then $1-2(1-u)(1-\gamma)\leq\wp_A$ since $\gamma\leq 1/2$, so $A_u'\cap A=A$, and since $2u(1-\gamma)\leq\wp_A$ then by Lemma \ref{lem:psi-flow}
$A_{2u(1-\gamma)}\supseteq A$, which implies that $A_u'=A_{2u(1-\gamma)}$.

Putting these all together, we have that
$$
A_u' = 
\begin{cases}
A_{2u(1-\gamma)} & \textrm{if } u\leq \frac{\wp_A}{2(1-\gamma)} \\
A & \textrm{if } u\in\left(\frac{\wp_A}{2(1-\gamma)},1-\frac{1-\wp_A}{2(1-\gamma)}\right] \\
A_{1-2(1-u)(1-\gamma)} & \textrm{if } u>1-\frac{1-\wp_A}{2(1-\gamma)}
\end{cases}
$$

Applying these relations when integrating for $\C_f'(A)$ leads to the lemma.
\end{proof}

One consequence of this is a conductance lower bound on congestion which scales properly. See Lemma \ref{lem:non-lazy-conductance} for details.

We finish this section by giving a more complicated re-scaling inequality which applies even when $\Psi(A)\neq\Q(A,A^c)$.

\begin{lemma} \label{lem:scaling}
Consider a Markov chain $\P$ with holding probability $\gamma\in[0,1/2]$, and let $f:\,[0,1]\to\R_+$ be a concave function. 
%
Let $\C_f'$ be the $f$-congestion for Markov kernel
$\P' = \frac{1-\gamma'}{1-\gamma}\P + \left(1-\frac{1-\gamma'}{1-\gamma}\right)\I$,
that is the Markov chain re-scaled to holding probability $1/2$.
Suppose $h:\,[0,1]\to[0,1]$ is a decreasing function with $\int_0^1 h(w)\,dw=\pi(A)$. Then
$$
2(1-\gamma)(1-\C_f'(A)) \geq 1-\C_f(A) \geq 2\gamma\,\left(1-\frac{\int_0^1 f(h(w))\,dw}{f(\pi(A))}\right)\,,
$$
if
\begin{eqnarray}
\int_0^{\gamma/2(1-\gamma)} (h(w)-\pi(A))\,dw&=&\int_{1-\gamma/2(1-\gamma)}^1 (\pi(A)-h(w))\,dw \label{eqn:condition1} \\
\int_0^{1/2} \pi(A_w')\,dw &=& \int_0^{1/2} h(w)\,dw \nonumber \\
\forall t\in[0,1]:\,\int_0^t \pi(A_w')\,dw &\geq& \int_0^t h(w)\,dw \label{eqn:condition2}
\end{eqnarray}
\end{lemma}

If $\int_0^{\gamma}+\int_{1-\gamma}^1 \pi(A_w')\,dw=2\gamma\pi(A)$ then by setting $h(w)=\pi(A_w')$ this implies the lower bound $1-\C_f(A)\geq 2\gamma(1-\C_f'(A))$. If this holds in general then the lemma becomes the much more consise $2(1-\gamma)(1-\C_f'(A))\geq 1-\C_f(A) \geq 2\min\{\gamma,1-\gamma\}(1-\C_f'(A))$, and there is no need for the function $h$. In Example \ref{ex:conductance-K_m} it is shown that the simple random walk on $K_m$ with holding probability $\gamma\geq 1/m$ has $1-\C_{\sqrt{a(1-a)}}(A)=\frac{m}{m-1}(1-\gamma)=2(1-\gamma)(1-\C_f'(A))$, and so the upper bound is sharp for every $\gamma>0$. The lower bound is sharp for the two point space because in Example \ref{ex:example1} it was found that $1-\C_{\sqrt{a(1-a)}}(A)=2\min\{\gamma,1-\gamma\}$. More generally, it is sharp when $A$ is a bipartition for the simple random walk on the cycle with an even number of vertices and holding probability $\gamma$.

\begin{proof}

\vspace{2ex}\noindent
--First, the upper bound. 

If $\gamma\geq 1/2$ then by Lemma \ref{lem:rescaling} the upper bound is an equality, so we may assume $\gamma<1/2$. 
By the relations shown in the proof of Lemma \ref{lem:rescaling}, and because we may assume without loss that $\wp_A\in[\gamma,1-\gamma]$, then
\begin{equation} \label{eqn:decompose}
A_u = \begin{cases}
A_{u/2(1-\gamma)}' & \textrm{if } u\leq\gamma \\
(A_{u/2(1-\gamma)}'\setminus A)\cup A_{1-(1-u)/2(1-\gamma)}' & \textrm{if }u\in[\gamma,1-\gamma] \\
A_{1-(1-u)/2(1-\gamma)}' & \textrm{if } u\geq 1-\gamma
\end{cases}
\end{equation}

Integrate over $u\notin[\gamma,1-\gamma]$:
$$
\int_0^{\gamma}+\int_{1-\gamma}^1 f(\pi(A))-f(\pi(A_u))\,du
 = 2(1-\gamma)\int_0^{\frac{\gamma}{2(1-\gamma)}}+\int_{1-\frac{\gamma}{2(1-\gamma)}}^1 f(\pi(A)-f(\pi(A_w'))\,dw
$$

Now integrate over $u\in[\gamma,1-\gamma]$. To do this, recall from equation \eqref{eqn:bigup_smalldown} that if $a+b=c+d$ and $a>c>d>b$ then $f(a)+f(b)\leq f(c)+f(d)$. Well, $\pi(A_{\frac{u}{2(1-\gamma)}}')=\pi(A)+\pi(A_{\frac{u}{2(1-\gamma)}}'\setminus A)$ and so 
$$
f(\pi(A_{\frac{u}{2(1-\gamma)}}'))+f(\pi(A_{1-\frac{1-u}{2(1-\gamma)}}')) \leq f(\pi(A))+f(\pi((A_{\frac{u}{2(1-\gamma)}}'\setminus A)\cup A_{1-\frac{1-u}{2(1-\gamma)}}'))\,.
$$
It follows that
\begin{eqnarray*}
\lefteqn{\int_{\gamma}^{1-\gamma} f(\pi(A))-f(\pi(A_u))\,du } \\
   &\leq& \int_{\gamma}^{1-\gamma} 
              \left(f(\pi(A))-f(\pi(A_{\frac{u}{2(1-\gamma)}}'))\right) 
            + \left(f(\pi(A))-f(\pi(A_{1-\frac{1-u}{2(1-\gamma)}}'))\right)\,du \\
 &=& 2(1-\gamma)\int_{\frac{\gamma}{2(1-\gamma)}}^{1-\frac{\gamma}{2(1-\gamma)}} f(\pi(A))-f(\pi(A_w'))\,dw
\end{eqnarray*}

Adding the integrals for $u\notin [\gamma,1-\gamma]$ and $u\in[\gamma,1-\gamma]$ gives the upper bound.

\vspace{2ex}\noindent
--Now, the lower bound.

The conditions on $h(u)$ were chosen so that if
$$
h_2(u) = \begin{cases}
h\left(\frac{u}{2(1-\gamma)}\right) & \textrm{if } u\leq\gamma \\
h\left(\frac{\gamma}{2(1-\gamma)}+\frac{u-\gamma}{2(1-\gamma)}\right)-\pi(A) + h\left(\frac 12+\frac{u-\gamma}{2(1-\gamma)}\right) & \textrm{if } u\in[\gamma,1-\gamma] \\
h\left(1-\frac{1-u}{2(1-\gamma)}\right) & \textrm{if } u> 1-\gamma
\end{cases}
$$
then $h_2$ is decreasing, $\int_0^1 h_2(u)\,du=\pi(A)$, and $\forall t\in[0,1]:\,\int_0^t \pi(A_u)\,du \geq \int_0^t h_2(u)\,du$ by \eqref{eqn:condition2} and \eqref{eqn:decompose}. It follows from Lemma \ref{lem:worst-case} that
$$
\C_f(A) \leq \frac{\int_0^1 f(h_2(u))\,du}{f(\pi(A))}\,.
$$


By Jensen's Inequality, the relation $\int_0^1 h(u)\,du=\pi(A)$ and \eqref{eqn:condition1},
$$
\int_{\gamma}^{1-\gamma} f(h_2(u))\,\frac{du}{1-2\gamma} 
  \leq f\left(\int_{\gamma}^{1-\gamma}h_2(u)\frac{du}{1-2\gamma}\right)
  = f(\pi(A))\,,
$$
and so 
\begin{eqnarray*}
\lefteqn{ \int_0^1 f(\pi(A))- f(h_2(u))\,du } \\
  &\geq& \int_0^{\gamma}+\int_{1-\gamma}^1 f(\pi(A)) - f(h_2(u))\,du \\
  &=& 2(1-\gamma)\,\int_0^{\frac{\gamma}{2(1-\gamma)}}+\int_{1-\frac{\gamma}{2(1-\gamma)}}^1 f(\pi(A)) - f(h(u))\,du \\
  &=& 2\gamma\int_0^1 f(\pi(A))-f(h_3(u))\,du \\
 &\geq& 2\gamma\int_0^1 f(\pi(A))-f(h(w))\,dw\,.
\end{eqnarray*}
In the final equality
$$
h_3(u) = \begin{cases}
h\left(\frac{\gamma}{1-\gamma} u\right)      & \textrm{when } u\leq 1/2 \\
h\left(1-\frac{\gamma}{1-\gamma}(1-u)\right) & \textrm{when } u>1/2 
\end{cases}
$$
The final inequality is an application of Lemma \ref{lem:worst-case}. To see this, first apply \eqref{eqn:condition1} to see that $\int_0^1 h_3(u)\,du=\pi(A)=\int_0^1 h(u)\,du$. If $t\leq 1/2$ then $\int_0^t h_3(u)\,du\geq\int_0^t h(u)\,du$ because  $h(u)$ is a decreasing function and so $h_3(u)\geq h(u)$. If $t>1/2$ then $\int_0^t h_3(u)\,du=\pi(A)-\int_t^1 h_3(u)\,du \geq \pi(A)-\int_t^1 h(u)\,du=\int_0^t h(u)\,du$, again because $h(u)$ is decreasing.
\end{proof}

\section{Examples} \label{sec:examples}

The purpose of this section is to demonstrate sharpness of bounds. We start with the elementary example of a walk on a complete graph, in which each bound is either sharp or at least asymptotically of the correct order. This is followed by a careful analysis of random walk on a cycle, in which we show fairly sharp total variation mixing time bounds. We finish by discussing the simple random walk on a directed non-lazy Eulerian graph, for which our methods appear to give the first proof of a mixing time bound.

First, we see that the conductance bounds are sharp.
\begin{example} \label{ex:example1}
Consider the uniform two-point space $\{0,1\}$ with transition kernel 
$\P(0,0)=\P(1,1)=\gamma\in[0,1]$ and $\P(0,1)=\P(1,0)=(1-\gamma)$.
Then $\~\Phi(A) = 2(1-\gamma)$, and so by Lemma \ref{lem:non-lazy-conductance}
$$
2(1-\gamma) \geq 1-\C_{\sqrt{a(1-a)}} \geq 2\min\{\gamma,1-\gamma\}\,.
$$ 
Hence $1-\C_{\sqrt{a(1-a)}}=2(1-\gamma)$ if $\gamma\geq 1/2$.

More generally, $1-\C_{\sqrt{a(1-a)}}\leq\~\phi=2\min\{\gamma,1-\gamma\}$ and so the upper and lower bound are equal and  $1-\C_{\sqrt{a(1-a)}}=2\min\{\gamma,1-\gamma\}$ for all $\gamma\in[0,1]$.
\end{example}

Theorem \ref{thm:profile} can lead to sharp bounds, even for holding probability under $1/2$.

\begin{example}  \label{ex:conductance-K_m}
Consider the random walk on the complete graph $K_m$ with $\P(x,y)=1/m$. Then $\forall A\subset V:\,\~\phi(A)=1$ and so $1\geq 1-\C_{\sqrt{a(1-a)}}(A)\geq 1-\sqrt{1-1^2}=1$.  Moreover, when $\pi(A)=1/2$ then $1\geq 1-\C_{a(1-a)}(A)\geq 1$ and $1\geq 1-\C_{a\log(1/a)}(A)\geq (2\log 2)^{-1}\approx 0.72$.  Therefore at least two of the three bounds in Theorem \ref{thm:profile} can be sharp.

By Lemma \ref{lem:rescaling} a rescaling argument can be used to extend this to sharp bounds for other holding probabilities, as long as $\Psi(A)=\Q(A,A^c)$. In particular, if $\gamma\geq 1/m$ then the walk on $K_m$ with $\P(x,x)=\gamma$ and $\P(x,y)=\frac{1-\gamma}{m-1},\,\forall y\neq x$ satisfies $\Psi(A)=\pi(A)\pi(A^c)\frac{m}{m-1}(1-\gamma)=\Q(A,A^c)$. Hence, if $\gamma=1/m$, and $\P'$ is the walk with holding probability $1/2$, then $1-\C_{\sqrt{a(1-a)}}'(A)=\frac{1-\C_{\sqrt{a(1-a)}}(A)}{2(1-1/m)}=\frac{1}{2(1-1/m)}$. More generally, if $\gamma\geq 1/m$ then
$$
1-\C_{\sqrt{a(1-a)}}(A) = 2(1-\gamma)(1-\C_{\sqrt{a(1-a)}}'(A)) =\frac{m}{m-1}\,(1-\gamma)\,.
$$
\end{example}

In fact, the $f$-congestion can be used to show sharp mixing time bounds, regardless of holding probability.

\begin{example} \label{ex:complete}
Given $\alpha\in[-\frac{1}{m-1},1]$ consider the walk on $K_m$ with $\P(x,y)=(1-\alpha)/m$ for all $y\neq x$
and $\P(x,x)=\alpha+(1-\alpha)/m$, that is, choose a point uniformly at random and
move there with probability $1-\alpha$, otherwise do nothing. 

The $n$ step distribution is $\P^n(x,x)=\frac 1m+\alpha^n\left(1-\frac 1m\right)$ and 
$\P^n(x,y) = \frac 1m - \frac{\alpha^n}{m}$ for all $y\neq x$. Therefore, when $\alpha\in[0,1]$ then 
$\D(\P^n(x,\cdot)\|\pi) = (1+o_m(1))\alpha^n\log m$ as $m\rightarrow\infty$. When $\alpha\in\left[\frac{-1}{m-1},\,1\right]$
then $\|\P^n(x,\cdot)-\pi\|_{TV} = |\alpha|^n (1-1/m)$ and $\|\P^n(x,\cdot)-\pi\|_{L^2(\pi)} = |\alpha|^n \sqrt{m-1}$.

Now for evolving sets. 

If $\alpha\in[0,1]$ then
$$
\pi(A_u) =
\begin{cases}
0      & \mbox{if $u\in(\alpha + (1-\alpha)\pi(A),1]$} \\
\pi(A) & \mbox{if $u\in((1-\alpha)\pi(A), \alpha+(1-\alpha)\pi(A)]$} \\
1      & \mbox{if $u\in[0,(1-\alpha)\pi(A)]$}
\end{cases}
$$
A quick calculation shows that $\C_{a(1-a)}=\C_{a\log(1/a)}=\C_{\sqrt{a(1-a)}}=\alpha$, and so Theorem \ref{cor:main} implies $\|\P^n(x,\cdot)-\pi\|_{TV} \leq \alpha^n\,(1-1/m)$, $\D(\P^n(x,\cdot)\|\pi) \leq \alpha^n\,\log m$
and $\|\P^n(x,\cdot)-\pi\|_{L^2(\pi)} \leq \alpha^n\sqrt{m-1}$. Total variation and $L^2$ bounds are correct, while relative entropy is asymptotically correct. 

When $\alpha\in\left[\frac{-1}{m-1},\,0\right)$ then 
$$
\pi(A_u) =
\begin{cases}
0        & if\ u > (1-\alpha)\pi(A), \\
\pi(A^c) & if\ u > \alpha + (1-\alpha)\pi(A), \\
1        & otherwise
\end{cases}
$$
This time $\C_{a(1-a)}=\C_{\sqrt{a(1-a)}}=-\alpha$ and so $\|\P^n(x,\cdot)-\pi\|_{TV} \leq (-\alpha)^n\,(1-1/m)$ and $\|\P^n(x,\cdot)-\pi\|_{L^2(\pi)} \leq (-\alpha)^t\sqrt{m-1}$, both exact.
\end{example}

%

A harder walk to bound is the simple random walk on the cycle $C_m$, that is $\P(x,x\pm 1)=1/2$. A bound must distinguish between the (periodic) walk on a cycle of even length, and the (convergent) walk on a cycle of odd length. 

\begin{example} \label{ex:cycles}
The walk on a cycle $C_m$ of even length has $\~\phi = 0$ because it is bipartite, with the worst set $A$ given by choosing $m/2$ alternating points around the cycle, and $B=A$ in the definition of $\Psi(A)$. Therefore $0=\~\phi\geq 1-\C_f\geq 0$ for all of the quantities dealt with in Theorem \ref{thm:profile}. Correctly, none of our bounds show mixing.

Now for the cycle $C_m$ of odd length. If $\pi(A)<1/2$ then 
$\Psi(A)\geq 1/2m$, with the worst sets given by points alternating around the cycle, as in 
the white vertices of Figure \ref{fig:odd-cycle}. Then $\Psi(A)=\Q(A,B)$ when $B$ contains those points at 
least distance two from $A$, one point adjacent to these and $A$, and the points in $A$, 
corresponding to the circled regions in Figure \ref{fig:odd-cycle}.
\begin{figure}[ht]
\begin{center}
\includegraphics[height=1.5in]{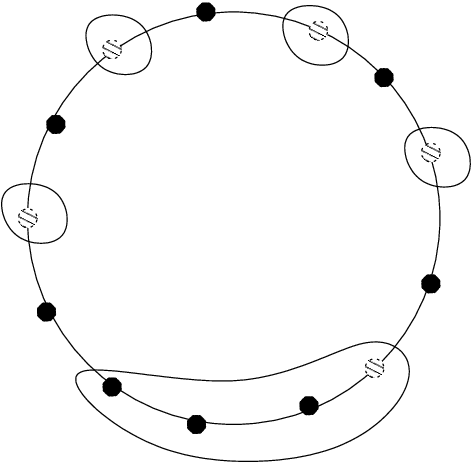}
\caption{Let $A$ be the white vertices and $B$ be the circled points. Then $\Psi(A)=\Q(A,B)=1/2m$.} 
   \label{fig:odd-cycle}
\end{center}
\end{figure}

Therefore 
$$
1-\C_{a(1-a)}(A) \geq 4\~\phi^2(A)\pi(A)\pi(A^c) \geq \frac{1}{m^2\pi(A)\pi(A^c)}
$$
By Theorem \ref{thm:profile-mixing} it follows that if $\epsilon\geq 1/2$ then
$$
\tau(\epsilon) \leq \int_{\pi_*}^{1-\epsilon} \frac{dx}{(1-x)(1-\C_{a(1-a)}(x))}
  \leq \int_{1/m}^{1-\epsilon} m^2\,x\,dx
  = \frac{m^2}{2}\,\left((1-\epsilon)^2-1/m^2\right)
$$
and so if $x\in V$ then
\begin{equation} \label{eqn:firstupper}
\|\P^n(x,\cdot)-\pi\|_{TV} \leq 1-\frac 1m\,\sqrt{1+2n}\quad \textrm{if}\ n\leq \frac{m^2}{8}-\frac 12\,.
\end{equation}
Standard techniques give poor bounds for large epsilon, such as $\epsilon>1/2$ above.

Bounds for $\epsilon<1/2$ can be obtained similarly, but better asymptotics can be derived by a slight modification of the argument. Observe that
\begin{eqnarray*}
\|\P^n(x,\cdot)-\pi\|_{TV} &\leq& \frac{1}{\pi(x)}\,\E \pi(S_n)(1-\pi(S_n)) 
   \leq \frac{1}{\pi(x)}\,\frac{\pi_*(1-\pi_*)}{\sin(3.14\,\pi_*)}\,\E \sin(3.14\,\pi(S_n)) \\
   &\leq& \frac{\sin(3.14\,\pi(x))}{\pi(x)}\,\frac{\pi_*(1-\pi_*)}{\sin(3.14\,\pi_*)}\,\C_{\sin(\pi a)}^n
   \leq (1-\pi_*)\,\C_{\sin(\pi a)}^n
\end{eqnarray*}
where $3.14$ is used to represent the number $\pi$. The choice of $\C_{\sin(\pi a)}$ is because if $\Psi(A)\geq C$ for some constant $C$ then $\C_f$ is minimized by $f(a)=\sin(\pi a)$ (see \cite{Mon07.1} for details).

Now, when $\wp_A=1/2$ then by Lemma \ref{lem:worst-case} and equation \eqref{eqn:lower_sets} it follows that $\C_{\sin(\pi a)}(A) \leq \cos(2\pi\Psi(A))$. On the cycle, if $A\subset V$ then $\pi(A_u)>\pi(A)$ when $u<1/2$, while $\pi(A_u)<\pi(A)$ when $u>1/2$, so $\wp_A=1/2$. Combined with
the earlier bound $\Psi(A)\geq 1/2m$ it follows that $\C_{\sin(\pi a)}(A) \leq \cos(\pi/m)$.
Then
\begin{equation} \label{eqn:better}
\|\P^n(x,\cdot)-\pi\|_{TV} \leq (1-\pi_*)\,\C_{\sin(\pi a)}^n = (1-1/m)\,\cos^n(\pi/m)\,.
\end{equation}

A fairly close lower bound holds as well. Let $\lambda_{max}=\max\{\lambda_2,\,|\lambda_m|\}$ be the second largest magnitude of an eigenvalue of $\P$. It is easily verified that $\cos\left(\frac{\pi(m-1)}{m}\right)$ is an eigenvalue with eigenvector $f(j)=\cos\left(\frac{2\pi(m-1)j}{m}\right)$, so
$\lambda_{max}\geq \left|\cos\left(\frac{\pi(m-1)}{m}\right)\right|=\cos(\pi/m)$.
But then
\begin{equation}\label{eqn:lower_bound}
\max_{x\in V} \|\P^n(x,\cdot)-\pi\|_{TV} \geq \frac 12\,\lambda_{max}^n \geq \frac 12\,\cos^n(\pi/m)\,.
\end{equation}
The first inequality is a general bound for time-reversible chains.


One bound that appears in the literature is 
$$
\frac 12\,\cos^n(\pi/m) \leq \max_{x\in V} \|\P^n(x,\cdot)-\pi\|_{TV} \leq e^{-\pi^2\,n/2m^2}
\quad if\ n\geq m^2/40\,.
$$
Our bound \eqref{eqn:better} is at most $(1-1/m)e^{-\pi^2\,n/2m^2}$, mildly better 
overall and with no conditions on $n$.
The old bound also required knowledge of the complete spectrum of the transition matrix. In contrast, we required only examination of edge expansion properties.
\end{example}

We finish with an example where our methods give the only known mixing time bounds, the simple random walk on a directed Eulerian graph.

\begin{example} \label{ex:eulerian}
Consider a directed Eulerian graph with vertex set $V$ and $m$ edges, that is, a strongly connected graph with in-degree=out-degree at each vertex. The simple random walk is a walk which chooses a neighboring vertex uniformly and then transitions there. This walk has $\P(x,y)=1/deg(x)$ if there is an edge from $x$ to $y$, and stationary distribution $\pi(x)=deg(x)/m$. It is known that the lazy simple random walk (i.e. $\P(x,x)=1/2$ and $\P(x,y)=1/2deg(x)$) has mixing time $\tau_2(\epsilon)=O( m^2\log(m/\epsilon))$, but nothing seems to be known about the non-lazy simple random walk even on undirected graphs.

Before stating a mixing bound we must exclude graphs on which the simple random walk does not converge. For instance, a bipartite graph. More generally, the walk is non-convergent if a directed graph has $k$ (equal sized) components such that a transition starting in component $i$ always goes to component $i+1\mod k$. The problem here is that a walk starting in one component has a neighborhood the same size as the original set, so it never grows to cover the entire space. If we let $N(A)=\{x\in V:\,\Q(A,x)>0\}$ denote the neighborhood of $A$, then the following weak expansion condition will suffice to rule out such situations:
\begin{equation} \label{eqn:expansion}
\forall A\subset V,\,\pi(A)\leq 1/2,\,\forall v\in V:\,\pi(N(A)\setminus v) \geq \pi(A)
\end{equation}
This just says that if any single vertex in the neighborhood of $A$ is removed, then the neighborhood is still at least as big as $A$. Note this cannot be satisfied if some vertex has only one outgoing edge, and so $\pi_*=\min_{v\in V}\pi(v)\geq 2/m$.

We now lower bound $\Psi(A)$. Suppose $A\subset V$ with $\pi(A)\leq 1/2$, and if $\Psi(A)=\Q(A,B)+(\pi(A^c)-\pi(B))\frac{\Q(A,v)}{\pi(v)}$. If $N(A)\subseteq B^c$ then $\pi(N(A)\setminus v)\leq \pi(B^c\setminus v)=1-\pi(B\cup v)<\pi(A)$, contradicting the expansion condition. Hence, $N(A)\cap B\neq\emptyset$ and so there are vertices $x\in A$, $y\in B$ with $\P(x,y)>0$. Then
$$
\Psi(A)\geq\Q(A,B) \geq \pi(x)\P(x,y) = \frac{deg(x)}{m}\,\frac{1}{deg(x)} = \frac{1}{m}\,.
$$
It follows that $\~\phi(r)\geq \frac{1}{m\,r(1-r)}$ if $r\leq 1/2$, and so from the convex version of equation \eqref{eqn:L2-modified-conductance} that
\begin{eqnarray*}
\tau_2(\epsilon) 
 &\leq& \left\lceil \int_{2/m}^{1/2} \frac{dr}{2r(1-r)\~\phi(r)^2/2}  + \int_{1/2}^{\frac{1}{1+\epsilon^2}} \frac{dr}{2r(1-r)\~\phi^2/2} \right\rceil \\
 &\leq& \left\lceil \frac{m^2}{12} + \frac{m^2}{8}\log\frac{1}{\epsilon}\right\rceil\,.
\end{eqnarray*}

The same argument can be used to improve on the classical $\tau_2(\epsilon)=O(m^2\log(m/\epsilon))$ bound for the lazy simple walk. Every lazy walk has $\Psi(A)=\Q(A,A^c)$, and so $\Psi(A)\geq 1/2m$ even without the expansion condition. It follows that $\~\phi(r)\geq \frac{1}{2m\,r(1-r)}$, and so the lazy simple random walk mixes in
$$
\tau_2(\epsilon) \leq \left\lceil \frac{m^2}{3}+\frac{m^2}{2}\log\frac{1}{\epsilon}\right\rceil\,.
$$

Note that the (lazy or non-lazy) simple random walk on a cycle with an odd number of vertices has $\tau_2(\epsilon)=\Theta(m^2\log\frac{1}{\epsilon})$, and so even for the lazy simple random walk our bounds are the first ones of the correct order.

A total variation bound can be found by integrating the appropriate total variation result of Theorems \ref{thm:profile-mixing} and \ref{thm:profile}. Instead, to give a taste of what improvements can be made, we note that in \cite{Mon06.6} the above technique is sharpened to show that the (non-lazy) simple random walk satisfies
$$
\tau_{TV}(\epsilon) \leq\left\lceil \frac{1}{-\log\cos\frac{2\pi}{m}}\log\frac{1-2/m}{\epsilon} \right\rceil
\approx \frac{m^2}{2\pi^2}\log\frac{1}{\epsilon}\,.
$$
This bound is exact for the simple random walk on a cycle with $3$ vertices (i.e. $K_3$ with $\alpha=-1/2$ in Example \ref{ex:complete}), while more generally equation \eqref{eqn:lower_bound} shows an extremely close lower bound for a cycle with an odd number of vertices:
$$
\tau_{TV}(\epsilon) \geq\left\lceil \frac{1}{-\log\cos\frac{2\pi}{m}}\log\frac{1}{2\epsilon} \right\rceil
\approx \frac{m^2}{2\pi^2}\log\frac{1}{2\epsilon}
$$

Numerous other improvements and generalizations are possible. See \cite{Mon06.6} in which we sharpen this analysis further, extend it to show bounds on other walks such as the max-degree walk, and also give near-optimal bounds for spectral gap and other quantities of interest.
\end{example}


\chapter{A comparison to previous isoperimetric bounds} \label{sec:congestion}

How do our new results compare to previous isoperimetric bounds? In this section we compare our new Evolving set mixing bounds to earlier Evolving Set bounds, to Spectral profile bounds, to Blocking Conductance results, and to Discrete Gradient methods.


\section{Evolving Sets}
Morris and Peres' used a more probabilistic argument than ours to show that if $x\in V$ and $S_0=\{x\}$ then
$$
\|\P^n(x,\cdot)-\pi\|_{L^2(\pi)}
   \leq \hat\E_n \frac{\min\{\sqrt{\pi(S_n)},\,\sqrt{1-\pi(S_n)}\}}{\pi(S_n)}\,,
$$
not a major difference but up to $\sqrt{2}$ times weaker than our bound in Theorem \ref{thm:main}. 
They did not have bounds on total variation or relative entropy.

Our rate of contraction $\C_{\sqrt{a(1-a)}}$ on $L^2$ distance is also better than the $\C_{\sqrt{a}}$ that they showed.
Let $f(x,y) = \sqrt{\frac{x}{y}} - \sqrt{\frac{x(1-x)}{y(1-y)}}$ with domain $x,y\in(0,1)$. This
is convex in $x$ because 
$\frac{d^2}{dx^2} f(x,y) = \frac{1-(1-x)^{3/2}\sqrt{1-y}}{4(x(1-x))^{3/2}\sqrt{y(1-y)}}\geq 0$. 
Then by Jensen's inequality,
$$
\int_0^1 \left( \sqrt{\frac{\pi(A_u)}{\pi(A)}}-\sqrt{\frac{\pi(A_u)(1-\pi(A_u))}{\pi(A)(1-\pi(A))}} \right)\,du
 \geq f\left(\int_0^1 \pi(A_u)\,du,\,\pi(A)\right)
  = f\left(\pi(A),\pi(A)\right)
  = 0,
$$
showing that $\C_{\sqrt{a(1-a)}}(A)\leq\C_{\sqrt{a}}(A)$.


\section{Spectral Profile}

Two isoperimetric bounds on mixing time are shown in the Spectral Profile paper \cite{GMT06.1}:
\begin{equation} \label{eqn:spectral-mixing}
\tau_2(\epsilon) \leq \left\lceil \int_{4\pi_*}^{4/\epsilon^2} \frac{4\,dr}{r\Phi_{\P\P^*}(r)^2}\right\rceil 
\quad\textrm{and}\quad
\tau_2(\epsilon) \leq \left\lceil \int_{4\pi_*}^{4/\epsilon^2} \frac{2\,dr}{\frac{\gamma}{1-\gamma}\,r\Phi(r)^2} \right\rceil 
\end{equation}
The holding probability $\gamma\in[0,1]$ is such that $\forall x\in V:\,\P(x,x)\geq\gamma$.

It will be shown below that the Evolving set $L^2$ bounds are at least as good as the bounds of \eqref{eqn:spectral-mixing}. However, Evolving set bounds have the advantage that they apply to other distances, such as total variation and relative entropy, for which the isoperimetric bounds on mixing via Spectral gap/profile are merely induced from the $L^2$ mixing bound. See \cite{Mon06.6} for an example where modified conductance is used to show a total variation mixing bound which is strictly better than the $L^2$ mixing bound. Comparison methods involving conductance are also available with Evolving sets, as will be shown in Section \ref{sec:evo-comparison}. On the other hand, we do not know of decomposition results which compare favorably to those available with spectral gap or log-Sobolev methods.

First, we show that bounding mixing time with modified conductance is no worse than using the multiplicative reversibilization $\P\P^*$ in \eqref{eqn:spectral-mixing}, but may give a substantial improvement. In particular, it is fairly simple to construct a distribution for $\Q(A,y)$ in the proof of Lemma \ref{lem:mod-vs-PP*} for which the upper bound is an equality, and likewise for the the lower bound, and so in the worst case scenario $\~\Phi_{\P\P^*}(r)^2=\~\phi(r)^4$, and the first bound of \eqref{eqn:spectral-mixing} may be nearly as bad as the square of the modified conductance mixing bound. 

\begin{lemma} \label{lem:mod-vs-PP*}
$$
\sqrt{\~\Phi_{\P\P^*}(A)} 
  \geq \~\phi(A) 
  \geq 1-\sqrt{1-\~\Phi_{\P\P^*}(A)} 
  \geq \frac 12\,\~\Phi_{\P\P^*}(A)\,.
$$
\end{lemma}

\begin{proof}
To simplify notation, in the definition of $\Psi(A)$ assume that the set $B$ satisfies $\pi(B)=\pi(A^c)$, i.e. $\Psi(A)=\Q(A,B)$. For the general case it suffices to split an appropriate vertex $v$, as in the proof of Lemma \ref{lem:non-lazy}.

To begin with, we need a few identities: 
$$
\Q(A,B)=\pi(A)-\Q(A,B^c) = \pi(A)-(\pi(B^c)-\Q(A^c,B^c)) = \Q(A^c,B^c)
$$
$$
\Q_{\P\P^*}(A,A^c) = \sum_{y\in V} \Q(A,y)\P^*(y,A^c)
  = \sum_{y\in V} \frac{\Q(A,y)}{\pi(y)}\left(1-\frac{\Q(A,y)}{\pi(y)}\right)\pi(y)
$$

First we bound the terms in the summation for $\Q_{\P\P^*}(A,A^c)$ for which $y\in B$.

Let $f(a)=a(1-a)$. The proof can be completed using Lemma \ref{lem:worst-case}, but using Jensen's Inequality is much simpler.
\begin{eqnarray*}
\sum_{y\in B} \frac{\Q(A,y)}{\pi(y)}\left(1-\frac{\Q(A,y)}{\pi(y)}\right)\frac{\pi(y)}{\pi(B)}
 &=& \sum_{y\in B} f\left(\frac{\Q(A,y)}{\pi(y)}\right)\,\frac{\pi(y)}{\pi(B)} \\
 &\leq& f\left(\sum_{y\in B} \frac{\Q(A,y)}{\pi(y)}\,\frac{\pi(y)}{\pi(B)}\right) \\
 &=& f\left(\frac{\Q(A,B)}{\pi(B)}\right) = \frac{\Psi(A)}{\pi(A^c)}\left(1-\frac{\Psi(A)}{\pi(A^c)}\right)
\end{eqnarray*}

To show a lower bound, note that if $y\in B$ and $v\notin B$ then $\frac{\Q(A,y)}{\pi(y)} \leq\frac{\Q(A,v)}{\pi(v)}$ by definition of set $B$, and so $\frac{\Q(A,y)}{\pi(y)} \leq \frac{\Q(A,B^c)}{\pi(B^c)}=\frac{\pi(A)-\Psi(A)}{\pi(A)}=1-\frac{\Psi(A)}{\pi(A)}$. Then 
\begin{eqnarray*}
\sum_{y\in B} \frac{\Q(A,y)}{\pi(y)}\left(1-\frac{\Q(A,y)}{\pi(y)}\right)\frac{\pi(y)}{\pi(B)}
 &\geq& \sum_{y\in B} \frac{\Q(A,y)}{\pi(y)}\,\frac{\Psi(A)}{\pi(A)}\,\frac{\pi(y)}{\pi(B)} \\
 &=&    \frac{\Psi(A)}{\pi(A)}\,\frac{\Q(A,B)}{\pi(B)} = \frac{\Psi(A)^2}{\pi(A)\pi(A^c)}
\end{eqnarray*}

To bound the terms over $B^c$ use the relation $\frac{\Q(A,y)}{\pi(y)}=1-\frac{\Q(A^c,y)}{\pi(y)}$ to re-write the sum:
$$
\sum_{y\in B^c} \frac{\Q(A,y)}{\pi(y)}\left(1-\frac{\Q(A,y)}{\pi(y)}\right)\frac{\pi(y)}{\pi(B^c)}
 = \sum_{y\in B^c} \frac{\Q(A^c,y)}{\pi(y)}\left(1-\frac{\Q(A^c,y)}{\pi(y)}\right)\frac{\pi(y)}{\pi(B^c)}
$$
Then follow the same steps as before, and apply the relation $\Q(A^c,B^c)=\Psi(A)$, to obtain
$$
\Psi(A) \left(1-\frac{\Psi(A)}{\pi(A)}\right)
\geq \sum_{y\in B^c} \frac{\Q(A^c,y)}{\pi(y)}\left(1-\frac{\Q(A^c,y)}{\pi(y)}\right)\pi(y) 
\geq \frac{\Psi(A)^2}{\pi(A^c)}\,.
$$

Adding the upper and lower bounds for the sums over $B$ and $B^c$ gives bounds on $\Q_{\P\P^*}(A,A^c)$:
$$
\Psi(A)(2-\~\phi(A)) \geq \Q_{\P\P^*}(A,A^c) \geq \frac{\Psi(A)^2}{\pi(A)\pi(A^c)}
$$
Dividing through by $\pi(A)\pi(A^c)$ and then re-arranging the inequalities completes the proof.
\end{proof}

The lemma induces mixing bounds in terms of $\~\Phi_{\P\P^*}(r)$ for total variation, relative entropy and $L^2$ distance. For instance,
\begin{equation} \label{eqn:PP*-mixing}
\tau_2(\epsilon) \leq \left\lceil \int_{4\pi_*}^{4/\epsilon^2} \frac{8\,dr}{r\~\Phi_{\P\P^*}(r)^2}\right\rceil 
\quad\textrm{and}\quad
\tau_2(\epsilon) \leq \left\lceil \frac{8}{\~\Phi_{\P\P^*}^2}\,\log\frac{1}{\epsilon\sqrt{\pi_*}} \right\rceil 
\end{equation}
This is not directly comparable to the Spectral profile bound, but it is never more than a factor two worse, and is strictly better when $x\~\Phi_{\P\P^*}\left(\frac{1}{1+x^2}\right)$ is convex as is often the case.

In a survey with Tetali \cite{MT06.1} we use a more specialized method based on an idea of \cite{Mor06.2}, which is applicable only to $1-\C_{\sqrt{a(1-a)}}(A)$, to show that
$$
1-\C_{\sqrt{a(1-a)}}(A) \geq 1-\sqrt[4]{1-\Phi_{\P\P^*}(A)^2}
  \geq \frac 14\,\Phi_{\P\P^*}(A)^2\,.
$$
This gives exactly the same mixing bound as the Spectral Profile result, and can be improved by a factor two when $x\Phi_{\P\P^*}^2\left(\frac{1}{1+x^2}\right)$ is convex.

Now, consider the second bound of \eqref{eqn:spectral-mixing}, with a holding probability. Modified conductance via Theorem \ref{thm:profile}, combined with Lemma \ref{lem:non-lazy}, gives a weak bound on $f$-congestion in terms of conductance for non-lazy walks. We now give a more direct argument improving substantially on this.

\begin{lemma} \label{lem:non-lazy-conductance}
Consider a Markov chain with holding probability $\gamma\in[0,1]$. If $A\subset V$ then
$$
\begin{array}{rclcl}
\vspace{1ex}
\~\Phi(A) &\geq& 1-\C_{a(1-a)}(A) 
  &\geq& \displaystyle  \left(\frac{2}{\max\{\gamma,1-\gamma\}}\right)\,\frac{\gamma}{1-\gamma}\,\~\Phi(A)^2\pi(A)\pi(A^c) \\
\vspace{1ex}
\~\Phi(A) &\geq& 1-\C_{a\log(1/a)}(A) &\geq& \displaystyle \left(\frac{1}{\max\{\gamma,1-\gamma\}}\right)\,\frac{\gamma}{1-\gamma}\,\frac{\Phi(A)^2}{\log(1/\pi(A))} \\
\~\Phi(A) &\geq& 1-\C_{\sqrt{a(1-a)}}(A) &\geq& \displaystyle \left(\frac{1}{4\max\{\gamma,1-\gamma\}}\right)\,\frac{\gamma}{1-\gamma}\,\~\Phi(A)^2
\end{array}
$$
\end{lemma}

\begin{proof}
The upper bounds follow from Theorem \ref{thm:profile} and the relation $\~\phi(A)\leq\~\Phi(A)$.

The lower bounds will be shown by using Lemma \ref{lem:scaling}. The lazy walk $\P'$ has $\~\phi_{\P'}(A)=\~\Phi_{\P'}(A)=\frac{1}{2(1-\gamma)}\~\Phi(A)$. 

The lower bounds with $\gamma>1/2$ follow immediately from Lemma \ref{lem:scaling} and Theorem \ref{thm:profile}.

Now, the lower bound for $\gamma<1/2$. By Theorem \ref{thm:conductance-congestion},
$$
\frac{\Q(A,A^c)}{2(1-\gamma)}=\Q_{\P'}(A,A^c)=\int_0^{1/2}(\pi(A_w')-\pi(A))\,dw=\int_{1/2}^1 (\pi(A)-\pi(A_w'))\,dw\,.
$$
Let $h(w)=\pi(A)+\frac{\Q(A,A^c)}{1-\gamma}$ if $w\leq 1/2$, and $h(w)=\pi(A)-\frac{\Q(A,A^c)}{1-\gamma}$ if $w>1/2$ . This satisfies the conditions of Lemma \ref{lem:scaling}. Theorem \ref{thm:profile} then completes the proof, for instance if $f(a)=\sqrt{a(1-a)}$ then
$$
1-\C_f(A) 
  \geq 2\gamma\,\left(1-\frac{\int_0^1 f(h(w))\,dw}{f(\pi(A))}\right)
  \geq 2\gamma\,\frac{\~\Phi_{\P'}(A)^2}{2}
  = \frac{\gamma}{4(1-\gamma)^2}\,\~\Phi(A)^2\,.
$$

\end{proof}

The lemma induces mixing bounds in terms of $\~\Phi(r)$ for total variation, relative entropy and $L^2$ distance. For instance,
\begin{equation} \label{eqn:holding-mixing}
\tau_2(\epsilon) \leq \left\lceil \int_{4\pi_*}^{4/\epsilon^2}
  \frac{4\max\{\gamma,1-\gamma\}}{\frac{\gamma}{1-\gamma}\,r\~\Phi(r)^2}\,dr
  \right\rceil 
\quad\textrm{and}\quad
\tau_2(\epsilon) \leq \left\lceil \frac{4\max\{\gamma,1-\gamma\}}{\frac{\gamma}{1-\gamma}\,\~\Phi^2}\,\log\frac{1}{\epsilon\sqrt{\pi_*}} \right\rceil\,. 
\end{equation}
This is not directly comparable to the Spectral profile bound, but it is never more than a factor two worse, and is strictly better when the walk is lazy (i.e. $\gamma=1/2$) or $x\~\Phi^2\left(\frac{1}{1+x^2}\right)$ is convex. 


\section{Blocking Conductance} \label{sec:blocking}

As discussed in the introduction, our methods give new insight into the mixing time bounds of Blocking conductance \cite{KLM06.1}. We note that the improved Average Conductance result of Fountoulakis and Reed \cite{FR07.1} is actually a special case of the Blocking Conductance total variation bound discussed below, so our discussion applies to their results as well.

In this section we work with the quantity $\Psi(A,t)$, first defined in equation \eqref{eqn:psiA}. That definition was only used for $t\leq\pi(A^c)$, but when $t>\pi(A^c)$ a different definition is more appropriate. In general, let
\begin{equation} \label{eqn:psiA2}
\Psi(A,t) = \min_{\substack{B\subset V,\,v\in V,\\\pi(B)\leq t,\,\pi(B\cup v)> t}}
              \Q(A,B)+(t-\pi(B))\,\frac{\Q(A,v)}{\pi(v)}
\end{equation}
if $t\leq\pi(A^c)$ and $\Psi(A,t)=\Psi(A^c,1-t)$ when $t>\pi(A^c)$. 

The Blocking Conductance theorem is the following:

\begin{theorem} \label{thm:blocking} [Blocking Conductance \cite{KLM06.1}]
Given a lazy, reversible, ergodic Markov chain then
$$
\tau_{TV}(\epsilon) \leq 15000\,\left(\int_{\pi_*}^{1/2} \hat{h}(x)\,dx+\hat{h}(1/2)\right)\,\log_2(1/2\epsilon)
$$
where $\hat{h}(x)$ can be any of the following:
\begin{enumerate}
\item 
$ \displaystyle
\forall x\in[0,1]:\,\hat{h}_{gl}(x)\geq \sup_{\substack{A\subset\hat{V},\\ \pi(A)\leq x}} \frac{1}{\pi(A)\,\psi_{gl}(A)} \quad where \quad \psi_{gl}(A)=\int_0^1 \frac{\Psi(A^c,t)}{\pi(A)^2}\,dt 
$
\item
$\displaystyle
\forall x\in[0,1]:\,\hat{h}_{mod}(x)\geq \sup_{\substack{A\subset\hat{V},\\ \pi(A)\leq x}} \frac{1}{x\,\psi_{mod}(A)} \quad where \quad 
  \psi_{mod}(A)=\int_0^1 \frac{\Psi(A^c,t)}{t\,\pi(A)}\,dt
$
\item
$\displaystyle
\forall x\in[0,1]:\,\hat{h}^+(x) \geq \sup_{\substack{A\subset\hat{V},\\ x/2\leq\pi(A)\leq x}} \frac{1}{x\,\psi^+(A)} \quad where \quad \psi^+(A)=\int_0^{\pi(A)} \frac{\Psi(A^c,t)}{\pi(A)^2}\,dt
$
\end{enumerate}
The state space $\hat{V}=[0,1]$ is the continuization of $V$, and is defined by associating to each $v\in V$ a disjoint interval of size $\pi(v)$, with ergodic flow such that if $dx\subset v_1$ and $dy\subset v_2$ then $\Q(dx,dy)=dx\,\P(v_1,v_2)\,\frac{dy}{\pi(v_2)}$.
\end{theorem}

The large coefficient is due to a conversion from one measure of mixing time to another, and the need for the continuization is because the theorem is proven in the continuous space setting. A discrete version of this is discussed in Section \ref{sec:blocking-comparison}.


To relate this to Evolving sets we first rewrite $\Psi(A,t)$ in terms of Evolving Sets \cite{Mon05.1}:

\begin{lemma} \label{lem:equivalence}
Given a finite irreducible Markov kernel and $A\subset V$ then
$$
\Psi(A^c,t) = 
\begin{cases}
\vspace{1ex}
\displaystyle \int_{w(t)}^1 (t-\pi(A_u))\,du & if\ t\leq \pi(A) \\
\displaystyle \int_0^{w(t)} (\pi(A_u)-t)\,du & if\ t\geq \pi(A) 
\end{cases}
$$
where $w(t)$ is any value satisfying $\inf\{y\,:\,\pi(A_y)\leq t\} \leq w(t) \leq \sup\{y\,:\,\pi(A_y)\geq t\}$.
\end{lemma}

\begin{proof}
We consider only the case that $t\leq\pi(A)$ since the case when $t>\pi(A)$ is similar.

By definition, if $v_1\in A_u$ and $v_2\notin A_u$ then $\frac{\Q(A,v_1)}{\pi(v_1)}>\frac{\Q(A,v_2)}{\pi(v_2)}$, and equivalently $\frac{\Q(A^c,v_1)}{\pi(v_1)}<\frac{\Q(A^c,v_2)}{\pi(v_2)}$. Hence, if $\pi(A_{w(t)})=t$ then $B=A_{w(t)}$ is the same set where the minimum occurs in the definition of $\Psi(A^c,t)$. If instead $\pi(A_{w(t)})>t$, then $B=\cup_{u>w(t)} A_u$ is the set where the minimum occurs in the definition of $\Psi(A^c,t)$, and if $v$ is any vertex in $A_{w(t)}\setminus B$ then $\Q(A,v)/\pi(v)=w(t)$. In both cases
$$
\Psi(A^c,t) = \Q(A^c,B) + (t-\pi(B))\frac{\Q(A^c,v)}{\pi(v)}\,.
$$

Let $B$ be as defined in the previous paragraph. Then, $B\subseteq A_{w(t)}$, and $A_u\subseteq B$ whenever $u>w(t)$, and so
\begin{eqnarray*}
\int_{w(t)}^1 (t-\pi(A_u))\,du &=& t(1-w(t)) - \sum_{y\in B} \left(\frac{\Q(A,y)}{\pi(y)}-w(t)\right)\pi(y) \\
  &=& t(1-w(t)) - (\Q(A,B)-w(t)\pi(B)) \\
  &=& (1-w(t))(t-\pi(B)) + \Q(A^c,B) \\
  &=& \Q(A^c,B) + (t-\pi(B))\,\frac{\Q(A^c,v)}{\pi(v)} = \Psi(A^c,t)\,.
\end{eqnarray*}
The first equality is because $\int_{x}^1 \pi(A_u)\,du=\sum_{y\in A_x} (Prob(y\in A_u)-x)\pi(y)$. The third equality uses $\Q(A,B)=\pi(B)-\Q(A^c,B)$. The fourth equality is because $\frac{\Q(A^c,v)}{\pi(v)} = 1-\frac{\Q(A,v)}{\pi(v)}=1-w(t)$ by our choice of $v$ and $w(t)$.
\end{proof}

The next step is to rewrite $f$-congestion quantities in terms of the $\psi(A)$ quantities appearing in the Blocking Conductance theorem.

\begin{lemma} \label{lem:blocking}
Let $\~\psi_{gl}(A)=\int_0^1 \frac{\Psi(A^c,t)}{\pi(A)^2\pi(A^c)^2}\,dt$ and $\~\psi^+(A)=\int_0^{\pi(A)} \frac{\Psi(A^c,t)}{\pi(A)^2\pi(A^c)^2}\,dt$. Then,
\begin{eqnarray*}
1-\C_{a(1-a)}(A) 
  &=& 2\pi(A)\pi(A^c)\~\psi_{gl}(A) \\
1-\C_{a\log(1/a)}(A) 
  &=& \frac{\psi_{mod}(A)}{\log \frac{1}{\pi(A)}} \\
1-\C_{\sqrt{a(1-a)}}(A) 
  &\geq& \frac 14\,\~\psi^+(A)\,. 
\end{eqnarray*}
The identity $1-\C_{\sqrt{a(1-a)}}(A)=1-\C_{\sqrt{a(1-a)}}(A^c)$ can be used when $\pi(A)>1/2$.
\end{lemma}

\begin{proof}
Start with the case of $1-\C_{a\log(1/a)}(A)$.

First, rewrite things a bit.
\begin{eqnarray*}
1-\C_{a\log(1/a)}(A) &=& \frac{\pi(A)\log(1/\pi(A))-\int_0^1 \pi(A_u)\log(1/\pi(A_u))\,du}{\pi(A)\log(1/\pi(A))} \\
  &=& \frac{1}{\log \frac{1}{\pi(A)}}\,\int_0^1 \int_{\pi(A_u)}^{\pi(A)} \frac{t-\pi(A_u)}{t\,\pi(A)}\,dt\,du
\end{eqnarray*}
The second equality applied the identity $\int_0^1 \pi(A_u)\,du = \pi(A)$.

Now to rewrite $\psi_{mod}(A)$ in terms of Evolving Sets. By Lemma \ref{lem:equivalence},
\begin{eqnarray} 
\int_0^1 \frac{\Psi(A^c,t)}{t\pi(A)}\,dt
  &=& \int_0^{\pi(A)} \int_{w(t)}^1 \frac{t-\pi(A_u)}{t\pi(A)}\,du\,dt 
      + \int_{\pi(A)}^1 \int_0^{w(t)} \frac{\pi(A_u)-t}{t\pi(A)}\,du\,dt \nonumber \\
  &=& \int_0^1 \int_{\pi(A_u)}^{\pi(A)} \frac{t-\pi(A_u)}{t\pi(A)}\,dt\,du \label{eqn:psi-mod}
\end{eqnarray}

The $1-\C_{a(1-a)}(A)$ result is shown similarly.

For the $1-\C_{\sqrt{a(1-a)}}(A)$ case we first re-write $\~\psi(A)^+$ in terms of Evolving Sets. Let $\wp_A=w(\pi(A))$, as in previous sections. Then
\begin{eqnarray} 
\int_0^{\pi(A)} \frac{\Psi(A^c,t)}{\pi(A)^2\pi(A^c)^2}\,dt
  &=& \int_0^{\pi(A)} \int_{w(t)}^1 \frac{t-\pi(A_u)}{\pi(A)^2\pi(A^c)^2}\,du\,dt \nonumber \\
  &=& \int_{\wp_A}^1 \int_{\pi(A_u)}^{\pi(A)} \frac{t-\pi(A_u)}{\pi(A)^2\pi(A^c)^2}\,dt\,du \nonumber \\
  &=& \frac 12\,\int_{\wp_A}^1 \frac{(\pi(A)-\pi(A_u))^2}{\pi(A)^2\pi(A^c)^2}\,du  \label{eqn:chi2-mod}
\end{eqnarray}
Finish with the inequality that for $x,y\in[0,1]$
$$
\sqrt{y(1-y)} \leq \sqrt{x(1-x)} + \frac{1-2x}{2\sqrt{x(1-x)}}\,(y-x) 
    - \frac{(y-x)^2}{2\,[x(1-x)]^{3/2}}\,\left(\delta_{y\leq x}\,\frac 14 + \delta_{y>x}\,x(1-x)\right)\,,
$$
substitute $y=\pi(A_u)$ and $x=\pi(A)$, and then integrate to obtain $\C_{\sqrt{a(1-a)}}(A)\leq 1-\frac 14\,\~\psi^+(A)$.
The relation follows from the inequality $\sqrt{z}\leq 1+\frac 12(z-1)-\frac 18(z-1)^2\delta_{z\leq 1}$ with $z=\frac{y(1-y)}{x(1-x)}$.
\end{proof}

When combined with Lemma \ref{lem:discreteCongestion} it follows, for instance, that
$$
\hat\E_{n+1} \log \frac{1}{\pi(S_{n+1})} - \hat\E_n\log \frac{1}{\pi(S_n)} = -\hat\E_n\psi_{mod}(S_n)\,.
$$
Hence the expectation of $\psi_{mod}(A)$ is exactly the rate at which the evolving set bound on relative entropy decreases. This shows that in a sense Blocking Conductance and Evolving Set bounds are both based on measuring the derivative of the distance with respect to time. Not surprisingly, the Evolving set mixing bounds then imply bounds of the Blocking Conductance form.

\begin{corollary} \label{cor:evolving-blocking}
Consider a finite (non-lazy, non-reversible) ergodic Markov chain. Then
\begin{eqnarray*}
\tau_{TV}(\epsilon) 
  &\leq& \left\lceil \frac 12\,h_{gl}(1/2) \log\frac{1-\pi_*}{\epsilon} \right\rceil \\
\tau_D(\epsilon) 
  &\leq& \left\lceil 2C\,\int_{\sqrt{\pi_*}}^{1/2} h_{mod}(x)\,dx + C\,h_{mod}(1/2)\,\log\frac{2}{\epsilon}\right\rceil \\
\tau_2(\epsilon) 
  &\leq& \left\lceil 4\,\int_{4\pi_*}^{1/2} h^+(x)\,dx + h^+(1/2)\log\frac{2\sqrt{2}}{\epsilon} \right\rceil
\end{eqnarray*}
where
$$
h_{gl}(x) = \max_{\substack{A\subset V,\\ \pi(A)\leq x}}\frac{1}{\pi(A)\psi_{gl}(A)},
\quad
h_{mod}(x) = \max_{\substack{A\subset V,\\ \pi(A)\leq x}}\frac{1}{x\psi_{mod}(A)},
\quad
h^+(x) = \max_{\substack{A\subset V,\\ \pi(A)\leq x}}\frac{1}{x\psi^+(A)}
$$
and $C$ is the optimal constant satisfying
$$
\forall r\geq \pi_*:\,\min_{\pi(A)\leq r}\frac{\psi_{mod}(A)}{\log(1/\pi(A))} \geq C^{-1} \min_{\pi(A)\leq r}\frac{\psi_{mod}(A)}{\log(1/r)}\,.
$$
\end{corollary}

\begin{proof}
For the total variation and $L^2$ bounds apply Corollary \ref{cor:main} and Theorem \ref{thm:profile-mixing} respectively to obtain mixing time bounds in terms of various $1-\C_f(A)$. Replacing the $f$-congestion by the appropriate $\psi(A)$ quantities from Lemma \ref{lem:blocking} then gives the results. However, the relative entropy case requires more work. This is because  $\forall r>1/2$ both $1-\C_{a(1-a)}(r)=1-\C_{a(1-a)}(1/2)$ and $1-\C_{\sqrt{a(1-a)}}(r)=1-\C_{\sqrt{a(1-a)}}(1/2)$, while $1-\C_{a\log(1/a)}(r)\neq 1-\C_{a\log(1/a)}(1/2)$ when $r>1/2$.

From Theorem \ref{thm:main}, it follows that if $g(a)=\min\{1+\log \frac{1}{2a},\,\frac{1-a}{a}(1+\log\frac{1}{2(1-a)})\}$ then
$$
\D(\P^n(x,\cdot)\|\pi) \leq \hat\E_n \log \frac{1}{\pi(S_n)}
 \leq \hat\E_n g(\pi(S_n))\,.
$$
By Theorem \ref{thm:profile-mixing}, and the relation $\C_{ag(a)}(r)=\C_{ag(a)}(1/2)$ for $r\geq 1/2$ (since $ag(a)=(1-a)g(1-a)$), the mixing time is then bounded by
$$
\tau_D(\epsilon) \leq
  \left\lceil \int_{\sqrt{e\pi_*/2}}^{1/2} \frac{2\,dr}{r(1+\log(1/2r))(1-\C_{ag(a)}(r))}
    + \frac{2\log\frac{2}{\epsilon}}{1-\C_{ag(a)}(1/2)}
  \right\rceil\,.
$$

Consider set $A\subset V$ with $\pi(A)\leq r\leq 1/2$. Then $\pi(A)g(\pi(A))=\pi(A)(1+\log\frac{1}{2\pi(A)})$ and $ag(a)\leq a(1+\log\frac{1}{2a})\,\forall a\in[0,1]$, and so
\begin{eqnarray*}
1-\C_{ag(a)}(A) &\geq& 1-\C_{a(1+\log(1/2a)}(A) \\
  &=&  \frac{\log(1/\pi(A))}{1+\log(1/2\pi(A))}\,(1-\C_{a\log(1/a)}(A)) \\
  &\geq& \frac{\log(1/r)}{1+\log(1/2r)}\,C^{-1}\,\frac{\psi_{mod}(r)}{\log(1/r)}\,.
\end{eqnarray*}

Substituting this into the bound on $\tau_D(\epsilon)$ given above completes the proof.
\end{proof}

The Corollary shows that as long as the bottlenecks get sufficiently worse as set size increases, then Evolving set bounds sharply improve on Blocking conductance results. To see this for the $h_{mod}$ case note that if $\min_{\pi(A)\leq r}\frac{\psi_{mod}(A)}{\log(1/\pi(A))} = \min_{\pi(A)\in[r/2,r]}\frac{\psi_{mod}(A)}{\log(1/\pi(A))}$, i.e. congestion decreases with set size, then it suffices to take $C=2$. The new bounds are, however, usually much better than the Blocking Conductance bounds because the laziness and reversibility requirements have been dropped, the bounds are given in terms of stronger measures of distance, and there is no need to work in a continuous state space.

This suggests that the Blocking Conductance method is the best method only in the case when bottlenecks are worst at small sets. This situation can arise when total variation mixing time is an order of magnitude faster than relative-entropy or $L^2$-mixing times. For instance, Fountoulakis and Reed \cite{FR07.1} use a version of Blocking Conductance to find the correct mixing time of walks on the giant component of the random graph $G_{n,p}$. The bottleneck condition also does not apply to certain walks used for estimating volume of convex bodies, or to Example \ref{ex:bottleneck} below. 

The interested reader can use the quantities calculated in Example \ref{ex:complete} to find that Corollary \ref{cor:evolving-blocking} is within a factor $4$ of being sharp for the walk on a complete graph. A ``convex'' version, based on Theorem \ref{thm:profile-mixing}, can be used to strengthen this to a factor $2$.

\begin{remark}
A straightforward generalization of work in \cite{Mon05.1} can be used to show that
$$
\psi_{gl}(A) \geq \frac 12\psi_{mod}(A) \geq 1-\C_{\sqrt{a}}(A) \geq \frac 14\psi^+(A) \geq \frac 14\phi(A)^2\,.
$$
Hence, these various $\psi(A)$ quantities are closely related to each other, and to modified conductance.
\end{remark}

\begin{remark} \label{rmk:blocking-vertex}
For a lazy walk a useful interpretation of $\psi^+(A)$ is given in \cite{KLM06.1}:
$$
\psi^+(A) 
 \geq \sup_{\lambda\leq\pi(A)} \min_{\substack{S\subset A,\\ \pi(S)<\lambda}} \frac{\lambda\Q(A\setminus S,A^c)}{\pi(A)^2} 
 \geq \frac 12\Phi^2(A)\,.
$$
When combined with Lemma \ref{lem:blocking} it follows that
$$
1-\C_{\sqrt{a(1-a)}}(A) \geq \frac 14 \sup_{\lambda\leq\pi(A)} \min_{\substack{S\subset A,\\ \pi(S)<\lambda}} \frac{\lambda\Q(A\setminus S,A^c)}{\pi(A)^2\pi(A^c)^2} 
 \geq \frac 18\~\Phi^2(A)\,.
$$

This can be interpreted as follows. Let $\lambda$ denote the maximal size of a ``blocking set'', such that if any set $S$ smaller than this is blocked from transitioning then it does not block too much of the ergodic flow $\Q(A,A^c)$. For instance,
$\Q(A\setminus S,A^c)=\Q(A,A^c)-\Q(S,A^c)\geq \Q(A,A^c)-\lambda/2$,
and so by setting $\lambda=\Q(A,A^c)$ then the first lower bound
on $\psi^+(A)$ implies the second.

See Remark \ref{rmk:spread-vertex} for a similar discussion involving the non-lazy case.
\end{remark}

\begin{example} \label{ex:bottleneck}
The $L^2$ mixing time can be slower than total variation mixing time when there is a bottleneck at a small set, in which case the difference between Theorem \ref{thm:blocking} and Corollary \ref{cor:evolving-blocking} may be real, and not simply an artifact of the method of proof. 

Consider the complete graph $K_m$ on $m$ vertices and attach an additional vertex $v$ by a single edge. We examine the lazy max-degree walk given by choosing a neighboring vertex with probability $1/2m$ each, and otherwise do nothing.

First, bound $\psi^+(A)$. If $A=\{v\}$ then let $\lambda=\pi(\{v\})=\frac{1}{m+1}$. The only set $\pi(S)<\lambda$ is $S=\emptyset$, and so $\psi^+(\{v\})\geq \frac{\lambda\,\Q(\{v\},K_m)}{\pi(\{v\})^2} =
\frac{1}{2m}$. If $A\neq\{v\}$ then $\Phi(A)\geq 1/8$, and so $\psi^+(A)\geq\frac{1}{128}$.

To bound mixing via Blocking Conductance, note that $\hat{h}^+(r)\leq \frac{2m}{x}$ if $r\leq \frac{1}{2(m+1)}$, while $\hat{h}^+(r)\leq \frac{128}{x}$ otherwise. Then, by Theorem \ref{thm:blocking},
$$
\tau_{TV}(\epsilon) = O(m\log(1/\epsilon))\,,
$$
which is of the correct order. 

For Evolving Sets, we can only say that $h^+(r)\leq \frac{2m}{x}$ for all $r$. Then, by Corollary \ref{cor:evolving-blocking},
$$
\tau_2(\epsilon)=O(m\log(m/\epsilon))
$$
which is again of the correct order.
\end{example}

\section{Comparison and Canonical paths for Evolving Sets} \label{sec:evo-comparison}

One of the most useful isoperimetric arguments for bounding mixing times has been the method of Canonical Paths, and in its more general form the method of Comparison. In this section we find similar results for Evolving Sets, although our results are somewhat weaker than might be hoped as they only allow us to compare the $f$-congestion $\C_f$ of one Markov chain with the Conductance Profile $\~\Phi(r)$ of another. Nevertheless, the results do serve to provide an overarching framework explaining why different versions of canonical path bounds are stronger in different situations, and in particular why $\rho_e$ (edge-congestion) can be multiplied by any of $\rho_e$ (edge-congestion), or $\rho_v$ (vertex-congestion) or $\ell$ (maximum path length) in order to bound spectral gap (and mixing time): 
$\lambda \geq c/\rho_e\rho_e,\,c/\rho_e\rho_v,\,c/\rho_e\ell$. 
Certainly $\rho_e$ is a measure of edge-congestion, but it is perhaps more surprising that the 
terms potentially multiplied by this are all measures of vertex-congestion; for a lazy chain $\rho_v\geq\rho_e/2$
and so $\rho_e$ bounds vertex congestion, certainly $\rho_v$ does as well,
and we will find that $\ell$ is a guarantor of good ``typical'' vertex congestion in some well defined sense.


In \cite{Mon07.1} we found that if only the worst case edge expansion (i.e., conductance or edge congestion) 
are known then the best evolving set bound that can be obtained for spectral gap is via $1-\C_{\sqrt{a(1-a)}}$, which we have seen to be a natural quantity for bounding $L^2$ or chi-square distance. 
The same idea suggests this to be the best that can be done for mixing time as well. This partially explains 
why all path bounds have been bounds on chi-square distance, since paths are generally used to show good
edge expansion. 


We start with an elementary comparison result, comparison of ergodic flows, to
illustrate the use of canonical paths.

\begin{theorem}[Comparison of Flows] \label{thm:flowcomparison}
Suppose $\MM$ and $\MM'$ are finite Markov chains on state space $V$
and edge sets $E$ and $E'$ respectively. To every edge $e'=(x,y)\in E'$ associate a path $\gamma_{xy}\subset E$ from $x$ to $y$. Let
$$
\rho_e = \max_{e=(u,v)\in E} \frac{1}{\pi(u)\P(u,v)}\,\sum_{\gamma_{xy} \ni e} \pi'(x)\,\P'(x,y)\,.
$$
Then, for every $A\subset V$,
$$
\Q(A,A^c) \geq \Q'(A,A^c)/\rho_e\ .
$$
\end{theorem}

\begin{proof}
For each edge $e'=(x,y)\in E'$ with $x\in A$ and $y\in A^c$, transport flow of $\pi'(x)\,\P'(x,y)$ 
along the path $\gamma_{xy}\subset E$ from $x$ to $y$, for a total of $\Q'(A,A^c)$
from $A$ to $A^c$. By definition of $\rho_e$ it follows that
if $e=(u,v)\in E$ with $u\in A$, $v\in A^c$ then
$\pi(u)\P(u,v) \geq \rho_e^{-1}\,\sum_{\gamma_{xy}\ni e} \pi'(x)\,\P'(x,y)$, and therefore
$$
\Q(A,A^c) = \sum_{\substack{e=(u,v)\in E,\\u\in A,\,v\in A^c}} \pi(u)\P(u,v) 
     \geq \frac{\sum_{x\in A,\,y\in A^c} \pi'(x)\,\P'(x,y)}{\rho_e} = \frac{\Q'(A,A^c)}{\rho_e}
$$
because if $(x,y)\in A\times A^c$ then $\gamma_{xy}$ must include some edge $(u,v)\in A\times A^c$.
\end{proof}

It follows that $\~\Phi(A)\geq \rho_e^{-1}\,\frac{\Q'(A,A^c)}{\pi(A)\pi(A^c)}$, and in particular 
when $\pi=\pi'$ then $\~\Phi(A)\geq \~\Phi'(A)/\rho_e$. By equation \eqref{eqn:P-mixing}, if $\MM$ is lazy and $\pi=\pi'$
then $1-\C_{\sqrt{a(1-a)}}(A)\geq 1-\sqrt{1-(\~\Phi'(A)/\rho_e)^2}$. Furthermore, since spectral gap determines the asymptotic rate of convergence, then also if $\MM$ is reversible then
\begin{equation} \label{eqn:new_sinclair}
\lambda \geq 1-\C_{\sqrt{a(1-a)}}\geq 1-\sqrt{1-(\~\Phi'/\rho_e)^2}\geq \~\Phi'^2/2\rho_e^2\,.
\end{equation}

As with comparison methods, by comparing to the complete graph we obtain a ``standard'' canonical path type bound. In this case let $\MM'$ have transitions $\P'(x,y)=\pi(y)$, so that $\Q'(A,A^c)=\pi(A)\pi(A^c)$ and $\~\Phi'(A)=1$. Then 
$$
\lambda\geq 1-\sqrt{1-1/\rho_e^2}\geq 1/2\rho_e^2\,,
$$
a factor of $4$ improvement over Jerrum and Sinclair's bound $\lambda\geq 1/8\rho_e^2$ \cite{JS88.1,Sin92.1}.

Our main result in this section is a comparison theorem in terms of
edge and vertex congestion. This applies to a wide range of distance measures and to
non-reversible Markov chains.

\begin{theorem} \label{thm:comparison}
Suppose $\MM$ and $\MM'$ are finite Markov chains on state space $V$
and edge sets $E$ and $E'$ respectively. To every edge $e'=(x,y)\in E'$ associate a path $\gamma_{xy}\subset E$. Let
$$
\rho_e = \max_{e\in E} \frac{1}{\Q(e)}\,\sum_{\gamma_{xy} \ni e} \pi'(x)\,\P'(x,y)
\quad \textrm{ and }\quad
\rho_v = \max_{v\in V} \frac{1}{\pi(v)}\,\sum_{\gamma_{xy} \ni v} \pi'(x)\,\P'(x,y)\,.
$$
If $\MM$ is lazy, $A\subset V$, and $f:\,[0,1]\rightarrow \R$ is a concave function with $f(\pi(A))\neq 0$ then
$$
1-\C_f(A) \geq 2\,\frac{\rho_v}{\rho_e}\,\left(1 - 
   \frac{f\left(\pi(A)+\frac{\Q'(A,A^c)}{\rho_v}\right) + f\left(\pi(A)-\frac{\Q'(A,A^c)}{\rho_v}\right)}{2\,f(\pi(A))} \right)\, .
$$

In the sum for $\rho_v$ if both paths $\gamma_{xy}$ and $\gamma_{yx}$ pass through the
vertex $v$ then it suffices to replace the sum of weights 
$\pi'(x)\,\P'(x,y)+\pi'(y)\,\P'(y,x)$ by
$\max\{\pi'(x)\,\P'(x,y),\,\pi'(y)\,\P'(y,x)\}$.
\end{theorem}

\begin{proof}
Given a set $A\subset V$, by Lemma \ref{lem:worst-case} to lower bound $1-\C_f(A)$
it suffices to construct the distribution of $\pi(A_u)$ which minimizes $\int_0^t \pi(A_u)\,du$.

Recall that $\int_0^{1/2} \pi(A_u)\,du=\frac{\pi(A)}{2}+\Q(A,A^c)$. Therefore, given the
correct distribution of $\pi(A_u)$ on $[0,1/2]$, and an underestimate of the ergodic flow
$\Q_0(A,A^c)\leq\Q(A,A^c)$, then
the distribution $\pi_0(A_u)=\pi(A)+(\pi(A_u)-\pi(A))\frac{\Q_0(A,A^c)}{\Q(A,A^c)}$ will certainly
satisfy $\forall t\in[0,1/2]:\,\int_0^t \pi_0(A_u)\,du\leq\int_0^t \pi(A_u)\,du$ and 
$\int_0^{1/2} \pi_0(A_u)\,du=\frac{\pi(A)}{2}+\Q_0(A,A^c)$. A similar argument holds for $t\in[1/2,1]$.
It follows that underestimating the ergodic flow by only considering that due to the paths
will only make the lower bound on $1-\C_f(A)$ too pessimistic.

Now, let $\Q_0$ denote the ergodic flow measured by the paths, so that 
$\forall v\in A^c:\,\Q_0(A,v)=\sum_{x\in A,y\in A^c} \delta_{\gamma_{xy}\ni v}\,\pi'(x)\P'(x,y)/\rho_e$ 
and likewise for $v\in A$, so $\Q_0(A,A^c)\geq\Q'(A,A^c)/\rho_e$. 
The vertex congestion implies that not too much of this passes through any specific vertex, 
and in particular
$\forall v\in A:\,\frac{\Q_0(A^c,v)}{\pi(v)} \leq \frac{\Q'(A,A^c)/\rho_e}{\Q'(A,A^c)/\rho_v}=\frac{\rho_v}{\rho_e}$,
while $\forall v\in A^c:\,\frac{\Q_0(A,v)}{\pi(v)} \leq \frac{\rho_v}{\rho_e}$ as well.

It follows that if $M:=\rho_v/\rho_e$ then $\forall u\in[M,1-M]:\,\pi(A_u)=\pi(A)$.
Subject to this constraint, in Lemma \ref{lem:worst-case} the integral $\int_0^t \pi(A_u)\,du$ is
minimized by
$$
\pi(A_u) =
\begin{cases}
\pi(A)+M^{-1}\,\Q(A,A^c) & \textrm{if } u<M \\
\pi(A)                   & \textrm{if } u\in[M,1-M] \\ 
\pi(A)-M^{-1}\,\Q(A,A^c) & \textrm{if } u>1-M
\end{cases}
$$

For these values of $\pi(A_u)$ integration shows that
$$
\int_0^1 f(\pi(A_u))\,du = M\,f\left(\pi(A)+\frac{\Q(A,A^c)}{M}\right) 
    + (1-2M)f(\pi(A)) + M\,f\left(\pi(A)-\frac{\Q(A,A^c)}{M}\right)
$$

By definition $M=\rho_v/\rho_e$, and by the remarks above $\Q(A,A^c)=\Q'(A,A^c)/\rho_e$,
which combined with the integral above gives the result.

The final comment on $\rho_v$ follows because if $v\in A$ then we need only consider congestion among paths 
entering $A$, and as only one of $\gamma_{xy}$ or $\gamma_{yx}$ will enter $A$ then there is no need to count
congestion due to both paths. Likewise if $v\in A^c$. 
\end{proof}

This can be used to show lower bounds for the various $f$-congestion quantities, and hence upper bounds on each notion of
mixing times.

\begin{corollary} [Comparison with Conductance Function] \label{cor:comparison}
Suppose $\MM$ and $\MM'$ are finite Markov chains on state space $V$
and edge sets $E$ and $E'$ respectively. If $A\subset V$ is a proper subset and $\pi=\pi'$ then 
\begin{eqnarray*}
1-\C_{\sqrt{a(1-a)}}(A) &\geq& 2\,\frac{\rho_v}{\rho_e}\,
                         \left(1 - \sqrt{1-\left(\frac{\~\Phi'(A)}{2\rho_v}\right)^2}\right) 
                  \geq  \frac 14\,\frac{\~\Phi'(A)^2}{\rho_v\rho_e} \\
1-\C_{a\log(1/a)}(A)        &\geq& \frac{1}{\log(1/\pi(A))}\,\frac{\Phi'(A)^2}{\rho_v\,\rho_e} \\
1-\C_{a(1-a)}(A)     &\geq& 2\pi(A)\pi(A^c)\,\frac{\~\Phi'(A)^2}{\rho_v\,\rho_e} 
\end{eqnarray*}
If $\pi\neq\pi'$ then replace $\~\Phi'(A)$ by $\frac{\Q'(A,A^c)}{\pi(A)\pi(A^c)}$ 
and likewise $\Phi'(A)$ by $\frac{\Q'(A,A^c)}{\pi(A)}$.
\end{corollary}

\begin{proof}
Using the previous lemma it is just a matter of simplification.

We will save ourselves the work of simplifying by instead reducing this to
a problem that was already solved previously. Observe that if
\begin{equation} \label{eqn:congestion-profile}
\pi(A_u) =
\begin{cases}
\pi(A) - \frac{\Q'(A,A^c)/2\rho_v}{1/2} & \textrm{if } u > 1/2 \\ 
\pi(A) + \frac{\Q'(A,A^c)/2\rho_v}{1/2} & \textrm{if } u \leq 1/2
\end{cases}
\end{equation}
then
$$
1-\C_f(A) = 1 - 
   \frac{f\left(\pi(A)+\frac{\Q'(A,A^c)}{\rho_v}\right) + f\left(\pi(A)-\frac{\Q'(A,A^c)}{\rho_v}\right)}{2\,f(\pi(A))}
$$
is exactly the same as the lower bound on $\frac{1-\C_f(A)}{2\rho_v/\rho_e}$ given in
Theorem \ref{thm:comparison}.

However, the distribution \eqref{eqn:congestion-profile} is the worst case bound
for $\pi(A_u)$ for a lazy chain with ergodic flows $\Q(A,A^c)\geq \Q'(A,A^c)/2\rho_v$
(see Figure \ref{fig:extreme}). We already determined that this leads to the bound 
$1-\C_{\sqrt{a(1-a)}}(A)\geq 1-\sqrt{1-\tilde\Phi(A)^2}$. Therefore,
$$
\frac{1-\C_{\sqrt{a(1-a)}}(A)}{2\rho_v/\rho_e} \geq 1-\sqrt{1-\left(\frac{\Q'(A,A^c)/2\rho_v}{\pi(A)\pi(A^c)}\right)^2}
$$
and the lower bound on $1-\C_{\sqrt{a(1-a)}}(A)$ follows immediately. The other bounds follow similarly.
\end{proof}

If $\MM'=\MM$, $\MM$ is reversible and $\forall x\in V:\,\P(x,x)\geq\gamma\in[1/2,1]$
then $\rho_v\leq(1-\gamma)\rho_e=1-\gamma$ and so 
$$
\~\Phi(A) \geq 1-\C_{\sqrt{a(1-a)}}(A) 
          \geq 2(1-\gamma)\,\left(1-\sqrt{1-\left(\frac{\~\Phi(A)}{2(1-\gamma)}\right)^2}\right)
          \geq \frac{\~\Phi(A)^2}{4(1-\gamma)}\,.
$$

Without reversibility we have only that $\rho_v\leq 2(1-\gamma)\rho_e$,
and a slightly weaker result is obtained.

Corollary \ref{cor:comparison} is all that is needed to prove Corollary \ref{cor:paths}.
In the particular case when $\MM'$ has transition probabilities $\P'(x,y)=\pi(y)$
then $\~\Phi'(A)=1$ and this reduces to a ``canonical paths'' theorem.
\begin{corollary}[Canonical Paths] \label{cor:paths}
Suppose $\MM$ is a finite ergodic lazy Markov chain on state space $V$ and edge set $E$, and $\Gamma$ is a set of paths
$\gamma_{xy}$ between every pair of distinct vertices $x,\,y\in V$. Define
$$
\rho_e = \max_{e\in E} \frac{1}{\Q(e)}\,\sum_{\gamma_{xy} \ni e} \pi(x)\,\pi(y)
\qquad and \qquad
\rho_v = \max_{v\in V} \frac{1}{\pi(v)}\,\sum_{\gamma_{xy} \ni v} \pi(x)\,\pi(y)\,.
$$
Then
$$
1-\C_{\sqrt{x(1-x)}} \geq 2\,\frac{\rho_v}{\rho_e}\,\left(1-\sqrt{1-1/4\rho_v^2}\right) \geq \frac{1}{4\,\rho_v\,\rho_e}
$$
and
$$
\tau(\epsilon) \leq 4\,\rho_v\,\rho_e\,\log \frac{1}{2\epsilon\sqrt{\pi_*}}\,.
$$
If $\MM$ is reversible then $\rho_e$ and $\rho_v$ can be taken as sums over
unordered pairs $(x,y)$ with undirected edges, and moreover
$$
\lambda \geq 1-\C_{\sqrt{x(1-x)}} \geq 2\,\frac{\rho_v}{\rho_e}\,\left(1-\sqrt{1-1/4\rho_v^2}\right) \geq \frac{1}{4\,\rho_v\,\rho_e}\,,
$$
where $\lambda = 1-\lambda_2$ is the {\em spectral gap} between $1$ and the second largest 
eigenvalue $\lambda_2$ of the transition matrix $\P$.
\end{corollary}

\begin{proof}
Let $\MM'$ be a walk on the state space $V$ 
with transition probabilities $\P'(x,y) = \pi(y)$. Then $\pi' = \pi$,
$\Q'(A,A^c) = \pi(A)\pi(A^c)$, to every edge in $\MM'$ (i.e. pair of vertices)
is associated a path $\gamma_{xy}$ given in the corollary, and the congestions are
exactly the $\rho_e$ and $\rho_v$ given in the corollary. The bounds on $\lambda$,
$1-\C_{\sqrt{a(1-a)}}$ and $\tau(\epsilon)$ then follow from Corollary \ref{cor:comparison}.
\end{proof}

As mentioned earlier, for a lazy chain $\rho_v\leq \rho_e/2$,
which reduces the result to $\lambda\geq 1-\C_{\sqrt{a(1-a)}} \geq 1-\sqrt{1-\rho_e^{-2}}$, which
we already know is sharp for the lazy walk on the uniform two-point space. Our extension shows that
the second edge congestion term in
the canonical paths bound $\tau = O(\rho_e\cdot \rho_e\,\log \pi_0^{-1})$ \cite{JS88.1,Sin92.1} 
should actually have been the smaller vertex congestion $\rho_v$. This can be a significant
improvement, as the following example shows.

\begin{example} \label{ex:matroid}
Feder and Mihail \cite{FM92.1} studied a random walk for sampling balanced matroids
and showed a result equivalent to $\rho_e \le n\,m$ and $\rho_v \le 2n$.
The Cheeger bound (Jerrum and Sinclair \cite{JS88.1}) implies a bound of
$\tau(\epsilon) \leq 8\rho_e^2\,\log \frac{1}{2\epsilon\sqrt{\pi_*}} 
  \leq 8m^2\,n^2\,\left(\frac{n}{2}\log m + \log(1/2\epsilon)\right)$. 
By Corollary \ref{cor:paths} our results show the stronger bound
$\tau(\epsilon) \leq 8\,m\,n^2\,\log \frac{1}{2\sqrt{\pi_*}\epsilon} 
    \leq 8\,m\,n^2\,\left(\frac{n}{2}\,\log m + \log (1/2\epsilon)\right)$,
exactly the same upper bound obtained by Feder and Mihail \cite{FM92.1} by using
a modified form of Poincar\'e. 

Also, in \cite{MS01.1} we showed that $\Phi(A) \geq \frac{\log_2(1/\pi(A))}{2m\,n}$, which by
equation \eqref{eqn:P-mixing}
implies that $\tau(\epsilon) \leq 4(\log 2)m^2\,n^2 + 8m^2\,n^2\,\log (1/2\epsilon)$,
not particularly good. However, $1-\C_{\sqrt{a(1-a)}}(A)\geq \Phi(A)^2/2$, and combining this with the 
canonical path lower bound on $1-\C_{\sqrt{a(1-a)}}(A)$ implies that
\begin{eqnarray*}
1-\C_{\sqrt{a(1-a)}}(A) 
  &\geq& \max\left\{\frac{1}{8m\,n^2},\,\frac{\log_2^2(1/\pi(A))}{8m^2n^2}\right\} \\
  &=&
\begin{cases}
\log_2^2(1/\pi(A))/8m^2n^2 & if\ \pi(A)\leq 2^{-\sqrt{m}} \\
1/8m\,n^2                  & if\ \pi(A) >   2^{-\sqrt{m}}
\end{cases}
\end{eqnarray*}
and this time the mixing time is $\tau(\epsilon) \leq 8(\log 2)m^{3/2}\,n^2 + 8m\,n^2\,\log (1/2\epsilon)$
when $2^{-\sqrt{m}}\geq m^{-n}$, an improvement over the canonical paths
bound for all simple balanced matroids (as $m\leq \binom{n}{2}$).

We note that the correct bound is still much smaller at $\tau(\epsilon)=O\left(m\,n\,\log \frac{n}{\epsilon}\right)$
\cite{JS02.1}.
\end{example}

It is apparent from the definitions that $\rho_v \geq \rho_e\,\min\{\P(x,y):\,x\neq y,\,\P(x,y)\neq 0\}$, 
since if $\vec{e}$ is the worst directed edge then $v$ can be taken as one of the endpoints. 
In the following example we show how a common enumeration process sometimes leads to this lower bound
being achieved.

\begin{example} \label{ex:permanent}
One of the first Markov chains analyzed via the canonical path method was a Markov chain of
Broder \cite{Br86.1} for approximating the permanent of a dense matrix, or equivalently counting
perfect matchings in a bipartite graph $G=K_{n,n}$ of minimum degree $n/2$. Jerrum and Sinclair \cite{JS88.1}
constructed canonical paths and used a clever enumeration process to show that $\rho_e \leq 12\,n^6$. 
We give only a rough sketch of where the $\rho_v$ computation differs from the $\rho_e$ computation,
and refer the reader to \cite{JS88.1} for further details.


List all cycles in $G$, assign them an ordering, and to each cycle $C$ fix a ``starting vertex'' $v_c$.
As done by Jerrum and Sinclair, to each pair of perfect matchings the path $\gamma_{IF}$ is given by 
considering the symmetric difference $I\oplus F$, and then unwinding the resulting cycles in the 
cycle ordering just given.  Let $M$ be a fixed vertex (a perfect or near-perfect matching) and suppose 
$\gamma_{IF}$ is a canonical path joining two perfect matchings and passing through $M$. 
If $M$ is near-perfect then let $e_{IM}$ denote the edge of $I$ incident with the starting vertex of the 
cycle being unwound when $M$ was reached, and $e_{FM}$ be the edge of $F$ that will be added at the next 
step of the unwinding. The encoding we use is
$$
\sigma_M(I,F) =
\begin{cases}
  I\oplus F\oplus M           & \textrm{if $M$ is perfect;} \\
  I\oplus F\oplus M - e_{IM} & \textrm{if $M+e_{IM}$ is perfect;} \\
  I\oplus F\oplus M - e_{IM} - e_{FM} & \textrm{otherwise}
\end{cases}
$$
Checking a few cases, as in \cite{JS88.1}, shows that $\sigma_M(I,F)$ is a perfect or near-perfect matching.

To see that this is injective we define the decoding process. If $\sigma_M(I,F) \oplus M$ is a sequence of
cycles then $I\oplus F = \sigma_M(I,F) \oplus M$. If $\sigma_M(I,F) \oplus M$ is a sequence of
cycles plus a path containing the remaining vertices then let $e$ be the edge joining the
endpoints of the path, and $I\oplus F = \sigma_M(I,F) \oplus M + e$. The only remaining case is  if 
$\sigma_M(I,F) \oplus M$ is a sequence of cycles, plus two paths; each path will have one end
in $M$ and one end not in $M$, join the paths to create a cycle by connecting the ends in $M$ to those not 
in $M$, and this gives $I\oplus F$. 

Given $I\oplus F$ then the order in which the unwinding occurred follows from
the cycle ordering. The matching $I$ contains all edges removed before $M$ was reached, plus all edges
in $M$ for the remaining cycles. The matching $F$ is the converse.

It follows that the enumeration given above measures all paths connecting two perfect matchings and
passing through $M$, not just those through an edge as in \cite{JS88.1}. The computation of Jerrum and
Sinclair then holds for $\rho_v$ as well, but without the need for the $\P(M,M')$ term in the denominator.
That is, $\rho_v \leq 12n^6\,\max \P(M,M') = 6\,n^4$ and therefore $\lambda \geq 1/288\,n^{10}$, an improvement over
the $\lambda \geq 1/1152\,n^{12}$ of \cite{JS88.1}.
\end{example}

The $O(n^2)$ improvement over the results of \cite{JS88.1} is nice and is the maximum
possible with our path bound, as $\rho_v \geq \rho_e\,\min\{\P(x,y):\,x\neq y,\,\P(x,y)\neq 0\}$. However, this is
not as large as the $O(n^5)$ improvement possible by use of a Poincar\'e bound $\lambda\geq 1/\rho_e\ell$, as
in \cite{Sin92.1,DS91.1}. In the following section a comparison theorem is proven in terms of maximum path length via Blocking conductance,
and hence will match the Poincar\'e bound, at least up to order of magnitude. Unfortunately, we have been 
unable to derive a similar theorem for the Evolving sets bounds. 

The reason our $\rho_v$ based bound is sometimes as good as using path length, sometimes not, and occasionally better, can be understood better by the following corollary.
\begin{corollary} \label{cor:path_length}
Let $\rho_e$ and $\rho_v$ be as in Corollary \ref{cor:paths}, and also let
$$
\rho_v^{ave} = \sum_{v\in V} \pi(v) \left[ \frac{1}{\pi(v)} \sum_{\gamma_{xy}\ni v} \pi(x)\,\pi(y) \right]
$$ 
be the average vertex congestion over the entire space $V$. Then
$$
\tau(\epsilon) \leq 2\,\left(\frac{\rho_v(1-\|\pi\|_2^2)}{\rho_v^{ave}}\right)\,
                    \rho_e\,\ell_{ave}\,\log \frac{1}{2\epsilon\sqrt{\pi_*}}
$$
where $\ell_{ave} = \frac{\sum_{x\neq y\in V} \pi(x)\,\pi(y)\,|\gamma_{xy}|}{\sum_{x\neq y\in V} \pi(x)\pi(y)}$ 
is the average length of the canonical paths.
\end{corollary}

\begin{proof}
Changing the order of summation gives
\begin{eqnarray*}
\rho_v^{ave} =& \displaystyle \sum_{x\neq y\in V} \sum_{v \in \gamma_{xy}} \pi(x)\,\pi(y) 
       &= \sum_{x\neq y\in V} \pi(x)\,\pi(y)\,(|\gamma_{xy}|+1)  \\
       =& \displaystyle (\ell_{ave}+1)\,\sum_{x\neq y\in V} \pi(x)\,\pi(y)
       &= (\ell_{ave}+1)(1-\|\pi\|_2^2)\ .
\end{eqnarray*}
Finish by multiplying the upper bound on $\tau(\epsilon)$ in Corollary \ref{cor:paths} by
$1 \leq \frac{2\ell_{ave}(1-\|\pi\|_2^2)}{\rho_v^{ave}}$.
\end{proof}

A similar result for the Comparison theorem also holds, but with the average path length
$$
\ell_{ave} = \frac{\sum_{(x,y)\in E'} \pi'(x)\,\P'(x,y)\,|\gamma_{xy}|}{\sum_{(x,y)\in E'} \pi'(x)\,\P'(x,y)}\ .
$$

This shows that when the canonical paths are well distributed among the vertices,
as in Example \ref{ex:matroid}, then Corollary \ref{cor:comparison} can be as strong 
as Poincar\'e bounds. The Corollary suggests that if the paths are short but concentrated 
on a few vertices then our results will be poor. On the other hand, if $\|\pi\|_2\approx 1$,
or if the paths are fairly well distributed among the vertices but there are some paths 
much longer than average then our vertex bounds may be better. 
The following example satisfies both of these conditions.
\begin{example}
Consider the lazy asymmetric walk on a line $P_n=[1\ldots n]$ given by
$\P(i,i+1)=\frac{c}{2}\leq 1/2$ and $\P(i,i-1)=\frac{1-c}{2}\leq 1/2$,
with $\P(1,1)=\frac 12+\frac{1-c}{2}$ and $\P(n,n)=\frac 12+\frac c2$.
The stationary distribution is $\pi(i)=\beta\,\alpha^i$ where
$\alpha = \frac{c}{1-c}$ and 
$\beta^{-1} = \sum_{i=1}^n \alpha^i = \frac{c}{1-2c}\,\left(1-\left(\frac{c}{1-c}\right)^n\right)$.

In particular, when $c<1/2$ then $1-\|\pi\|_2^2 \xrightarrow{n\rightarrow\infty} 2c$, and
if $c\approx 0$ then Corollary \ref{cor:path_length} suggests that using $\rho_v$ should be much
better than using path length $\ell$.
There is an obvious choice of canonical paths given by $\gamma_{ij}=i,\,i+1,\,\ldots,j$
when $i<j$, and vice-versa when $i>j$. It is easily verified that
$\ell=n-1$, $\rho_e\rightarrow 2$ and $\rho_v\rightarrow 1$ when $c\rightarrow 0^+$ and $n\rightarrow\infty$.
Then $\lambda \geq 1/2\rho_e^2 = 1/8$ and $\lambda \geq 1/4\rho_v\rho_e = 1/8$
give the same bounds, whereas a Poincar\'e inequality implies only that
$\lambda \geq 1/\rho_e\ell = 1/2(n-1)$, far worse. The correct value is $\lambda\rightarrow 1/2$.
\end{example}

\section{An aside into Blocking Conductance and comparison} \label{sec:blocking-comparison}

One advantage of the Spectral Profile method over that of Evolving sets is that it is fairly straightforward to show that one can compare Spectral Profile of two Markov chains, whereas in the Evolving set case we do not know how to compare $f$-congestions directly. Because of the similarity between Evolving Set and Blocking Conductance bounds we showe here an unpublished comparison method for Blocking Conductance, in the hope that it may help give insight into future work in proving such a result for Evolving Sets as well.

Blocking Conductance only applies to reversible chains and the constant factors are significantly weaker than with our evolving set results, so we skip giving proofs in terms of $\rho_e$ and $\rho_v$ as in the previous section, although such bounds are not hard to show. Instead, we only consider 
a Poincar\'e type comparison theorem 
because this is a situation in which we may genuinely improve on known results. In particular, when studying problems for which total-variation mixing time is faster than $L^2$ mixing, such as the application of Blocking Conductance to studying $G(n,p)$ given by Fountoulakis and Reed \cite{FR07.2}, a Blocking Conductance comparison theorem might show stronger bounds than the $L^2$-based comparison of Spectral Profile method.

Our results will be shown with an earlier form of Blocking Conductance from
our Ph.D. Dissertation \cite{Mon02.1}, since it seems best suited for our purposes here.

Given a finite Markov chain $\MM$ with state space $V$ of cardinality $n=|V|$,
let $(V,<)$ denote an ordering of the space, say as $1,\,2,\ldots,\, n$, 
let $A_i$ denote the subsets $A_i=[1\ldots i]$, and $A_i^+=A_{i+1}=[1\ldots i+1]$.
Then the blocking conductance theorem says:

\begin{theorem} \label{thm:finite-blocking}
Given a lazy, reversible, ergodic, finite Markov chain $\MM$ there exists some ordering $(V,<)_{\MM}$
of the vertices such that if $S_0=[1\ldots k]$ is the largest $A_i$ such that $\pi([1\ldots k])\leq 3/4$, and if 
$S_0\supset S_1 \supset \cdots \supset S_m=[1]$ is any nested sequence of $A_i$,
then the mixing time is bounded by
$$
\tau(\MM) \leq K'\, \sum_{i=1}^{m} \frac{\pi(S_i)}{\Q(S_{i+1}^+,S_i^c)} \,,
$$
where $K'=1376$ arises from converting between different measures of mixing time.
\end{theorem}

Since the ordering is not known in advance it is necessary to 
maximize the right hand side over all possible orderings. We use the notation 
\begin{equation} \label{eqn:nested2}
\tau_{BC}(\MM) = K'\,\max_{orderings\ (V,<)} \min_{\substack{nested\ sequences\\ S_0\supset S_1\supset\cdots\supset S_m=[1]}} \sum_{i=1}^{m} \frac{\pi(S_i)}{\Q(S_{i+1}^+,S_i^c)} \,,
\end{equation}
that is, the best possible upper bound on mixing time when the ordering is not known. In Remark \ref{rmk:integral-blocking} some more conventional forms of this relation are discussed.

Our main result of this section is then:
\begin{theorem}[Comparison with Blocking Conductance] \label{thm:comparison2}
Suppose that $\MM$ and $\MM'$ are finite Markov chains with the same set of vertices $V$,
the same stationary distribution $\pi$, and suppose
that to every edge $e'=(x,y)\in E'$ is associated a path $\gamma_{xy}\subset E$.
Let
$$
\rho_e = \max_{e\in E} \frac{1}{\Q(e)}\,\sum_{\gamma_{xy} \ni e} \pi(x)\,\P'(x,y)
$$
be a measure of edge congestion and $\ell =\max_{(x,y)\in E'} |\gamma_{xy}|$ be the length 
of the longest path. 

Then, 
$$
\tau_{BC}(\MM) \leq 4\,\rho_e\,\ell\,\tau_{BC}(\MM')
$$
where $\tau_{BC}$ is given by equation \eqref{eqn:nested2}.
\end{theorem}

\begin{proof}
Fix some ordering $(V,<)$ of the state space $V$, as discussed in the Preliminaries. 
Let $S_0\supset S_1 \supset \cdots \supset S_m=[1]$ be the nested sequence of
sets $A_i=[1\ldots i]$ which minimizes 
$$
\sum_{i=1}^{m} \frac{\pi(S_i)}{\Q(S_{i+1}^+,S_i^c)} \,,
$$
that is the sum in the upper bound for $\tau_{BC}(\MM)$. Observe that every set $S$ with 
$S_1^+\subseteq S\subseteq S_0$ satisfies
$\Q(S,S^c) \geq \Q'(S,S^c)/\rho_e \geq \Q'(S_1^+,S_0^c)/\rho_e$.

Let $B_0=S_0$, and let $B_1$ be the smallest set $A_i=[1\ldots i]\subset B_0$ such that 
$\Q(B_1^+,B_0^c) \geq \Q'(S_1^+,S_0^c)/2\rho_e$.
Likewise, let $B_2$ be the smallest $A_i\subset B_1$ such that $\Q(B_2^+,B_1^c) \geq \Q'(S_1^+,S_0^c)/2\rho_e$.
Continue until either $B_i\subseteq S_1$, or until $B_{2\ell }$ is defined. 

Suppose that sets up to $B_{2\ell}$ were defined and $B_{2\ell}\supset S_1$. 
Ergodic flow from set $B_{i+1}$ then satisfies the relation $\Q(B_{i+1},B_i^c) < \Q'(S_1^+,S_0^c)/2\rho_e$, since $B_{i+1}$ was defined to be the smallest
initial sequence $[1\ldots k]$ such that $\Q(B_{i+1}^+,B_i^c)\geq\Q'(S_1^+,S_0^c)/2\rho_e$. By definition
of $\rho_e$ the ergodic flow $\Q(B_{i+1},B_i^c)$ therefore contains under half the canonical paths from $S_1^+$ to 
$S_0^c$, and so more than half these paths pass through $B_i\setminus B_{i+1}$. It follows that if 
$Cut'(S_1^+,S_0^c)$ denotes the edges passing from $S_1^+$ to $S_0^c$ in $\MM'$, then
\begin{eqnarray*}
\sum_{(x,y)\in Cut'(S_1^+,S_0^c)} \pi(x)\,\P'(x,y)\,(|\gamma_{xy}|-1) 
 &\geq& \sum_{i=0}^{2\ell -1} 
        \sum_{\substack{(x,y)\in Cut'(S_1^+,S_0^c),\\ \gamma_{xy} \cap \left(B_i\setminus B_{i+1}\right) \neq \emptyset}} 
            \pi(x)\,\P'(x,y) \\
  &>&   2\,\ell \,\Q'(S_1^+,S_0^c)/2 \\
  &\geq& \sum_{(x,y)\in Cut'(S_1^+,S_0^c)} \pi(x)\,\P'(x,y)\,|\gamma_{xy}| \ .
\end{eqnarray*}
This gives a contradiction, so it follows that $B_{2\ell}\subseteq S_1$.

By construction, if sets up to $B_{2\ell}$ are defined then
\begin{equation} \label{eq:covering}
\sum_{i=0}^{2\ell -1} \frac{\pi(B_i)}{\Q(B_{i+1}^+,B_i^c)}
   \leq 4\,\rho_e\,\ell \,\frac{\pi(S_0)}{\Q'(S_1^+,S_0^c)} \ .
\end{equation}
If fewer $B_i$ were required then the sum is even smaller. Also,
it can be assumed that the last $B_i$ is equal to $S_1$, because if say
$B_{2\ell}\subsetneq S_1$ then increasing the size to $B_{2\ell}=S_1$ 
can only decrease the sum in the left side of (\ref{eq:covering}), which reinforces 
the inequality. 

Repeating this process for the other sets $S_k$ leads to a sequence of sets in 
$\MM$ which is at most $4\,\rho_e\,\ell$ 
times worse than the $S_k$ sequence in $\MM'$. This holds for any ordering $(V,<)$ and so 
the theorem follows.
\end{proof}

Just as with Diaconis and Saloff-Coste's comparison theorem \cite{DSC93.1},
comparison with the trivial chain $\P'(x,y)=\pi(y)$ gives the familiar bound 
$\tau(\MM)=O(\rho_e\ell\,\log \frac{1}{\pi_*})$. Therefore, at least as far as order of
magnitude is concerned, Blocking Conductance matches the canonical path bounds found in most applications. However, the comparison method may be superior to other methods when the total variation mixing time is smaller than the $L^2$ mixing time, such as with $G_{n,p}$ \cite{FR07.2}.

\begin{remark}\label{rmk:integral-blocking}
It is easier to understand the Blocking conductance theorem in an alternate form. First, a definitition. For $A\subset V$ let
$$
{\mathcal B}(A) = \sup_{\alpha \leq \pi(A)} \min_{\substack{B\subset A,\\ \pi(B)\geq \pi(A)-\alpha}}
                       \frac{\alpha\,\Q(B,A^c)}{\pi(A)^2}\ .
$$
The quantity $\Q(B,\,A^c)/\pi(B)$ is a lower bound on the
probability of leaving $A$ in a single step, conditioned on the initial point being drawn
from $B$. So ${\mathcal B}(A)$ roughly measures the size of a potential vertex bottleneck
relative to the size of $A$, times the probability of stepping over this bottleneck.
For example, if $\alpha = \frac 12\,\Q(A,A^c)$ then $\Q(B,\,A^c) \geq \frac 12\,\Q(A,A^c)$
when $\pi(B)\geq\pi(A)-\alpha$, and so ${\mathcal B}(A) \geq \frac 14\,\Phi(A)^2$.
The name {\em Blocking conductance} comes from this connection to conductance and the fact
that $\pi(B)\geq\pi(A)-\alpha$ and so $\alpha$ measures the size of a subset of $A$ which is 
``blocked'' when computing ${\mathcal B}(A)$.

Given an ordering $(V,<)$ a good sequence of $S_i$ is easily enough constructed from 
${\mathcal B}(\cdot)$. Given $S_i$ then let $S_{i+1}$ be such that $\pi(S_i\setminus S_{i+1}) \approx \alpha$.
The following is a slight improvement on a result of the author \cite{Mon02.1}.

\begin{theorem}
If $\MM$ is a finite, irreducible, reversible, lazy Markov chain then
$$
\tau(\MM) \leq \tau_{BC}(\MM) \leq \frac 32\,K'\,\int_{\pi_0/2}^{3/4} \frac{dx}{x\,\phi(x)}\,
$$
where $\K'=1376$ arises from converting between different measures of mixing time
and the {\em blocking conductance function} $\phi$ is given, for $x\in[\pi_0/2,3/4]$, by
$$
\phi(x) = \min_{\substack{\pi(A)\in [x,2x],\\\pi(A)\leq 3/4}} {\mathcal B}(A)\,.
$$
\end{theorem}
\end{remark}


\section{Spread and the discrete gradients} \label{sec:more_bounds}

The quantity we consider here was proposed by Kannan, Lov\'asz and Montenegro \cite{KLM06.1} in the context of Blocking Conductance, in an alternative form by Morris and Peres \cite{MP05.1}, and also used by Montenegro \cite{Mon05.1}. As a bit of motivation, recall that earlier isoperimetric bounds, such as those involving conductance or modified conductance, were shown by explicitly constructing the worst case for $\pi(A_u)$ and then applying Lemma \ref{lem:worst-case}. We now show that useful relations can be derived even when Lemma \ref{lem:worst-case} appears not to be appicable. 
We examine only the relation between the Evolving set bounds and the spread. The interested reader can see Montenegro \cite{Mon05.1} for an examination of the relation between spread, the discrete gradients of Houdr\'e and Tetali \cite{HT04.1}, and spectral gap.

Two isoperimetric quantities will be used, extending earlier definitions. 
Recall from definition \eqref{eqn:psiA2} that $\Psi(A,t)$ is the smallest flow from from $A$ into a subset of $V$ of size $t$, so in particular for a lazy walk $\Psi(A,\pi(A^c))=\Q(A,A^c)$.
\begin{definition}
If $A\subset V$ then the {\em internal and external spread} are given by
$$
\~\psi^+(A) = \int_0^{\pi(A)} \frac{\Psi(A^c,t)}{\pi(A)^2\,\pi(A^c)^2}\,dt
\quad and \quad
\~\psi^-(A) = \int_0^{1-\pi(A)} \frac{\Psi(A,t)}{\pi(A)^2\,\pi(A^c)^2}\,dt\,.
$$
Quantities $\psi^{\pm}(A)$ are defined similarly but without $\pi(A^c)^2$ in the
denominators.
\end{definition}


The spread turns out to fairly closely bound many natural choices of $\C_f$.
This was first observed in \cite{Mon05.1} where the connection between
$\C_{\sqrt{a}}$ and spread were studied for lazy Markov chains in order to relate 
Blocking Conductance and Evolving Set results.

\begin{theorem} \label{thm:gradients}
Given a finite irreducible Markov kernel and $A\subset V$ with $\pi(A)\leq 1/2$ then
\begin{eqnarray*}
\~\psi^+(A)\pi(A^c) + \~\psi^-(A) \,\max\left\{\pi(A),\,\frac 14\right\}
\geq& \displaystyle 1-\C_{\sqrt{a(1-a)}}(A) 
&\geq \frac 14\,\~\psi^+(A)+\~\psi^-(A)\,\pi(A)\pi(A^c)  \\
\frac{2\psi^+(A)+\psi^-(A)}{\log(1/\pi(A))}
\geq& \displaystyle 1-\C_{a\log(1/a)}(A)
&\geq \frac{\psi^+(A)}{\log(1/\pi(A))} + \~\psi^-(A)\,\pi(A)\pi(A^c)
\end{eqnarray*}
$$
and\ 1-\C_{a(1-a)}(A) = 2\,\left(\~\psi^+(A)+\~\psi^-(A)\right)\,\pi(A)\pi(A^c)\,.
$$
The bounds for $1-\C_{a\log(1/a)}(A)$ and $1-\C_{a(1-a)}(A)$ hold when $\pi(A)>1/2$ as well.
\end{theorem}

It follows that $1-\C_{a\log(1/a)}(A),\,1-\C_{\sqrt{a(1-a)}}(A)\geq \frac 12\,(1-\C_{a(1-a)}(A))$.

Observe that $\Psi(A^c,\pi(A)-t)\geq \Psi(A^c)-t=\Psi(A)-t$, and so
$$
\~\psi^+(A)
 \geq \frac{\int_0^{\Psi(A)} (\Psi(A)-t)\,dt}{\pi(A)^2\pi(A^c)^2}
  =   \frac 12\,\~\phi(A)^2\,.
$$
Likewise, $\~\psi^-(A)\geq\frac 12\,\~\phi(A)^2$.
Therefore, at least up to a small multiplicative factor, this supercede's Theorem \ref{thm:profile}.

\begin{remark} \label{rmk:spread-vertex}
As discussed in Remark \ref{rmk:blocking-vertex}, for a lazy walk spread incorporates measures of both edge and vertex expansion, whereas (modified) conductance involves only edge expansion. This can be generalized to a non-lazy walk as well. To see this, given set $A\subset V$ let $\beta(A)=\max\{t:\,\Psi(A^c,\pi(A)-t)\geq \Psi(A^c)/2\}$.
Then $\~\psi^+(A) \geq \frac{\beta(A)\,\Psi(A)}{2\,\pi(A)^2\pi(A^c)^2}$ and so
$$
\~\phi \geq 1-\C_{\sqrt{a(1-a)}}(A) \geq \min_{\pi(A)\leq 1/2} \frac 14\,\~\psi^+(A) 
\geq \min_{\pi(A)\leq 1/2} \frac 18\,\frac{\beta(A)}{\pi(A)\pi(A^c)}\,\~\phi(A)\,.
$$
The quantity $\beta(A)$ is a notion of vertex expansion that measures how large a set of vertices must 
be ``blocked'' so that only half of the flow $\Psi(A^c)$ remains. 
Therefore $\~\psi^+(A)$ can be thought of as a product of edge and vertex expansion.
Since $\frac{\pi(A)}{2}\geq\beta(A)\geq\frac{\Psi(A^c)}{2}$ then the above lower bound may be the same order as the upper bound, and is at worst $4$ times weaker than
our modified conductance lower bound of Theorem \ref{thm:profile}. However, when there are many boundary
vertices and $\beta(A)\gg\Psi(A^c)$ then this can be substantially better.
\end{remark}

\begin{proof}
The proof involves working with the spread written in a form involving Evolving Sets:
$$
\~\psi^+(A) = \frac 12\,\int_{\wp_A}^1 \frac{(\pi(A_u)-\pi(A))^2}{\pi(A)^2\pi(A^c)^2}\,du
,\quad
\~\psi^-(A) = \frac 12\,\int_0^{\wp_A} \frac{(\pi(A_u)-\pi(A))^2}{\pi(A)^2\pi(A^c)^2}\,du
$$
The first of these was derived in equation \eqref{eqn:chi2-mod}; the $\~\psi^-(A)$ relation follows in exactly the same way.

The equality for $1-\C_{a(1-a)}(A)$ follows immediately from the form in Lemma \ref{lem:blocking} and
the definitions of $\~\psi^{\pm}(A)$. Alternatively, as an introduction to the method used in the remainder of the proof, start with the identity
$$
\forall x,y\in[0,1]:\,y(1-y) = x(1-x) + (1-2x)\,(y-x) - (y-x)^2\,.
$$
Letting $y=\pi(A_u)$, $x=\pi(A)$ and integrate over $u\in[0,1]$. Finish by applying the Martingale property that $\int_0^1 (\pi(A_u)-\pi(A))\,du=0$, and dividing by $\pi(A)\pi(A^c)$.


The inequality required to show the lower bound on $1-\C_{a\log(1/a)}(A)$ is
\begin{equation} \label{eqn:C_{D}-lowerbound}
y\log \frac 1y \leq x\log \frac 1x + \left(\log \frac 1x - 1\right)\,(y-x)
- \frac{(y-x)^2}{2}\,\left(\delta_{y\leq x}\,\frac 1x + \delta_{y>x}\,\frac{\log \frac 1x - (1-x)}{(1-x)^2/2}\right)
\end{equation}
for all $x,\,y\in[0,1]$ (proven below). Let $y=\pi(A_u)$, $x=\pi(A)$, recall $\int_0^1(\pi(A_u)-\pi(A))\,du=0$, and use the form of $\C_{a\log(1/a)}(A)$ derived in Lemma \ref{lem:blocking},
\begin{eqnarray*}
\C_{a\log(1/a)}(A) &=& \int_0^1 \frac{\pi(A_u)\log (1/\pi(A_u))}{\pi(A)\log (1/\pi(A))}\,du \\
&\leq& 1 - \frac{\psi^+(A)}{\log (1/\pi(A))} - 2\~\psi^-(A)\,\pi(A)\pi(A^c)\,\left(\frac {1}{\pi(A^c)} - \frac{1}{\log(1/\pi(A))}\right)
\end{eqnarray*}
To finish apply the inequality 
$\forall x\in[0,1]:\,\frac{2}{1-x}-\frac{2}{\log(1/x)}\geq 1$ with $x=\pi(A)$.

The upper bound on $1-\C_{a\log(1/a)}(A)$ follows similarly, but with
$$
y\log \frac 1y \geq x\log \frac 1x + \left(\log \frac 1x - 1\right)\,(y-x)
- \frac{(y-x)^2}{2x}\,(1 + \delta_{y\leq x})
$$
for all $x,\,y\in[0,1]$.

For the lower bound on $1-\C_{\sqrt{a(1-a)}}(A)$ use the inequality
$$
\sqrt{y(1-y)} \leq \sqrt{x(1-x)} + \frac{1-2x}{2\sqrt{x(1-x)}}\,(y-x) 
- \frac{(y-x)^2}{2\,[x(1-x)]^{3/2}}\,\left(\delta_{y\leq x}\,\frac 14 + \delta_{y>x}\,x(1-x)\right)
$$
for $x\in[0,1/2],\,y\in[0,1]$. This relation follows from the inequality $\sqrt{z}\leq 1+\frac 12(z-1)-\frac 18(z-1)^2\delta_{z\leq 1}$ with $z=\frac{y(1-y)}{x(1-x)}$. A calculation as done for $\C_{a\log(1/a)}(A)$ gives the result.

The upper bound on $1-\C_{\sqrt{a(1-a)}}(A)$ uses
$$
\mbox{$
\sqrt{y(1-y)} \geq \sqrt{x(1-x)} + \frac{1-2x}{2\sqrt{x(1-x)}}\,(y-x)
- \frac{(y-x)^2}{2\,[x(1-x)]^{3/2}}\,\left(\delta_{y\leq x}\,(1-x)+\frac{\delta_{y\in(x,1-x)}}{4}+\delta_{y\geq 1-x}\,x\right)
$}
$$
for $x\in[0,1/2],\,y\in[0,1]$. This relation follows from the inequality  
$\sqrt{z} \geq 1 + \frac 12\,(z-1) - \frac 12\,(z-1)^2\,\delta_{z\leq 1} - \frac 18\,(z-1)^2\,\delta_{z\geq 1}$, again with $z=\frac{y(1-y)}{x(1-x)}$. Consolidate a bit via the relation $\frac{\delta_{y\in(x,1-x)}}{4}+\delta_{y\geq 1-x}\,x\leq\delta_{y\geq x}\max\{x,1/4\}$. Finish again as in $\C_{a\log(1/a)}(A)$.

We finish with a proof of equation (\ref{eqn:C_{D}-lowerbound}). Consider
$$
g(x,y)
  = y\log \frac 1y - x\log \frac 1x - \left(\log \frac 1x-1\right)(y-x)+\frac{(y-x)^2}{2}\,(h_1(x)\delta_{y\leq x}+h_2(x)\delta_{y>x}) 
$$
where $h_1(x)=1/x$ and $h_2(x)=\frac{\log \frac 1x - (1-x)}{(1-x)^2/2}$.
We will show that $g(x,y)\leq 0$ for $x,y\in[0,1]$ by showing that $g(x,y)$ is increasing
for $y\in[0,x]$, while for $y\in(x,1]$ it decreases and then increases. It is easily verified
that the inequality holds at $y=x$, $y\rightarrow x^+$ and $y=1$, so the result then follows.
First, calculate a few derivatives.
\begin{eqnarray*}
\frac{dg}{dy}     &=& \log \frac{x}{y} + (y-x)\,\left(h_1(x)\delta_{y<x}+h_2(x)\delta_{y>x}\right) \\
\frac{d^2g}{dy^2} &=& -\frac{1}{y} + h_1(x)\delta_{y<x}+h_2(x)\delta_{y>x} \\
\frac{d^3g}{dy^3} &=& \frac{1}{y^2} > 0
\end{eqnarray*}
Consider $y\in[0,x)$. Since $h_1(x)=1/x$ then $\left.\frac{d^2g}{dy^2}\right|_{y<x}\leq 0$ and so
$g(x,y)$ is concave in $y$. But $\left.\frac{dg}{dy}\right|_{y\rightarrow x^-} = 0$ and 
so $g(x,y)$ is increasing in $y$ for $y\in[0,x]$, as desired. 
Now consider $y\in(x,1]$.
The third derivative is positive, so $\frac{d^2g}{dy^2}$ is increasing in $y$.
But $\left.\frac{d^2g}{dy^2}\right|_{y\rightarrow x^+}\leq 0$ and $\left.\frac{d^2g}{dy^2}\right|_{y=1}\geq 0$
and so $g$ is initially concave and transitions to convex. Since $\left.\frac{dg}{dy}\right|_{y\rightarrow x^+}=0$
then $g$ decreases, then if it transitions to convex then it may increase later.
This completes the proof of equation (\ref{eqn:C_{D}-lowerbound}).
\end{proof}

The bounds are sharp. 

The worst case for the lower bound of $1-\C_{a\log(1/a)}(A)$ is
when the flow leaves a small sliver of $A$ and flows uniformly into $A^c$,
that is $\Psi(A^c,t)=\left(t-(\pi(A)-\Q(A,A^c))\right)^+$ when $t\leq \pi(A)$
and $\Psi(A^c,t) = \frac{1-t}{1-\pi(A)}\,\Q(A,A^c)$ when $t\geq \pi(A)$.
For the upper bound this is reversed, with the flow leaving uniformly from $A$ and 
concentrated in a sliver of $A^c$. This is also sharp on $1-\C_{\sqrt{a(1-a)}}(A)$  

The upper bound on $1-\C_{\sqrt{a(1-a)}}(A)$ is sharp for all set sizes, despite the odd looking $\max\{1/4,x\}$
term. When $\pi(A)\in[1/4,1/2]$ then
look at the walk on $K_n$ discussed in Example \ref{ex:complete}, with $\gamma=1/n$, that is $\alpha=0$.
Then $1-\C_{\sqrt{a(1-a)}}(A)=1$, $\~\psi^+(A)=1/2\pi(A^c)$ and $\~\psi^-(A)=1/2\pi(A)$ and so the bound becomes
$1-\C_{\sqrt{a(1-a)}}(A)\leq 1/2+\max\{1/8\pi(A),\,1/2\}$ which is correct when $\pi(A)\geq 1/4$.
When $\pi(A)<1/4$ then consider a Markov chain with transitions to be defined. Let $A$ be a set of size $x=\pi(A)$, 
let $B\subset A^c$ with $\pi(B)=\epsilon$, and let the transition kernel satisfy
$\forall v\in A:\,\P(v,B)=1-\P(v,A)=\epsilon/2x$ while $\forall v\in B:\,\P(v,A)=\P(v,B)=1/2$. 
Then it is easily computed that
$\~\psi^+(A)=\frac{\epsilon}{4\,z(1-z)^2}$, $\~\psi^-(A)=\frac{\epsilon^2}{4\,x^2(1-x)^2}$,
and $1-\C_{\sqrt{a(1-a)}}(A)=\frac 12 + \frac{\epsilon}{2x} - \frac 12\,\sqrt{(1+\epsilon/x)(1-\epsilon/(1-x))}$.
As $\epsilon\rightarrow 0^+$ then $(1-\C_{\sqrt{a(1-a)}}(A)-\~\psi^+(A)\,(1-x))/\~\psi^-(A)\rightarrow 1/4$
for all $z\in[0,1/2]$, and so $\~\psi^-(A)$ must always be multiplied by at least $1/4$ in the upper bound.

Another instance of sharpness is the lazy walk on the line of even length. It is easily checked that
$\~\psi^+(A),\,\~\psi^-(A) \geq 1/8n^2\,\pi(A)^2\pi(A^c)^2$, with equality when 
$A$ is an initial interval of the line. Then the lower bound on $1-\C_{\sqrt{a(1-a)}}$ is
achieved at $\pi(A)=1/2$, with $1-\C_{\sqrt{a(1-a)}} \geq 1/n^2$. The correct bound is
$1-\C_{\sqrt{a(1-a)}} = \frac 12\,\left(1-\sqrt{1-4/n^2}\right)\xrightarrow{n\rightarrow\infty} 1/n^2$,
and our bound was correct.

As mentioned earlier, the spread incorporates measures of both edge and vertex expansion, whereas (modified) conductance involves only edge expansion. Hence an improvement will be likely when vertex expansion is much larger than edge expansion. The most extreme example of this is a walk on the complete graph.

\begin{example}
Consider the lazy Markov chain on the complete graph $K_n$ given by choosing a vertex uniformly at random,
and moving there with probability $1/2$.

Conductance based bounds tend to be decent when considering the $L^2$ distance, but may be poor for other distances. In contrast, the $\chi^2$ bound is fine. Likewise, in this case the Cheeger inequality of Theorem \ref{thm:profile} is fine for the bound on $1-\C_{\sqrt{a(1-a)}}(A)$. In contrast, $1-\C_{a(1-a)}(A) = 1-\C_{a\log(1/a)}(A)= 1/2$ but the Cheeger inequalities of Theorem \ref{thm:profile} shows only 
$1-\C_{a(1-a)}(A)\geq \pi(A)\pi(A^c)$ and $(1-\pi(A))^2/2\log(1/\pi(A))$, both of which go to $0$ as $\pi(A)\rightarrow 0^+$.

We now use the spread. It is clear that $\Psi(A^c,t)=\frac{\pi(A^c)}{2}t$ if $t\leq\pi(A)$, while $\Psi(A,t)=\frac{\pi(A)}{2}t$ if $t\leq\pi(A^c)$, and so $\~\psi^+(A)=\frac{1}{4\pi(A^c)}$ while $\~\psi^-(A)=\frac{1}{4\pi(A)}$. The lower bounds are now within a factor of two, with $1-\C_{a(1-a)}(A) = 1/2$ and $1-\C_{a\log(1/a)}(A) \geq \frac 14$.
%
%
\end{example}


Compare this to the lazy random walk on the cycle of odd length (see Example \ref{ex:cycles}). 
In this case $\Q(A,A^c)=1/n$ and $\beta(A)=1/n$, so edge and vertex expansion are of similar orders. 
Therefore it is not surprising that the lower bounds on $\C_{a(1-a)}$ and $\C_{a\log(1/a)}$ given by
Theorem \ref{thm:profile} will be the correct order. 

Another case where vertex expansion is high is a walk on a product space. See \cite{Mon05.1} for a proof that the lazy walk on a Boolean cube $2^d$ has $\~\psi^+(A)=\Omega\left(\frac{1}{d\log d}\right)$, and so $1-\C_{\sqrt{a(1-a)}}(A)=\Omega\left(\frac{1}{d\log d}\right)$, a substantial improvement on the Cheeger inequality bound of $\Omega(1/d^2)$, although when it comes to mixing time this is still not as good as what can be shown by log-Sobolev or spectral methods.

%
%


\bibliographystyle{plain}
\cleardoublepage
\addcontentsline{toc}{chapter}{Bibliography}
\bibliography{../../references}

\end{document}